\documentclass{article}

\usepackage{graphicx}
\usepackage{float}
\usepackage{amsfonts}
\usepackage{amsthm}
\usepackage{amsmath}
\usepackage{amssymb}
\usepackage{mathrsfs}
\usepackage{tikz-cd}
\usepackage{stmaryrd}
\usepackage{enumitem}
\usepackage{placeins}
\usepackage{pdflscape}
\usepackage{pdflscape}
\usepackage[colorlinks=true,citecolor=blue]{hyperref}
\usepackage{adjustbox}

\usepackage{comment}

\theoremstyle{definition}
\newtheorem{Definition}{Definition}[subsection]

\theoremstyle{plain}
\newtheorem{Theorem}[Definition]{Theorem}

\theoremstyle{plain}
\newtheorem{Proposition}[Definition]{Proposition}

\theoremstyle{plain}
\newtheorem{Lemma}[Definition]{Lemma}

\theoremstyle{plain}
\newtheorem{AppLemma}{Lemma}[section]

\theoremstyle{definition}

\theoremstyle{plain}
\newtheorem{AppProposition}[AppLemma]{Proposition}

\theoremstyle{plain}
\newtheorem{AppTheorem}[AppLemma]{Theorem}

\theoremstyle{plain}

\theoremstyle{remark}

\theoremstyle{plain}
\newtheorem{AppQuestion}[AppLemma]{Question}

\theoremstyle{plain}
\newtheorem{Corollary}[Definition]{Corollary}

\theoremstyle{plain}
\newtheorem{Conjecture}[Definition]{Conjecture}

\theoremstyle{definition}
\newtheorem{Example}[Definition]{Example}

\theoremstyle{remark}
\newtheorem{Remark}[Definition]{Remark}

\theoremstyle{plain}
\newcommand{\thistheoremname}{}
\newtheorem*{genericthm*}{\thistheoremname}
\newenvironment{namedthm*}[1]
  {\renewcommand{\thistheoremname}{#1}%
   \begin{genericthm*}}
  {\end{genericthm*}}

\newlength{\diagramwidth}

\title{Higher Witt Groups for 2-Categories I: Centralizers}
\author{Hao Xu}
\date{March 2025}

\begin{document}

\maketitle

\begin{abstract}

In this article, we investigate monoidal, braided, sylleptic centralizers of monoidal, braided, sylleptic 2-functors. We specifically focus on multifusion 2-categories and show that monoidal, braided, sylleptic centralizers are multifusion again, via studying the corresponding enveloping algebras. We provide a characterization of the non-degeneracy condition for monoidal, braided, and sylleptic fusion 2-categories, via vanishing of their centers. Applying Double Centralizer Theorems, we establish the relationship between monoidal, braided, symmetric local modules and free modules. In particular, we obtain factorization properties of non-degenerate monoidal, braided, and sylleptic fusion 2-categories. Main results in this article will be used to study higher Witt equivalences of non-degenerate monoidal, braided, sylleptic 2-categories in the sequential articles.
\end{abstract}

{\hypersetup{linkcolor=black}\tableofcontents}

\section{Introduction}

\subsection*{Motivations}

Centers and centralizers are the most basic notions in the study of associative algebras and their representation theories. It is well known that centers of algebras are Morita invariant \cite{Mor58}, and by definition they are commutative. Centralizers generalize the notion of center to algebra homomorphisms, but they are in general not commutative and only subalgebras of the targets of homomorphisms.

\textit{Higher Algebra.} Center and centralizer are defined via universal properties in \cite{Lur17} for $\mathbb{E}_n$ algebras in symmetric monoidal $(\infty,1)$-categories. The first difference between classical algebra theory and higher algebra theory comes from Eckmann-Hilton arguments. In the world of Abelian groups, by iteration of algebra structures, we get a trichotomy of Abelian groups, associative algebras, commutative algebras. By contrast, in the world of chain complexes or spectra, there is an entire hierarchy of commutativity, denoted by $\mathbb{E}_n$ algebras for positive integer $n$, together with a limiting case of $\mathbb{E}_\infty$ algebras. Taking centers of $\mathbb{E}_n$ algebras only increases the level of commutativity by one, so there is a whole family of notations for $\mathbb{E}_n$ centers of $\mathbb{E}_n$ algebras. By carefully defining $\mathbb{E}_n$ algebras and their modules using $\infty$-operads, \cite{Lur17} proved that an $\mathbb{E}_n$ center of an $\mathbb{E}_n$ algebra is equipped with a canonical $\mathbb{E}_{n+1}$ algebra structure, while an $\mathbb{E}_n$ centralizer of an $\mathbb{E}_n$ algebra homomorphism is again an $\mathbb{E}_n$ algebra.

\textit{Stabilization Hypothesis.} The study of $\mathbb{E}_n$ structures also appeared in earlier works of algebraic topology \cite{May} and higher category theory \cite{BD}. For weak $n$-categories, one can equip them with multiple monoidal structures. The $1$-tuple monoidal $n$-categories are the categorification of associative algebras and monoidal categories. The $2$-tuple monoidal $n$-categories are the categorification of braided monoidal categories, whose commutativity is exactly between monoidal categories and symmetric monoidal categories, and play an important role in the study of representations of quantum groups \cite{Dri88}. The Baez-Dolan stabilization hypothesis states that for $k \geq n+2$, a $k$-tuple monoidal $n$-category is symmetric monoidal (i.e. $\mathbb{E}_\infty$-monoidal). The low-dimensional cases are summarized as follows.
\begin{center}
\begin{tabular}{ |c|c|c|c| } 
 \hline
 {} & {Abelian Group} & {1-Category} & {2-Category} \\ \hline
 {$\mathbb{E}_1$} & {Algebra} & {Monoidal Cat.} & {Monoidal 2-Cat.} \\ 
 {$\mathbb{E}_2$} & {Comm. Algebra} & {Braided Cat.} & {Braided 2-Cat.} \\ 
 {$\mathbb{E}_3$} & {-} & {Symmetric Cat.} & {Sylleptic. 2-Cat.} \\ 
 {$\mathbb{E}_4$} & {-} & {-} & {Symmetric 2-Cat.} \\ 
 \hline
\end{tabular}
\end{center}

Heuristically, one has the following list of centers and centralizers for 2-categories. Rigorous definitions of the $\mathbb{E}_n$ centers of 2-categories can be found in \cite{Cr}, where they are also proven to be $\mathbb{E}_{n+1}$ monoidal.

\begin{center}
\begin{tabular}{ |c|c|c| } 
 \hline
 {2-cat. $\mathfrak{M}$} & {endo-hom $\mathbf{End}(\mathfrak{M})$} & {$\mathbb{E}_1$ 2-cat.} \\ 
 {2-functor $\mathfrak{M} \xrightarrow{F} \mathfrak{N}$} & {hom $\mathbf{Fun}(\mathfrak{M},\mathfrak{N})$} & {2-cat.} \\ \hline
 {$\mathbb{E}_1$ 2-cat. $\mathfrak{C}$} & {$\mathbb{E}_1$ center $\mathscr{Z}_1(\mathfrak{C})$} & {$\mathbb{E}_2$ 2-Cat.} \\ 
 {$\mathfrak{C} \xrightarrow{F} \mathfrak{D}$} & {$\mathbb{E}_1$ centralizer $\mathscr{Z}_1(F)$} & {$\mathbb{E}_1$ 2-Cat.} \\ \hline
 {$\mathbb{E}_2$ 2-cat. $\mathfrak{B}$} & {$\mathbb{E}_2$ center $\mathscr{Z}_1(\mathfrak{B})$} & {$\mathbb{E}_3$ 2-Cat.} \\ 
 {$\mathfrak{A} \xrightarrow{F} \mathfrak{B}$} & {$\mathbb{E}_2$ centralizer $\mathscr{Z}_1(F)$} & {$\mathbb{E}_2$ 2-Cat.} \\ \hline
 {$\mathbb{E}_3$ 2-cat. $\mathfrak{S}$} & {$\mathbb{E}_3$ center $\mathscr{Z}_1(\mathfrak{S})$} & {$\mathbb{E}_\infty$ 2-Cat.} \\ 
 {$\mathfrak{S} \xrightarrow{F} \mathfrak{T}$} & {$\mathbb{E}_3$ centralizer $\mathscr{Z}_1(F)$} & {$\mathbb{E}_3$ 2-Cat.} \\ 
 \hline
\end{tabular}
\end{center}

\textit{Tensor Categories.} Tensor categories are categorification of algebras, and \cite{EGNO} is our standard reference. Drinfeld introduced the notion of $\mathbb{E}_1$ center for a tensor category, or known as Drinfeld center in the literature, in his unpublished note on the quantum double of Hopf algebra. The definition of Drinfeld center first appeared in \cite{JS91}, and also in \cite{Maj91} independently. In \cite{Ost2}, Ostrik showed that Drinfeld center is invariant under Morita equivalence. Morita equivalence of tensor categories, introduced by \cite{FRS} in the context of rational 2-dimensional conformal field theory and by \cite{Mu1} in the context of subfactors, is reformulated in \cite{EO03} using Morita dual category for a module category \cite{Ost1,Ost2}.

In contrast, the notion of $\mathbb{E}_1$ centralizer is less known in the literature. It is used, for example in \cite[Section 3.6]{DGNO}, on fusion subcategories in a braided fusion category, to prove the Double Centralizer Theorem of Müger centralizers (see also Appendix \ref{sec:DoubleCentralizerTheorem}). Another special case of $\mathbb{E}_1$ centralizer, called Morita dual tensor categories \cite{Mu1,Ost2,EO03}, is more renowned in the literature. More precisely, a module category $\mathcal{M}$ over a monoidal category is equivalent to a monoidal functor $H^\mathcal{M}:\mathcal{C} \to \mathbf{End}(\mathcal{M})$, and the Morita dual $\mathcal{C}^*_\mathcal{M} := \mathbf{End}_\mathcal{C}(\mathcal{M})$ is equivalent to the $\mathbb{E}_1$ centralizer $\mathcal{Z}_1(H^\mathcal{M})$. On the other hand, given any monoidal functor $F:\mathcal{C} \to \mathcal{D}$, one can realize the $\mathbb{E}_1$ centralizer as $\mathcal{Z}_1(F) \simeq \mathbf{End}_{\mathcal{C}|\mathcal{D}}({}_{\langle F \rangle} \mathcal{D})$, where $\mathcal{D}$ is viewed as $(\mathcal{C},\mathcal{D})$-bimodule with left $\mathcal{C}$-action induced by $F$ and right $\mathcal{D}$-action induced by right translation. Essentially, there is no information lost if we restrict from $\mathbb{E}_1$ centralizers to Morita duals.

\textit{Braided Tensor Categories.} The notion of $\mathbb{E}_2$ center of a braided tensor category is introduced by M{\"u}ger \cite{Mu0}. Later M{\"u}ger introduced $\mathbb{E}_2$ centralizer for a tensor subcategory of a braided tensor category in \cite{Mu3}. In general, given a braided tensor functor $F:\mathcal{A} \to \mathcal{B}$, it factors uniquely as a surjective one followed by an embedding: $\mathcal{A} \twoheadrightarrow \mathcal{C} \hookrightarrow \mathcal{B}$, where $\mathcal{C}$ is the tensor subcategory of $\mathcal{B}$ spanned by the image of $F$. Then the M{\"u}ger centralizer of $F$ is just the M{\"u}ger centralizer of $\mathcal{C}$ in $\mathcal{B}$.

The Morita theoretic interpretation of M{\"u}ger center and M{\"u}ger centralizer requires to study monoidal module categories over braided monoidal categories \cite{DGNO}. As a generalization of Witt equivalence for metric groups, Davydov et al. developed the theory of Witt equivalence of braided fusion categories in \cite{DMNO,DNO}. Brochier et al. further combined these two ingredients: in \cite{BJS}, they defined a Morita 4-category with braided tensor categories as objects, central monoidal bimodule categories as 1-morphisms, central bimodule categories as 2-morphisms; in \cite{BJSS}, they showed that non-degenerate braided fusion categories yield invertible objects in this Morita 4-category, and two non-degenerate braided fusion categories are equivalent as objects in this 4-category if they are Witt equivalent in the sense of \cite{DMNO}.

\textit{Fusion 2-Categories.} Fusion 2-categories are introduced as a categorification of fusion 1-categories by Douglas and Reutter \cite{DR}, where they used the notion further to construct 4-dimensional TQFTs via state sum construction. Décoppet developed the theory of fusion 2-categories in \cite{D1,D2,D3,D5}. A detailed theory of Morita equivalence for algebras within a fusion 2-category and Morita equivalence between fusion 2-categories has been accomplished in \cite{D4,D7,D8,D9,D10}.

Drinfeld center of a monoidal 2-category is introduced in \cite{BN}. Décoppet showed that the Drinfeld center of a fusion 2-category is a braided fusion 2-category, and Morita equivalent fusion 2-categories have braided equivalent Drinfeld centers \cite{D9}. Morita dual 2-category, introduced in \cite{D8}, categorified the 1-categorical notion. This is an example of $\mathbb{E}_1$ centralizer of a monoidal 2-functor. In Section 3, we prove that any $\mathbb{E}_1$ centralizer $\mathscr{Z}_1(F)$ of a monoidal 2-functor $\mathfrak{C} \xrightarrow{F} \mathfrak{D}
$ can be written as the Morita dual of $\mathfrak{C} \boxtimes \mathfrak{D}^{1mp}$ with respect to an associated module 2-category ${}_{\langle F \rangle} \mathfrak{D}$, hence the two theories are essentially the same.

\textit{Witt Equivalences.} $\mathbb{E}_1$ Morita invariance of Drinfeld centers of fusion 2-categories implies $\mathbb{E}_2$ morita invariance of Müger centers of braided fusion 1-categories. Given a braided fusion 1-category $\mathcal{B}$, its module 1-categories form a fusion 2-category $\mathbf{Mod}(\mathcal{B})$ under the relative Deligne tensor product $\boxtimes_\mathcal{B}$ \cite{ENO09,DSPS14}. Fusion 2-categories arise in such a way is called connected fusion 2-categories \cite{DR}. In \cite[Section 3]{D9}, Décoppet proved that two braided fusion 1-categories $\mathcal{A}$ and $\mathcal{B}$ are Witt equivalent (in the sense of \cite{DNO}, provided that their Müger centers are equivalent: $\mathcal{Z}_2(\mathcal{A}) \simeq \mathcal{Z}_2(\mathcal{B})$) if and only if connected fusion 2-categories $\mathbf{Mod}(\mathcal{A})$ and $\mathbf{Mod}(\mathcal{B})$ are Morita equivalent (in the sense of \cite{D8}).

\textit{Topological Field Theory and Cobordism Hypothesis.} Finally, let us mention some connections to mathematical physics. Topological Field Theories (abbreviated as TQFTs) are introduced by Atiyah \cite{Ati}, inspired by Segal's work on axiomization of two dimensional conformal field theories (abbreviated as CFTs) \cite{Seg}. Baez and Dolan \cite{BD} generalized the definition to symmetric $n$-categories by extending dimensions down to the point, and Lurie \cite{Lur08} enhanced the definition to symmetric $(\infty,n)$-categories by incorporate diffeomorphisms into the category of cobordisms. The well-known Cobordism Hypothesis (proved in \cite{Lur08,AF}) states that a TQFT is determined by the valued assigned to the point, which is fully dualizable in the target.

Using Cobordism Hypothesis, one can construct many TQFTs from higher Morita categories. For 2-dimensional cases, \cite{SP} provides a detailed analysis of framed, oriented and unoriented TQFTs taken values in the 2-category of algebras, bimodules and bimodule morphisms. In 3-dimension, \cite{DSPS13} constructs a Morita 3-category of finite tensor categories, exact bimodule categories and so on. They showed that fully dualizable objects in this 3-category are exactly multifusion categories. In 4-dimension, \cite{BJS} constructs a Morita 4-category of braided tensor categories, central monoidal bimodule categories, central bimodule categories and so on. They showed that braided multifusion categories are fully dualizable in this 4-category. More generally, \cite{D10} shows that multifusion 2-categories and their finite semisimple bimodule 2-categories form a Morita 4-category, where objects are all fully dualizable.

In future works, we would like to construct a Morita 5-category of braided multifusion 2-categories and a Morita 6-category of sylleptic multifusion 2-categories, and show that objects are fully dualizable. We will establish various Witt equivalence relations, and interpret them under the Morita theoretic contexts.

\textit{Topological Orders.} The theory of fusion 2-categories also have applications in the study of topological orders, a physical notion which was introduced in \cite{Wen90} and is one of the most active fields of research in condensed matter physics (see \cite{Wen19} for a recent review). It is hypothesized that low-energy effective theory for topological orders are TQFTs (which originated from high-energy physics). In (2+1)D, topological defects of codimension 2 and higher form a modular tensor category \cite{LW, Kit06,KK}, those of codimension 1 and higher form a fusion 2-category \cite{KK,DR}. More generally, all topological defects in an $(n+1)$D topological order form a fusion $n$-category \cite{KW, GJF,JF,KLWZZ,KZ20,KZ21}. The ongoing program of classifying all topological orders (or SPT/SET orders) demands us to develop the Morita theory of higher fusion categories, where higher centers and centralizers serve as building blocks.
%Topological quantum phases of matter, on the other hand, are studied in condensed matter physics. After the success of Landau’s paradigm to describe phase transitions via symmetry breaking, physicists found more general types of phases, possibly beyond Landau’s paradigm. Integer and fractional quantum Hall effects provide such archetypal examples. Short- and long-range entanglement between spined particles are speculated to be the explanation behind these exotic quantum phases of matters. Among all of them, the simplest instances are given by quantum phases whose energy spectrum has a gap above ground states, and whose local observables are trivial. The deformation classes of such phases are called topological orders \cite{Wen90} (see also \cite{Wen19} for a recent review).

% Despite there exist no interesting local quantum symmetries in topological orders, we still have global algebraic structures from fusion of topological defects of all codimensions. On the physical level, there has been classifications of topological orders and their equivariant companions, called symmetry protected/enriched topological (SPT/SET) orders \cite{LKW,JF,KLWZZ,KZ20,KZ21}. 

\subsection*{Main Results}

Let $\mathfrak{C},\mathfrak{D}$ be multifusion 2-categories, $\mathfrak{A}, \mathfrak{B}$ be braided multifusion 2-categories, $\mathfrak{S},\mathfrak{T}$ be sylleptic multifusion 2-categories. We have proven that $\mathbb{E}_1$, $\mathbb{E}_2$, $\mathbb{E}_3$ centralizers between multifusion 2-categories are multifusion.

\renewcommand{\thistheoremname}{Results A}

\begin{genericthm*}
(Proposition \ref{prop:MonoidalCentralizerIsMultifusion}) For any monoidal 2-functor $F:\mathfrak{C} \to \mathfrak{D}$, its Drinfeld centralizer $\mathscr{Z}_1(F)$ is multifusion. 

(Proposition \ref{prop:BraidedCentralizerIsBraidedMultifusion}) For any braided 2-functor $F:\mathfrak{A} \to \mathfrak{B}$, its braided centralizer $\mathscr{Z}_2(F)$ is braided multifusion. 

(Proposition \ref{prop:SyllepticCentralizerIsSyllepticMultifusion}) For any sylleptic 2-functor $F:\mathfrak{S} \to \mathfrak{T}$, its sylleptic centralizer $\mathscr{Z}_3(F)$ is sylleptic multifusion. 
\end{genericthm*}

We prove the above results by characterizing the centers as Morita duals of various enveloping algebras.

\renewcommand{\thistheoremname}{Results B}

\begin{genericthm*}
(Corollary \ref{cor:DrinfeldCenterAsBimoduleFunctors}) $\mathscr{Z}_1(F) \simeq \mathbf{End}_{\mathfrak{C}|\mathfrak{D}}({}_{\langle F \rangle} \mathfrak{D})$.

(Lemma \ref{lem:BraidedCentralizerAsMoritaDual}) $\mathscr{Z}_2(F) \simeq \mathbf{End}_{\int_{S^1}{F}}(\mathfrak{B})$.

(Remark \ref{rmk:S2EnvelopingAlgebra}) $\mathscr{Z}_3(F) \simeq \mathbf{End}_{\int_{S^2}{F}}(\mathfrak{S})$.
\end{genericthm*}

Let $A$ be a separable algebra (Definition \ref{def:algebra}, Definition \ref{def:SeparableAlgebra}) in a multifusion 2-category $\mathfrak{C}$, $B$ be a separable braided algebra (Definition \ref{def:braidedalgebra}) in a braided multifusion 2-category $\mathfrak{B}$, $S$ be a symmetric algebra (Definition \ref{def:symmetricalgebra}) in a sylleptic multifusion 2-category $\mathfrak{S}$. 

\renewcommand{\thistheoremname}{Results C}

\begin{genericthm*}
Induction 2-functors $Ind^\pm:\mathbf{Mod}_\mathfrak{B}(B) \to \mathbf{Bimod}_\mathfrak{B}(B)$ induce multifusion structure on $\mathbf{Mod}_\mathfrak{B}(B)$ (Corollary \ref{cor:MonoidalStructureOnModulesByInduction2Functor}). 

Similarly, induction 2-functors $Ind^\pm:\mathbf{Mod}_\mathfrak{B}(B) \to \mathbf{Mod}^{\mathbb{E}_1}_\mathfrak{B}(B)$ induce braided multifusion structure on $\mathbf{Mod}_\mathfrak{S}(S)$ (Corollary \ref{cor:BraidingOnModulesByInduction2Functor}).
\end{genericthm*}

Meanwhile, there are canonical 2-functors $\mathfrak{C} \to \mathbf{Mod}_\mathfrak{C}(A)$, $\mathfrak{B} \to \mathbf{Mod}_\mathfrak{B}(B)$, $\mathfrak{S} \to \mathbf{Mod}_\mathfrak{S}(S)$ sending objects to their corresponding free modules. These free modules 2-functors are endowed with central structures:

\renewcommand{\thistheoremname}{Results D}

\begin{genericthm*}
Lemma \ref{lem:SendingFreeModulesToEndohomOfModules} constructs a monoidal 2-functor \[\mathfrak{C}^{1mp} \to \mathbf{End}(\mathbf{Mod}_\mathfrak{C}(A)).\] 

Lemma \ref{lem:SendingFreeModulesToDrinfeldCenterOfModules} constructs a braided 2-functor \[\mathfrak{B}^{2mp} \to \mathscr{Z}_1(\mathbf{Mod}_\mathfrak{B}(B)).\]

Lemma \ref{lem:freeE3localmodule} constructs a sylleptic 2-functor \[\mathfrak{S}^{3mp} \to \mathscr{Z}_2(\mathbf{Mod}_\mathfrak{S}(S)).\]
\end{genericthm*}

For separable algebra $A$ in $\mathfrak{C}$, its $\mathbb{E}_1$ local module is defined to be an $(A,A)$-bimodule in $\mathfrak{C}$ (Definition \ref{def:bimodule}). For separable braided algebra $B$ in $\mathfrak{B}$, its $\mathbb{E}_2$ local module was introduced in \cite{ZLZHKT,Pom,DX} (Definition \ref{def:E2localmodule}). For separable symmetric algebra $S$ in $\mathfrak{S}$, we introduce $\mathbb{E}_3$ local module (Definition \ref{def:E3localmodule}). Then we characterize local modules as centralizers of free modules: 

\renewcommand{\thistheoremname}{Results E}

\begin{genericthm*}
Theorem \ref{thm:E1LocalModulesAsE1Centralizer} establishes a monoidal 2-equivalence \[\mathbf{Mod}^{\mathbb{E}_1}_\mathfrak{C}(A) \simeq \mathscr{Z}_1(\mathfrak{C}^{1mp} \to \mathbf{End}(\mathbf{Mod}_\mathfrak{C}(A))). \] 

Theorem \ref{thm:E2LocalModulesAsE2Centralizer} establishes a braided 2-equivalence \[\mathbf{Mod}^{\mathbb{E}_2}_\mathfrak{B}(B) \simeq \mathscr{Z}_2(\mathfrak{B}^{2mp} \to \mathscr{Z}_1(\mathbf{Mod}_\mathfrak{B}(B))). \] 

Theorem \ref{thm:E3LocalModulesAsE3Centralizer} establishes a sylleptic 2-equivalence \[\mathbf{Mod}^{\mathbb{E}_3}_\mathfrak{S}(S) \simeq \mathscr{Z}_3(\mathfrak{S}^{3mp} \to \mathscr{Z}_2(\mathbf{Mod}_\mathfrak{S}(S))).\]
\end{genericthm*}

Conversely, by Double Centralizer Theorems (Appendix \ref{sec:DoubleCentralizerTheorem}), one can realize free modules as the centralizers of local modules. 

\renewcommand{\thistheoremname}{Results F}

\begin{genericthm*}
In Theorem \ref{thm:DrinfeldCenterAndE1LocalModules}, for separable algebra $A$ in an indecomposable multifusion 2-category $\mathfrak{C}$, one has \[\mathscr{Z}_1(\mathfrak{C}) \simeq \mathscr{Z}_1(\mathfrak{C}^{1mp} \boxtimes \mathbf{Mod}^{\mathbb{E}_1}_\mathfrak{C}(A) \to \mathbf{End}(\mathbf{Mod}_\mathfrak{C}(A))).\]

In Theorem \ref{thm:E2CenterAndE2LocalModules}, for separable braided algebra $B$ in a braided fusion 2-category $\mathfrak{B}$, one has \[\mathscr{Z}_2(\mathfrak{B}) \simeq \mathscr{Z}_2(\mathfrak{B}^{2mp} \boxtimes \mathbf{Mod}^{\mathbb{E}_2}_\mathfrak{B}(B) \to \mathscr{Z}_1(\mathbf{Mod}_\mathfrak{B}(B))).\]

In Theorem \ref{thm:E3CenterAndE3LocalModules}, for separable symmetric algebra $S$ in a sylleptic fusion 2-category $\mathfrak{S}$, one has \[\mathscr{Z}_3(\mathfrak{S}) \simeq \mathscr{Z}_3(\mathfrak{S}^{3mp} \boxtimes \mathbf{Mod}^{\mathbb{E}_3}_\mathfrak{S}(S) \to \mathscr{Z}_2(\mathbf{Mod}_\mathfrak{S}(S))).\]
\end{genericthm*}

As corollaries, one has the factorization properties (which generalize the factorization property of modular tensor categories \cite[Theorem 7.10]{Mu2}).

\renewcommand{\thistheoremname}{Results G}

\begin{genericthm*}
(Corollary \ref{cor:FactorizationOfNondegenerateFusion2Cateory}) Any non-degenerate fusion 2-category $\mathfrak{C}$ (see Definition \ref{def:NonDegenerateFusion2Category}) has $\mathscr{Z}_1(\mathfrak{C}) \simeq \mathfrak{C}^{1mp} \boxtimes \mathfrak{C}.$ \vspace*{1mm}

(Corollary \ref{cor:FactorizationOfNondegenerateBraidedFusion2Cateory}) Any non-degenerate braided fusion 2-category $\mathfrak{B}$ (see Definition \ref{def:NonDegenerateBraidedFusion2Category}) has $\mathscr{Z}_2(\mathfrak{B}) \simeq \mathfrak{B}^{2mp} \boxtimes \mathfrak{B}.$ \vspace*{1mm}

(Corollary \ref{cor:FactorizationOfNondegenerateSyllepticFusion2Cateory}) Any non-degenerate sylleptic fusion 2-category $\mathfrak{S}$ (see Definition \ref{def:NonDegenerateSyllepticFusion2Category}) has $\mathscr{Z}_3(\mathfrak{S}) \simeq \mathfrak{S}^{3mp} \boxtimes \mathfrak{S}.$
\end{genericthm*}

Furthermore, we obtain the partial generalizations of \cite[Corollary 4.5]{Sch}.

\renewcommand{\thistheoremname}{Results H}

\begin{genericthm*}
(Corollary \ref{cor:E0CenterOfModuleCategoryIsE1LocalModuleInE0Center}) $\mathbf{Mod}^{\mathbb{E}_1}_{\mathbf{End}(\mathfrak{C})}(A) \simeq \mathbf{End}(\mathbf{Mod}_\mathfrak{C}(A)).$ \vspace*{1mm}

(Corollary \ref{cor:E1CenterOfModuleCategoryIsE2LocalModuleInE1Center}) $\mathbf{Mod}^{\mathbb{E}_2}_{\mathscr{Z}_1(\mathfrak{B})}(B) \simeq \mathscr{Z}_1(\mathbf{Mod}_\mathfrak{B}(B)).$ \vspace*{1mm}

(Corollary \ref{cor:E2CenterOfModuleCategoryIsE3LocalModuleInE2Center}) $\mathbf{Mod}^{\mathbb{E}_3}_{\mathscr{Z}_2(\mathfrak{S})}(S) \simeq \mathscr{Z}_2(\mathbf{Mod}_\mathfrak{S}(S)).$
\end{genericthm*}

\subsection*{Outline}

In Section 2, we recall basic notions for use in the later sections. We introduce the graphical calculus in semistrict monoidal 2-categories, and its extensions to braided, sylleptic and symmetric 2-categories, and define ($\mathbb{E}_1$-, $\mathbb{E}_2$-, $\mathbb{E}_3$-) monoidal 2-functors, ($\mathbb{E}_1$-, $\mathbb{E}_2$-) monoidal 2-natural transforms, and monoidal modifications between ($\mathbb{E}_1$-, $\mathbb{E}_2$-, $\mathbb{E}_3$-, $\mathbb{E}_\infty$-) monoidal 2-categories. Then we define ($\mathbb{E}_1$-, $\mathbb{E}_2$-, $\mathbb{E}_3$-) algebras, modules, and ($\mathbb{E}_1$-, $\mathbb{E}_2$-, $\mathbb{E}_3$-) local modules in monoidal, braided and sylleptic 2-categories. For separable ($\mathbb{E}_1$-, $\mathbb{E}_2$-, $\mathbb{E}_3$-) algebras in an ($\mathbb{E}_1$-, $\mathbb{E}_2$-, $\mathbb{E}_3$-) multifusion 2-category, we cite many results from \cite{D4,D7,D8,DY,DX} on 2-categories of modules and local modules. We generalize the induction 2-functors \cite{DY} from $\mathbb{E}_1$ local modules to $\mathbb{E}_2$ local modules. Finally, we introduce the notion of module 2-categories and enriched 2-categories \cite{D4}.

In Section 3, we first define $\mathbb{E}_1$ center of a monoidal 2-category \cite{BN} and $\mathbb{E}_1$ centralizer of a monoidal 2-functor. To establish that $\mathbb{E}_1$ centralizers of multifusion 2-categories are again multifusion, we characterize $\mathbb{E}_1$ centralizers as Morita duals. Moving on, we show that taking $\mathbb{E}_1$ centers is additive and multiplicative. Then we define non-degeneracy of multifusion 2-categories via vanishing of their $\mathbb{E}_1$ centers. Lastly, we demonstrate that free modules and $\mathbb{E}_1$ local modules centralize each other in the endo-hom of the 2-category of modules and prove the factorization property of non-degenerate fusion 2-categories.

In Section 4, we first define $\mathbb{E}_2$ center of a braided 2-category \cite{Cr} and $\mathbb{E}_2$ centralizer of a braided 2-functor. To establish that $\mathbb{E}_2$ centralizers of braided multifusion 2-categories are again multifusion, we characterize $\mathbb{E}_2$ centralizers as Morita duals to $S^1$-enveloping algebras, which are multifusion 2-categories themselves. Moving on, we show that taking $\mathbb{E}_2$ centers is additive and multiplicative. Then we define non-degeneracy of braided fusion 2-categories via vanishing of their $\mathbb{E}_2$ centers. Lastly, we demonstrate that free modules and $\mathbb{E}_2$ local modules centralize each other in the $\mathbb{E}_1$ center of the 2-category of modules and prove the factorization property of non-degenerate braided fusion 2-categories.

In Section 5, we first define $\mathbb{E}_3$ center of a sylleptic 2-category \cite{Cr} and $\mathbb{E}_3$ centralizer of a sylleptic 2-functor. Then we prove that $\mathbb{E}_3$ centralizers of sylleptic multifusion 2-categories are again multifusion, since they are monoidal sub-2-categories of the target sylleptic multifusion 2-categories. Moving on, we show that taking $\mathbb{E}_3$ centers is additive and multiplicative. Next, we define non-degeneracy of sylleptic fusion 2-categories via vanishing of their $\mathbb{E}_2$ centers. However, we observe certain pathology of this notion, compared to the $\mathbb{E}_1$ and $\mathbb{E}_2$ analogies. So we denote the above non-degeneracy of sylleptic fusion 2-categories as weak non-degeneracy, and compare it with a stronger version of non-degeneracy via higher condensation theory \cite{GJF,JF,KZ20,KZ21}. Lastly, we demonstrate that free modules and $\mathbb{E}_3$ local modules centralize each other in the $\mathbb{E}_2$ center of the 2-category of modules and prove the factorization property of non-degenerate sylleptic fusion 2-categories.

Appendix \ref{sec:proofE3LocalModulesAsE3Centralizer} contains figures for the proof of Theorem \ref{thm:E3LocalModulesAsE3Centralizer}. Appendix \ref{sec:DoubleCentralizerTheorem} is a brief survey of double $\mathbb{E}_n$ centralizers in higher fusion categories, established upon \cite{GJF,JF,KZ20,KZ21}.

\subsection*{Acknowledgements}

I would like to express my sincere gratitude to Thibault Décoppet for his invaluable assistance and insightful discussions during the initial stages of this project. I am particularly grateful for his outstanding series of work on the foundations of algebras and modules in 2-categories \cite{D4,D7,D8,D9,DY,DX,D10}, which contributed to the majority of the diagrams presented in Section 2 and 3. I would like to thank Zhi-Hao Zhang for many fruitful discussions, in particular about centers of sylleptic fusion 2-categories. The author was supported by DAAD Graduate School Scholarship Programme (57572629) and DFG Project 398436923.

I would like to thank Liang Kong, Zhi-Hao Zhang and Hao Zheng for feedback on the first version of this paper. After submission of the preprint to arXiv, another preprint \cite{KZZZ} appeared at the same time, which contains a physical discussion on a similar topic. I would like to thank the authors of \cite{KZZZ} for pointing out the connection to their work and many valuable comments on the interpretation of our results. 

A mistake was realized in the Appendix B of the first version, which is mainly related to the proof of double centralizing property of $\mathbb{E}_2$ centralizers. That section has been taken out in this revised version, and now all the proofs involving double centralizing properties follow from the higher condensation theory (in Appendix B, which was the old Appendix C). On the other hand, the claimed reciprocity may still hold in a weaker context (see also \cite[Remark 4.2.17]{KZZZ}), and we will investigate this in future works.

\section{Preliminaries}

\subsection{Graphical Calculus}

We work within a monoidal 2-category $\mathfrak{C}$ with monoidal product $\Box$ and monoidal unit $I$ in the sense of \cite{SP}. Thanks to the coherence theorem of \cite{Gur}, we may assume without loss of generality that $\mathfrak{C}$ is strict cubical in the sense of Definition 2.26 of \cite{SP}. We use the graphical calculus of \cite{GS}, as described in \cite{D4} and \cite{D7}. We often omit the symbol $\Box$ from our notations, and we use the symbol $1_f$ to denote an identity $(n+1)$-morphism for an $n$-morphism $f$. We may abbreviate the subscripts depending on the context. Following \cite{GJF}, we use $\Omega \mathfrak{C}$ to denote the endomorphism category of the monoidal unit of monoidal 2-category $\mathfrak{C}$.

The interchanger is depicted using the string diagram below on the left, and its inverse by that on the right: $$\begin{tabular}{c c c c}
\includegraphics[width=20mm]{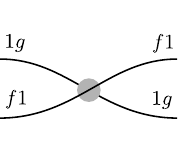},\ \ \ \   & \ \ \ \  \includegraphics[width=20mm]{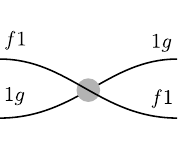}.
\end{tabular}$$ In particular, the lines correspond to 1-morphisms, and the coupons to 2-morphisms. The regions represent objects, which are uniquely determined by the 1-morphisms. Further, our string diagrams are read from top to bottom, which yields the compositions of 1-morphisms, and then from left to right.

For our purposes, it is also necessary to recall the graphical conventions related to 2-natural transformations from \cite{GS}. These will only be used for the braiding, which will be introduced below. Let $F,G:\mathfrak{A}\rightarrow \mathfrak{B}$ be two (weak) 2-functors, and let $\tau:F\Rightarrow G$ be 2-natural transformation. That is, for every object $A$ in $\mathfrak{A}$, we have a 1-morphism $\tau_A:F(A)\rightarrow G(A)$, and for every 1-morphism $f:A\rightarrow B$ in $\mathfrak{A}$, we have a 2-isomorphism $$\begin{tikzcd}[sep=tiny]
F(A) \arrow[ddd, "F(f)"']\arrow[rrr, "\tau_A"]  &                                        &    & G(A) \arrow[ddd, "G(f)"]  \\
 &  &    & \\
  &  &  &  \\
F(B)\arrow[rrr, "\tau_B"']\arrow[Rightarrow, rrruuu, "\tau_f", shorten > = 2ex, shorten < = 2ex]                                            &                                        &    &  G(B), 
\end{tikzcd}$$
These 2-isomorphisms have to satisfy obvious coherence relations. In our graphical language, we will depict the 2-isomorphism $\tau_f$ using the following diagram on the left, and its inverse using the diagram on the right: $$\begin{tabular}{c c c c}
\includegraphics[width=20mm]{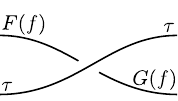},\ \ \ \   & \ \ \ \  \includegraphics[width=20mm]{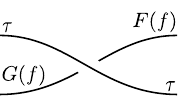}.
\end{tabular}$$ 

In the present article, we will for the most part work within $\mathfrak{B}$ a braided monoidal 2-category in the sense of \cite{SP}. Thanks to the coherence theorem of \cite{Gur2}, we may assume that $\mathfrak{B}$ is a semi-strict braided monoidal 2-category. In particular, $\mathfrak{B}$ comes equipped with a braiding $b$, which is an adjoint 2-natural equivalence given on objects $A,B$ in $\mathfrak{B}$ by $$b_{A,B}:A \, \Box \, B\rightarrow B \, \Box \, A.$$ Its pseudo-inverse will be denoted by $b^{\bullet}$. Further, there are two invertible modifications $R$ and $S$, which are given on the objects $A,B,C$ of $\mathfrak{B}$ by
\begin{center}
\begin{tabular}{@{}c c@{}}

$\begin{tikzcd}
ABC \arrow[rr, "b"] \arrow[rd, "b1"'] & {} \arrow[d, Rightarrow, "R"]          & BCA, \\
                                      & BAC \arrow[ru, "1b"'] &    
\end{tikzcd}$

&

$\begin{tikzcd}
ABC \arrow[rr, "b_2"] \arrow[rd, "1b"'] & {} \arrow[d, Rightarrow, "S"]     & CAB \\
                                        & ACB \arrow[ru, "b1"'] &    
\end{tikzcd}$
\end{tabular}
\end{center}

\noindent where the subscript in $b_2$ records were the braiding occurs. To avoid any possible confusion, we will systematically write $b$ instead of a would be $b_1$ as this can too easily be confused with $b1$. Further, these modifications are subject to the following relations, which are taken from section 2.1.1 of \cite{DY}:

\begin{enumerate}
\item [a.] For every objects $A,B,C,D$ in $\mathfrak{B}$, we have
\end{enumerate}

\newlength{\calculus}

\settoheight{\calculus}{\includegraphics[width=37.5mm]{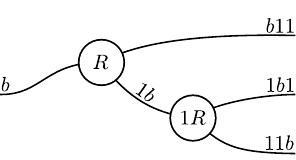}}

\begin{equation}\label{eqn:braidingaxiom1}
\begin{tabular}{@{}ccc@{}}

\includegraphics[width=37.5mm]{Pictures/Preliminaries/Braided2Category/braidingaxiom1.pdf} & \raisebox{0.45\calculus}{$=$} &
\includegraphics[width=37.5mm]{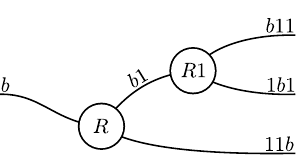}

\end{tabular}
\end{equation}

\begin{enumerate}
\item [] in $Hom_{\mathfrak{B}}(ABCD, BCDA)$,
\item [b.] For every objects $A,B,C,D$ in $\mathfrak{B}$, we have
\end{enumerate}

\settoheight{\calculus}{\includegraphics[width=37.5mm]{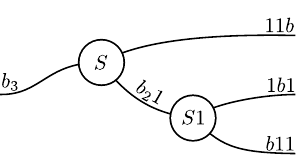}}

\begin{equation}\label{eqn:braidingaxiom2}
\begin{tabular}{@{}ccc@{}}

\includegraphics[width=37.5mm]{Pictures/Preliminaries/Braided2Category/braidingaxiom3.pdf} & \raisebox{0.45\calculus}{$=$} &
\includegraphics[width=37.5mm]{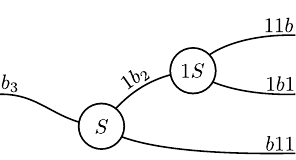}

\end{tabular}
\end{equation}

\begin{enumerate}
\item [] in $Hom_{\mathfrak{B}}(ABCD, DABC)$,
\item [c.] For every objects $A,B,C,D$ in $\mathfrak{B}$, we have
\end{enumerate}

\settoheight{\calculus}{\includegraphics[width=37.5mm]{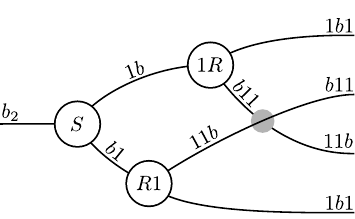}}

\begin{equation}\label{eqn:braidingaxiom3}
\begin{tabular}{@{}ccc@{}}

\includegraphics[width=45mm]{Pictures/Preliminaries/Braided2Category/braidingaxiom5.pdf} & \raisebox{0.45\calculus}{$=$} &
\includegraphics[width=37.5mm]{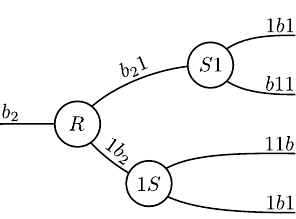}

\end{tabular}
\end{equation}

\begin{enumerate}
\item [] in $Hom_{\mathfrak{B}}(ABCD, CDAB)$,
\item [d.] For every objects $A,B,C$ in $\mathfrak{B}$, we have
\end{enumerate}

\settoheight{\calculus}{\includegraphics[width=45mm]{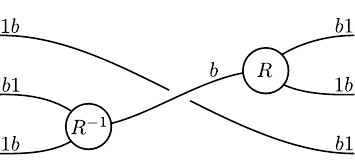}}

\begin{equation}\label{eqn:braidingaxiom4}
\begin{tabular}{@{}ccc@{}}

\includegraphics[width=45mm]{Pictures/Preliminaries/Braided2Category/braidingaxiom7.pdf} & \raisebox{0.45\calculus}{$=$} &
\includegraphics[width=45mm]{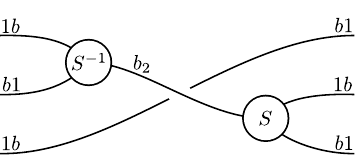}

\end{tabular}
\end{equation}

\begin{enumerate}
\item [] in $Hom_{\mathfrak{B}}(ABC, CBA)$,
\item [e.] For every object $A$ in $\mathfrak{B}$, the adjoint 2-natural equivalences $$b_{A,I}:A \, \Box \, I\rightarrow I \, \Box \, A\textrm{ and } b_{I,A}:I \, \Box \, A \rightarrow A \, \Box \, I$$ are the identity adjoint 2-natural equivalences,

\item [f.] For every objects $A,B,C$ in $\mathfrak{B}$, the 2-isomorphisms $R_{A,B,C}$ and $S_{A,B,C}$ are the identity 2-isomorphism whenever either $A$, $B$, or $C$ is equal to $I$.
\end{enumerate}

We would also like to recall sylleptic monoidal 2-categories from \cite{SP}. Suppose $\mathfrak{S}$ is a semi-strict braided monoidal 2-category with braiding $(b,R,S)$ given as above. A syllepsis on $\mathfrak{S}$ is an invertible modification $\sigma$

$$\begin{tikzcd}
AB \arrow[rr, equal] \arrow[rd, "b"'] & {} \arrow[d, Rightarrow, "\sigma"]          & AB, \\
                                      & BA \arrow[ru, "b"'] &    
\end{tikzcd}$$ satisfying:

\begin{enumerate}
\item [a.] For every objects $A,B,C$ in $\mathfrak{S}$, we have
\end{enumerate}

\settoheight{\calculus}{\includegraphics[width=35mm]{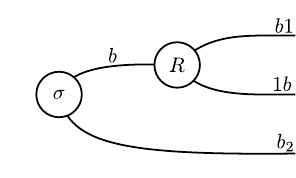}}

\begin{equation}\label{eqn:syllepsisaxiom1}
\begin{tabular}{@{}ccc@{}}

\includegraphics[width=35mm]{Pictures/Preliminaries/Sylleptic2Category/syllepsis1left.pdf} & \raisebox{0.45\calculus}{$=$} &
\includegraphics[width=49mm]{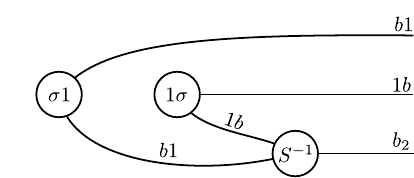}

\end{tabular}
\end{equation}

\begin{enumerate}
\item [] in $Hom_{\mathfrak{S}}(ABC, ABC)$,
\item [b.] For every objects $A,B,C$ in $\mathfrak{S}$, we have
\end{enumerate}

\settoheight{\calculus}{\includegraphics[width=35mm]{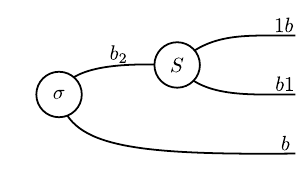}}

\begin{equation}\label{eqn:syllepsisaxiom2}
\begin{tabular}{@{}ccc@{}}

\includegraphics[width=35mm]{Pictures/Preliminaries/Sylleptic2Category/syllepsis2left.pdf} & \raisebox{0.45\calculus}{$=$} &
\includegraphics[width=49mm]{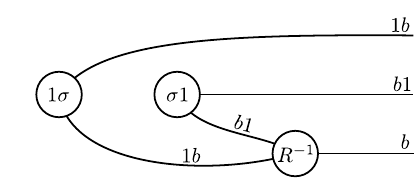}

\end{tabular}
\end{equation}

\begin{enumerate}
\item [] in $Hom_{\mathfrak{S}}(ABC, ABC)$,
\item [c.] For every objects $A,B$ in $\mathfrak{S}$, we have $\sigma_{A,B}$ equals the identity 2-morphism whenever $A$ or $B$ equals $I$.
\end{enumerate}

A braided monoidal 2-category $(\mathfrak{S},b,R,S)$ together with a syllepsis $\sigma$ is called a sylleptic monoidal 2-category.

Finally, recall from \cite{SP} that a sylleptic monoidal 2-category $\mathfrak{S}$ is symmetric if its syllepsis satisfies the condition 
\begin{equation} \label{eqn:symmetricsyllepsis}
    \sigma_{B,A} \circ b_{A,B} = b_{B,A} \circ \sigma_{A,B},
\end{equation} for any objects $A,B$ in $\mathfrak{S}$.

We will also work under the linear setting, following \cite{DR,GJF}. For simplicity, let us assume that the ground field $\Bbbk$ is algebraically closed of characteristic zero. In detail, a multifusion 2-category is a rigid\footnote{Rigidity means every object in the monoidal 2-category has a left and right dual, and every 1-morphism admits a left and right adjoint.} monoidal 2-category whose underlying 2-category is finite semisimple, and the monoidal structure is compatible with the linear structure. If moreover the monoidal unit is simple as an object, then we say it is a fusion 2-category. See \cite{D5,D2} for a generalization beyond the assumption of $\Bbbk$ being algebraically closed of characteristic zero. We use $\pi_0 \mathfrak{C}$ to denote the set of connected components for simples objects in finite semisimple 2-category $\mathfrak{C}$.

Finally, recall that a left dual for object $x$ in monoidal 2-category $\mathfrak{C}$ consists of an object ${}^\lor x$ together with unit 1-morphism $i_x:I \to x \, \Box \, {}^\lor x$, counit 1-morphism $e_x:{}^\lor x \, \Box \, x \to I$ and 2-isomorphisms witnessing zigzag conditions:
    $$\Xi_x:(1_x \, \Box \, e_x) \circ (i_x \, \Box \, 1_x) \to 1_x, $$
    $$ \Phi_x:(e_x \, \Box \, 1_{{}^\lor x}) \circ (1_{{}^\lor x} \, \Box \, i_x) \to 1_{{}^\lor x}.$$
    Furthermore, we can assume this left dual is coherent \cite{Pstr}, i.e. $\Xi_x$ and $\Phi_x$ satisfy the two coherence conditions:

\settoheight{\calculus}{\includegraphics[width=30mm]{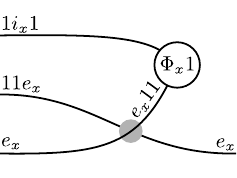}}

\begin{equation}\label{eqn:coherentleftdual1}
\begin{tabular}{@{}ccc@{}}

\includegraphics[width=30mm]{Pictures/Preliminaries/CoherentDual/CoherentLeftDual1left.pdf} & \raisebox{0.45\calculus}{$=$} &
\includegraphics[width=30mm]{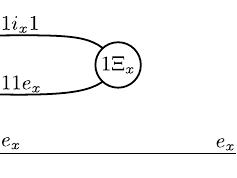},

\end{tabular}
\end{equation}

\settoheight{\calculus}{\includegraphics[width=30mm]{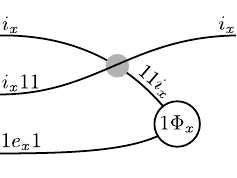}}

\begin{equation}\label{eqn:coherentleftdual2}
\begin{tabular}{@{}ccc@{}}

\includegraphics[width=30mm]{Pictures/Preliminaries/CoherentDual/CoherentLeftDual2left.pdf} & \raisebox{0.45\calculus}{$=$} &
\includegraphics[width=30mm]{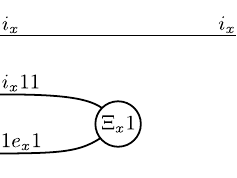}.

\end{tabular}
\end{equation}
    One can similarly define a coherent right dual of object $x$ by simply reversing the direction of monoidal product $\Box$.

In \cite[Appendix A]{DX}, we see that taking left or right dual objects in a rigid monoidal 2-category $\mathfrak{C}$ can be promoted into a monoidal functor $\mathfrak{C} \to \mathfrak{C}^{1mp,1op}$.

We use $\mathfrak{M}^{1op}$ to denote the 2-category obtained by reversing the direction of 1-arrows in the 2-category $\mathfrak{M}$. We use $\mathfrak{C}^{1mp}$ to denote the monoidal 2-category obtained by reversing the direction of monoidal product in the monoidal 2-category $\mathfrak{C}$. Similarly, we use $\mathfrak{B}^{2mp}$ to denote the braided monoidal 2-category obtained by reversing the direction of braiding in the braided monoidal 2-category $\mathfrak{B}$, and we use $\mathfrak{S}^{3mp}$ to denote the sylleptic monoidal 2-category obtained by reversing the direction of syllepsis in the sylleptic monoidal 2-category $\mathfrak{S}$. More specifically, suppose $\mathfrak{S}$ has sylleptic monoidal structure $(\Box,b,R,S,\sigma)$, then $\mathfrak{S}^{3mp}$ has the same underlying braided monoidal structure $(\Box,b,R,S,\sigma)$ and its syllepsis $\overline{\sigma}$ is determined on objects $x$ and $y$ as follows:
\settoheight{\calculus}{\includegraphics[width=21mm]{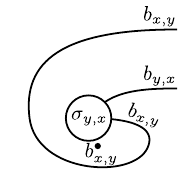}}

\begin{equation}\label{eqn:reversedsyllepsis}
\begin{tabular}{@{}ccccc@{}}

\raisebox{0.45\calculus}{$\overline{\sigma}_{x,y}:=$} & \includegraphics[width=21mm]{Pictures/Preliminaries/Sylleptic2Category/ReversedSyllepsis1.pdf} & \raisebox{0.45\calculus}{$=$} &
\includegraphics[width=21mm]{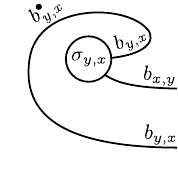} \raisebox{0.45\calculus}{.} &

\end{tabular}
\end{equation}

\subsection{2-Functors}

In this section, we shall recall various notions of 2-functors, 2-natural transforms and modifications between monoidal, braided, sylleptic and symmetric 2-categories. Readers may find more detailed discussions and references in \cite[Chapter 2.3]{SP}.

Suppose $\mathfrak{C}$ and $\mathfrak{D}$ are two semi-strict monoidal 2-categories.
\begin{Definition} \label{def:Monoidal2Functor}
    A monoidal 2-functor from $\mathfrak{C}$ to $\mathfrak{D}$ consists of:
    \begin{enumerate}
        \item[1.] An underlying 2-functor $F:\mathfrak{C} \to \mathfrak{D}$;

        \item[2.] A 2-natural equivalence $\chi_{x,y}:F(x) \, \Box^\mathfrak{D} \, F(y) \to F(x \, \Box^\mathfrak{C} \, y)$ given on objects $x,y$ in $\mathfrak{C}$;

        \item[3.] An isomorphism $\iota:I^\mathfrak{D} \to F(I^\mathfrak{C})$;

        \item[4.] Invertible modifications
    \end{enumerate}
    \[\begin{tikzcd}[column sep=45pt,row sep=25pt]
            {F(x) \, \Box^\mathfrak{D} \, F(y) \, \Box^\mathfrak{D} \, F(z)}
                \arrow[r,"\chi_{x,y} \, \Box^\mathfrak{D} \, F(z)"]
                \arrow[d,"F(x) \, \Box^\mathfrak{D} \, \chi_{y,z}"']
            & {F(x \, \Box^\mathfrak{C} \, y) \, \Box^\mathfrak{D} \, F(z)}
                \arrow[d,"\chi_{xy,z}"]
            \\ {F(x) \, \Box^\mathfrak{D} \, F(y \, \Box^\mathfrak{C} \, z)}
                \arrow[r,"\chi_{x,yz}"']
                \arrow[ur,Rightarrow,shorten <=20pt, shorten >=20pt,"\omega_{x,y,z}"]
            & {F(x \, \Box^\mathfrak{C} \, y \, \Box^\mathfrak{C} \, z)}
        \end{tikzcd}, \]
        \[\begin{tikzcd}[column sep=40pt]
            {F(I^\mathfrak{C}) \, \Box^\mathfrak{D} \, F(x)}
                \arrow[rd,"\iota \, \Box^\mathfrak{D} \, F(x)"]
                \arrow[dd,equal]
            & {}
            & {F(x) \, \Box^\mathfrak{D} \, F(I^\mathfrak{C})}
                \arrow[ld,"F(x) \, \Box^\mathfrak{D} \, \iota"']
                \arrow[dd,equal]
            \\ {}
            & {F(x)}
                \arrow[dl,"\chi_{I,x}"]
                \arrow[dr,"\chi_{x,I}"']
                \arrow[l,Rightarrow,"\gamma_x"',shorten <= 30pt, shorten >= 30pt]
            & {}
                \arrow[l,Rightarrow,"\delta_x"',shorten <= 30pt, shorten >= 30pt]
            \\ {F(I^\mathfrak{C}) \, \Box^\mathfrak{D} \, F(x)}
            & {}
            & {F(x) \, \Box^\mathfrak{D} \, F(I^\mathfrak{C})}
        \end{tikzcd},\]
        \begin{enumerate}
            \item[] given on objects $x,y,z$ in $\mathfrak{C}$;
        \end{enumerate}
        subject to the following conditions:
        \begin{enumerate}
            \item[a.] For any objects $x,y,z,w$ in $\mathfrak{C}$, the equation holds
        \end{enumerate}

        \settoheight{\calculus}{\includegraphics[width=28mm]{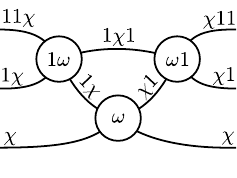}}

        \begin{equation}\label{eqn:Monoidal2FunctorAssociativity}
        \begin{tabular}{@{}ccc@{}}
    
        \includegraphics[width=28mm]{Pictures/Preliminaries/Monoidal2Functor/pentagon1left.pdf} & \raisebox{0.45\calculus}{$=$} &
        \includegraphics[width=28mm]{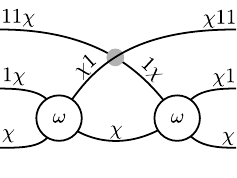} 
        
        \end{tabular}
        \end{equation}

        \begin{enumerate}
            \item[] in $Hom_\mathfrak{D}(F(x) \, \Box^\mathfrak{D} \, F(y) \, \Box^\mathfrak{D} \, F(z) \, \Box^\mathfrak{D} \, F(w),G(x \, \Box^\mathfrak{C} \, y \, \Box^\mathfrak{C} \, z \, \Box^\mathfrak{C} \, w ))$;
        \end{enumerate}

        \begin{enumerate}
            \item[b.] For any objects $x,y$ in $\mathfrak{C}$, the equation holds
        \end{enumerate}

        \settoheight{\calculus}{\includegraphics[width=28mm]{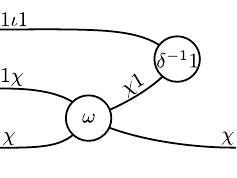}}

        \begin{equation}\label{eqn:Monoidal2FunctorUnitality}
        \begin{tabular}{@{}ccc@{}}
    
        \includegraphics[width=28mm]{Pictures/Preliminaries/Monoidal2Functor/triangle1left.pdf} & \raisebox{0.45\calculus}{$=$} &
        \includegraphics[width=14mm]{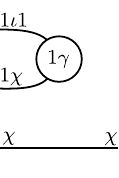} 
        
        \end{tabular}
        \end{equation}

        \begin{enumerate}
            \item[] in $Hom_\mathfrak{D}(F(x) \, \Box^\mathfrak{D} \, F(y),G(x \, \Box^\mathfrak{C} \, y))$.
        \end{enumerate}
\end{Definition}

\begin{Definition} \label{def:Monoidal2NatTrans}
    Let $F$ and $G$ be two monoidal 2-functors given from $\mathfrak{C}$ to $\mathfrak{D}$. A monoidal 2-natural transformation from $F$ to $G$ consists of:
    \begin{enumerate}
        \item[1.] An underlying 2-natural transformation $\eta:F \to G$;

        \item[2.] Invertible modification
        \[\begin{tikzcd}[column sep=40pt,row sep=30pt]
            {F(x) \, \Box^\mathfrak{D} \, F(y)}
                \arrow[r,"\eta_x \, \Box^\mathfrak{D} \, F(y)"]
                \arrow[d,"\chi^F_{x,y}"']
            & {G(x) \, \Box^\mathfrak{D} \, F(y)}
                \arrow[r,"G(x) \, \Box^\mathfrak{D} \eta_y"]
            & {G(x) \, \Box^\mathfrak{D} \, G(y)}
                \arrow[d,"\chi^G_{x,y}"]
                \arrow[lld,Rightarrow,"\Pi_{x,y}"',shorten <=70pt, shorten >=50pt]
            \\ {F(x \, \Box^\mathfrak{C} \, y)}
                \arrow[rr,"\eta_{xy}"']
            & {}
            & {G(x \, \Box^\mathfrak{C} \, y)}
        \end{tikzcd},\] and a 2-isomorphism
        \[\begin{tikzcd}[row sep=30pt]
            {I^\mathfrak{D}}
                \arrow[rr,"\iota^G"]
                \arrow[rd,"\iota^F"']
            & {}
            & {G(I^\mathfrak{C})}
            \\ {}
            & {F(I^\mathfrak{C})}
                \arrow[ur,"\eta_I"']
                \arrow[u,Rightarrow,shorten <=10pt, shorten >=5pt,"M",]
            & {}
        \end{tikzcd};\]
    \end{enumerate}
    subject to the following conditions:
    \begin{enumerate}
        \item[a.] For any objects $x,y,z$ in $\mathfrak{C}$, the equation holds
    \end{enumerate}

    \settoheight{\calculus}{\includegraphics[width=35mm]{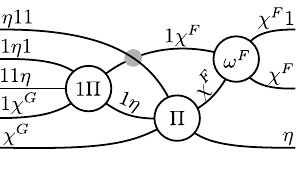}}

    \begin{equation}\label{eqn:Monoidal2NatTransAssociativity}
    \begin{tabular}{@{}ccc@{}}
    
    \includegraphics[width=35mm]{Pictures/Preliminaries/Monoidal2Functor/hexagonleft.pdf} & \raisebox{0.45\calculus}{$=$} &
    \includegraphics[width=35mm]{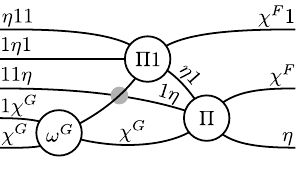} 
        
    \end{tabular}
    \end{equation}

    \begin{enumerate}
        \item[] in $Hom_\mathfrak{D}(F(x) \, \Box^\mathfrak{D} \, F(y) \, \Box^\mathfrak{D} \, F(z),G(x \, \Box^\mathfrak{C} \, y \, \Box^\mathfrak{C} \, z))$;
    \end{enumerate}
        
    \begin{enumerate}
        \item[b.] For any object $x$ in $\mathfrak{C}$, the equation holds
    \end{enumerate}

    \settoheight{\calculus}{\includegraphics[width=28mm]{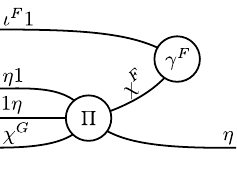}}

    \begin{equation}\label{eqn:Monoidal2NatTransUnitality1}
    \begin{tabular}{@{}ccc@{}}
    
    \includegraphics[width=28mm]{Pictures/Preliminaries/Monoidal2Functor/square1left.pdf} & \raisebox{0.45\calculus}{$=$} &
    \includegraphics[width=28mm]{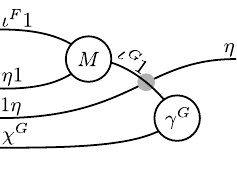} 
        
    \end{tabular}
    \end{equation}

    \begin{enumerate}
        \item[] in $Hom_\mathfrak{D}(F(x),G(x))$;
    \end{enumerate}
    
    \begin{enumerate}
        \item[c.] For any object $x$ in $\mathfrak{C}$, the equation holds
    \end{enumerate}

    \settoheight{\calculus}{\includegraphics[width=35mm]{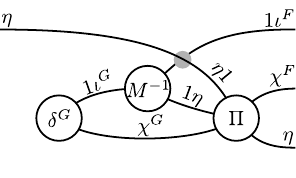}}

    \begin{equation}\label{eqn:Monoidal2NatTransUnitality2}
    \begin{tabular}{@{}ccc@{}}
    
    \includegraphics[width=35mm]{Pictures/Preliminaries/Monoidal2Functor/square2left.pdf} & \raisebox{0.45\calculus}{$=$} &
    \includegraphics[width=21mm]{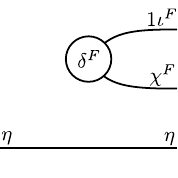} 
        
    \end{tabular}
    \end{equation}

    \begin{enumerate}
        \item[] in $Hom_\mathfrak{D}(F(x),G(x))$.
    \end{enumerate}
\end{Definition}

\begin{Definition} \label{def:MonoidalModification}
    Let $\eta$ and $\phi$ be two monoidal 2-natural transforms from $F$ to $G$, which are monoidal 2-functors given from $\mathfrak{C}$ to $\mathfrak{D}$. A monoidal modification from $\eta$ to $\phi$ consists of an underlying modification $\theta:\eta \to \phi$ subject to the two equations:
    \begin{enumerate}
        \item[a.] For any objects $x,y$ in $\mathfrak{C}$, the equation holds
    \end{enumerate}

    \settoheight{\calculus}{\includegraphics[width=28mm]{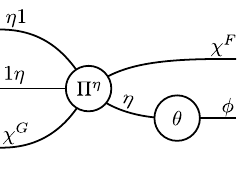}}

    \begin{equation}\label{eqn:MonoidalModificationAssociativity}
    \begin{tabular}{@{}ccc@{}}
    
    \includegraphics[width=28mm]{Pictures/Preliminaries/Monoidal2Functor/pentagon2left.pdf} & \raisebox{0.45\calculus}{$=$} &
    \includegraphics[width=28mm]{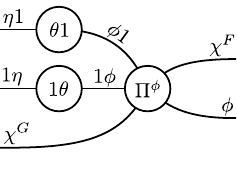} 
        
    \end{tabular}
    \end{equation}

    \begin{enumerate}
        \item[] in $Hom_\mathfrak{D}(F(x) \, \Box^\mathfrak{D} \, F(y),G(x \, \Box^\mathfrak{C} \, y))$;
    \end{enumerate}

    \begin{enumerate}
        \item[b.] The following equation holds:
    \end{enumerate}

    \settoheight{\calculus}{\includegraphics[width=28mm]{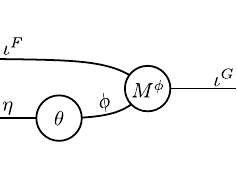}}

    \begin{equation}\label{eqn:MonoidalModificationUnitality}
    \begin{tabular}{@{}cccc@{}}
    
    \includegraphics[width=28mm]{Pictures/Preliminaries/Monoidal2Functor/triangle2left.pdf} & \raisebox{0.45\calculus}{$=$} &
    \includegraphics[width=21mm]{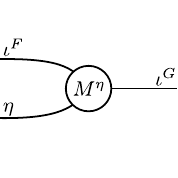} & \raisebox{0.45\calculus}{.}
        
    \end{tabular}
    \end{equation}
\end{Definition}

Suppose $\mathfrak{A}$ and $\mathfrak{B}$ are two braided 2-categories. 

\begin{Definition} \label{def:Braided2Functor}
    A braided 2-functor from $\mathfrak{A}$ to $\mathfrak{B}$ consists of:
    \begin{enumerate}
        \item An underlying 2-functor $F:\mathfrak{A} \to \mathfrak{B}$;

        \item A monoidal 2-functor structure $(\chi,\iota,\omega,\gamma,\delta)$ on $F$;

        \item An invertible modification
        \[\begin{tikzcd}[row sep=30pt, column sep=45pt]
            {F(x) \, \Box^\mathfrak{B} \, F(y)}
                \arrow[d,"\chi_{x,y}"']
                \arrow[r,"b^\mathfrak{B}_{F(x),F(y)}"]
            & {F(y) \, \Box^\mathfrak{B} \, F(x)}
                \arrow[d,"\chi_{y,x}"]
            \\ {F(x \, \Box^\mathfrak{A} \, y)}
                \arrow[r,"F(b^\mathfrak{A}_{x,y})"']
                \arrow[ur,Rightarrow,shorten <=20pt, shorten >=25pt,"u_{x,y}"]
            & {F(y \, \Box^\mathfrak{A} \, x)}
        \end{tikzcd}\]
    \end{enumerate}

    \begin{enumerate}
        \item[] given on objects $x,y$ in $\mathfrak{A}$;
    \end{enumerate}
    subject to the extra conditions:
    \begin{enumerate}
        \item[a.] For any objects $x,y,z$ in $\mathfrak{A}$, the equation holds
    \end{enumerate}

    \settoheight{\calculus}{\includegraphics[width=35mm]{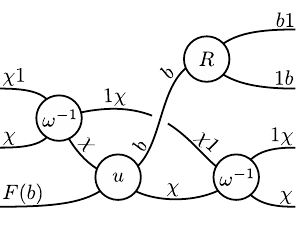}}

    \begin{equation}\label{eqn:Braided2Functor1}
    \begin{tabular}{@{}ccc@{}}
    
    \includegraphics[width=35mm]{Pictures/Preliminaries/Monoidal2Functor/braided2functor1left.pdf} & \raisebox{0.45\calculus}{$=$} &
    \includegraphics[width=49mm]{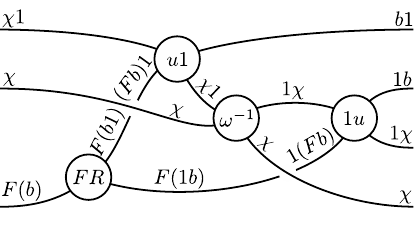}
        
    \end{tabular}
    \end{equation}
    
    \begin{enumerate}
        \item[] in $Hom_\mathfrak{B}(F(x) \, \Box^\mathfrak{B} \, F(y) \, \Box^\mathfrak{B} \, F(z),F(y \, \Box^\mathfrak{A} \, z \, \Box^\mathfrak{A} \, x))$;
    \end{enumerate}

    \begin{enumerate}
        \item[b.] For any objects $x,y,z$ in $\mathfrak{A}$, the equation holds
    \end{enumerate}

    \settoheight{\calculus}{\includegraphics[width=35mm]{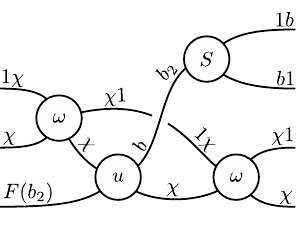}}

    \begin{equation}\label{eqn:Braided2Functor2}
    \begin{tabular}{@{}ccc@{}}
    
    \includegraphics[width=35mm]{Pictures/Preliminaries/Monoidal2Functor/braided2functor2left.pdf} & \raisebox{0.45\calculus}{$=$} &
    \includegraphics[width=49mm]{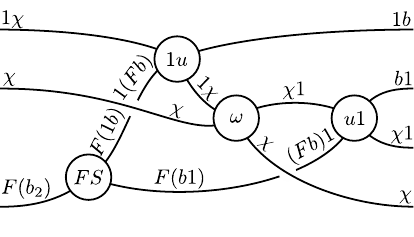}
        
    \end{tabular}
    \end{equation}
    
    \begin{enumerate}
        \item[] in $Hom_\mathfrak{B}(F(x) \, \Box^\mathfrak{B} \, F(y) \, \Box^\mathfrak{B} \, F(z),F(z \, \Box^\mathfrak{A} \, x \, \Box^\mathfrak{A} \, y))$.
    \end{enumerate}
\end{Definition}

\begin{Definition} \label{def:Braided2NatTrans}
    Let $F$ and $G$ be two braided 2-functors given from $\mathfrak{A}$ to $\mathfrak{B}$. A braided 2-natural transformation from $F$ to $G$ consists of:
    \begin{enumerate}
        \item[1.] An underlying 2-natural transformation $\eta:F \to G$;

        \item[2.] A monoidal 2-natural transform structure $(\Pi,M)$ on $\eta$;
    \end{enumerate}
    subject to the extra condition:
    \begin{enumerate}
        \item[a.] For any objects $x,y$ in $\mathfrak{A}$, the equation holds
    \end{enumerate}

    \settoheight{\calculus}{\includegraphics[width=35mm]{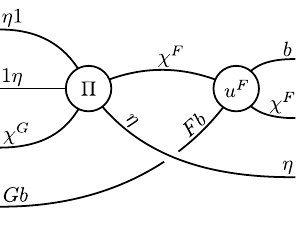}}

    \begin{equation}\label{eqn:Braided2NatTrans}
    \begin{tabular}{@{}ccc@{}}
    
    \includegraphics[width=35mm]{Pictures/Preliminaries/Monoidal2Functor/braided2nattransleft.pdf} & \raisebox{0.45\calculus}{$=$} &
    \includegraphics[width=35mm]{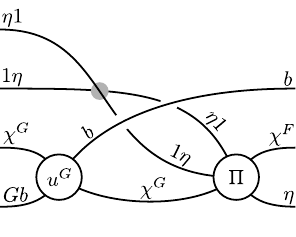}
        
    \end{tabular}
    \end{equation}
    
    \begin{enumerate}
        \item[] in $Hom_\mathfrak{B}(F(x) \, \Box^\mathfrak{B} \, F(y),G(x \, \Box^\mathfrak{A} \, y))$.
    \end{enumerate}
\end{Definition}

Suppose $\mathfrak{S}$ and $\mathfrak{T}$ are two sylleptic 2-categories. 

\begin{Definition} \label{def:Sylleptic2Functor}
    A sylleptic 2-functor from $\mathfrak{S}$ to $\mathfrak{T}$ consists of:
    \begin{enumerate}
        \item An underlying 2-functor $F:\mathfrak{S} \to \mathfrak{T}$;

        \item A braided 2-functor structure $(\chi,\iota,\omega,\gamma,\delta,u)$ on $F$;
    \end{enumerate}
    subject to an extra condition:
    \begin{enumerate}
        \item[a.] For any objects $x,y$ in $\mathfrak{S}$, the equation holds
    \end{enumerate}

    \settoheight{\calculus}{\includegraphics[width=35mm]{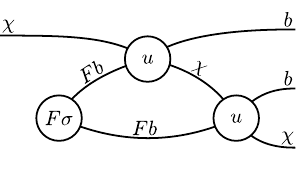}}

    \begin{equation}\label{eqn:Sylleptic2Functor}
    \begin{tabular}{@{}ccc@{}}
    
    \includegraphics[width=35mm]{Pictures/Preliminaries/Monoidal2Functor/sylleptic2functorleft.pdf} & \raisebox{0.45\calculus}{$=$} &
    \includegraphics[width=21mm]{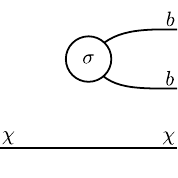}
        
    \end{tabular}
    \end{equation}
    
    \begin{enumerate}
        \item[] in $Hom_\mathfrak{T}(F(x) \, \Box^\mathfrak{T} \, F(y),F(x \, \Box^\mathfrak{S} \, y))$.
    \end{enumerate}
\end{Definition}

\begin{Lemma}
    One has the following 3-categories:
    \begin{enumerate}
        \item A 3-category consists of 
        \begin{itemize}
            \item monoidal 2-categories as objects,

            \item monoidal 2-functors as 1-morphisms,

            \item monoidal 2-natural transforms as 2-morphisms,

            \item monoidal modifications as 3-morphisms;
        \end{itemize}

        \item A 3-category consists of 
        \begin{itemize}
            \item braided 2-categories as objects,

            \item braided 2-functors as 1-morphisms,

            \item braided 2-natural transforms as 2-morphisms,

            \item monoidal modifications as 3-morphisms;
        \end{itemize}
        
        \item A 3-category consists of 
        \begin{itemize}
            \item sylleptic 2-categories as objects,

            \item sylleptic 2-functors as 1-morphisms,

            \item braided 2-natural transforms as 2-morphisms,

            \item monoidal modifications as 3-morphisms;
        \end{itemize}
        
        \item A 3-category consists of 
        \begin{itemize}
            \item symmetric 2-categories as objects,

            \item sylleptic 2-functors as 1-morphisms,

            \item braided 2-natural transforms as 2-morphisms,

            \item monoidal modifications as 3-morphisms.
        \end{itemize}
    \end{enumerate}
\end{Lemma}

\subsection{Algebras}

Let $\mathfrak{C}$ be a strict cubical monoidal 2-category. We recall the definition of an algebra in $\mathfrak{C}$ from \cite{D7}. The definition of an algebra in an arbitrary monoidal 2-category using our graphical conventions may be found in \cite{D4}.

\begin{Definition}\label{def:algebra}
An algebra in $\mathfrak{C}$ consists of:
\begin{enumerate}
    \item An object $A$ of $\mathfrak{C}$;
    \item Two 1-morphisms $m:A \, \Box \, A \rightarrow A$ and $i:I \rightarrow A$;
    \item Three 2-isomorphisms
\end{enumerate}
\begin{center}
\begin{tabular}{@{}c c c@{}}
$\begin{tikzcd}[sep=small]
A \arrow[rrrr, equal] \arrow[rrdd, "i1"'] &  & {} \arrow[dd, Rightarrow, "\lambda"', near start, shorten > = 1ex] &  & A \\
                                   &  &                           &  &   \\
                                   &  & AA, \arrow[rruu, "m"']     &  &  
\end{tikzcd}$

&

$\begin{tikzcd}[sep=small]
AAA \arrow[dd, "1m"'] \arrow[rr, "m1"]    &  & AA \arrow[dd, "m"] \\
                                            &  &                      \\
AA \arrow[rr, "m"'] \arrow[rruu, Rightarrow, "\mu", shorten > = 2.5ex, shorten < = 2.5ex] &  & A,                   
\end{tikzcd}$

&

$\begin{tikzcd}[sep=small]
                                  &  & AA \arrow[rrdd, "m"] \arrow[dd, Rightarrow, "\rho", shorten > = 1ex, shorten < = 2ex] &  &   \\
                                  &  &                                             &  &   \\
A \arrow[rruu, "1i"] \arrow[rrrr,equal] &  & {}                                          &  & A,
\end{tikzcd}$

\end{tabular}
\end{center}

satisfying:

\begin{enumerate}
\item [a.] We have:
\end{enumerate}

\newlength{\prelim}

\settoheight{\prelim}{\includegraphics[width=52.5mm]{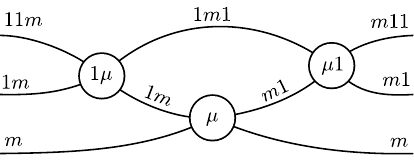}}

\begin{equation}\label{eqn:algebraassociativity}
\begin{tabular}{@{}ccc@{}}

\includegraphics[width=52.5mm]{Pictures/Preliminaries/Algebra/associativity1.pdf} & \raisebox{0.45\prelim}{$=$} &
\includegraphics[width=40mm]{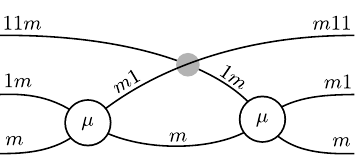},

\end{tabular}
\end{equation}

\begin{enumerate}
\item [b.] We have:
\end{enumerate}

\settoheight{\prelim}{\includegraphics[width=22.5mm]{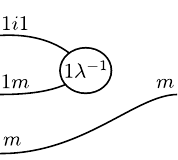}}

\begin{equation}\label{eqn:algebraunitality}
\begin{tabular}{@{}ccc@{}}

\includegraphics[width=22.5mm]{Pictures/Preliminaries/Algebra/unitality1.pdf} & \raisebox{0.45\prelim}{$=$} &

\includegraphics[width=37.5mm]{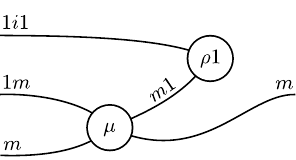}.

\end{tabular}
\end{equation}
\end{Definition}

\begin{Definition}\label{def:algebrahomomorphism}
Let $A$ and $B$ be two algebras in $\mathfrak{C}$. An algebra 1-morphism $f:A\rightarrow B$ consists of a 1-morphism $f:A\rightarrow B$ in $\mathfrak{C}$, together with two invertible 2-morphisms
 
\begin{center}
    \begin{tabular}{@{}c c@{}}
    $\begin{tikzcd}[sep=small]
    A A \arrow[rr, "f1"] \arrow[dd, "m^A"']  &  & B A \arrow[rr,"1f"]  & & B B \arrow[dd,"m^B"] \\
                                       &  &                        &  &   \\
                                       A \arrow[rrrr, "f"'] 
                                       \arrow[rrrruu, Rightarrow, "\psi^f", shorten > = 4ex, shorten < = 4ex]  & &  &  &  B,
    \end{tikzcd}$
    
    &
    
    $\begin{tikzcd}[sep=small]
                                      &  & A \arrow[rrdd, "f"] \arrow[dd, Rightarrow, "\eta^f", shorten > = 1ex, shorten < = 2ex] &  &   \\
                                      &  &                                             &  &   \\
    I \arrow[rruu, "i^A"] \arrow[rrrr,"i^B"'] &  & {}                                          &  & B,
    \end{tikzcd}$
\end{tabular}
\end{center}
satisfying:

\begin{enumerate}
\item [a.] We have:
\end{enumerate}

\settoheight{\prelim}{\includegraphics[width=47.5mm]{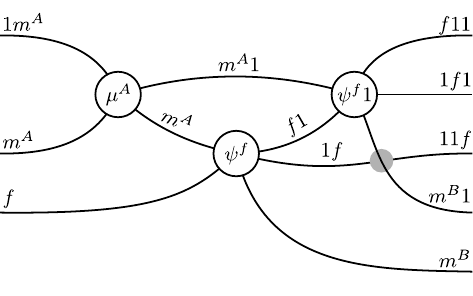}}

\begin{center}
\begin{equation} \label{eqn:algebra1morphismcoh1}
\begin{tabular}{@{}ccc@{}}

\includegraphics[width=47.5mm]{Pictures/Preliminaries/Algebra/algebra1morcoh1left.pdf} & \raisebox{0.45\prelim}{$=$} &

\includegraphics[width=47.5mm]{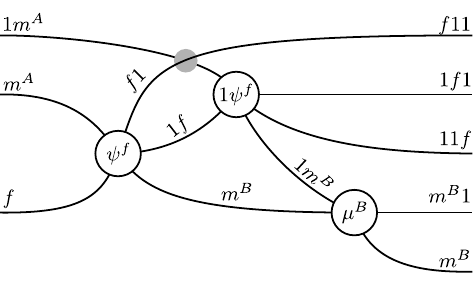},

\end{tabular}
\end{equation}
\end{center}

\begin{enumerate}
\item [b.] We have:
\end{enumerate}

\settoheight{\prelim}{\includegraphics[width=60mm]{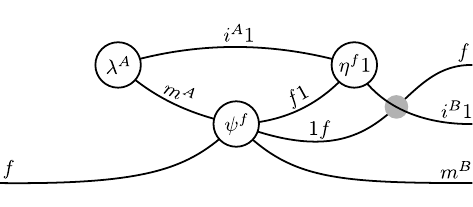}}

\begin{center}
\begin{equation} \label{eqn:algebra1morphismcoh2}
\begin{tabular}{@{}ccc@{}}

\includegraphics[width=60mm]{Pictures/Preliminaries/Algebra/algebra1morcoh2left.pdf} & \raisebox{0.45\prelim}{$=$} &

\includegraphics[width=30mm]{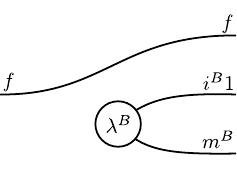},

\end{tabular}
\end{equation}
\end{center}

\begin{enumerate}
\item [c.] We have:
\end{enumerate}

\settoheight{\prelim}{\includegraphics[width=60mm]{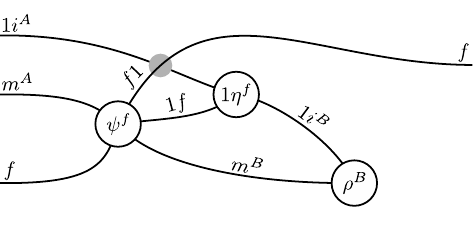}}

\begin{center}
\begin{equation} \label{eqn:algebra1morphismcoh3}
\begin{tabular}{@{}ccc@{}}

\includegraphics[width=60mm]{Pictures/Preliminaries/Algebra/algebra1morcoh3left.pdf} & \raisebox{0.45\prelim}{$=$} &

\includegraphics[width=30mm]{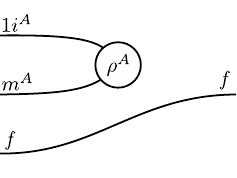}.

\end{tabular}
\end{equation}
\end{center}
\end{Definition}

\begin{Definition} \label{def:algebra2morphism}
Let $A$ and $B$ be two algebras, $f$ and $g$ be two algebra 1-morphisms from $A$ to $B$ in $\mathfrak{C}$. An algebra 2-morphism $\gamma$ consists of a 2-morphism $\gamma:f \to g$ satisfying:
\begin{enumerate}
\item [a.] We have:
\end{enumerate}

\settoheight{\prelim}{\includegraphics[width=45mm]{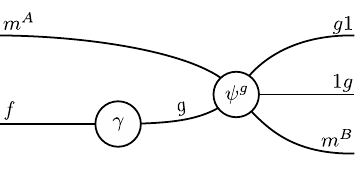}}

\begin{center}
\begin{equation} \label{eqn:algebra2morphismcoh1}
\begin{tabular}{@{}ccc@{}}

\includegraphics[width=45mm]{Pictures/Preliminaries/Algebra/algebra2morcoh1left.pdf} & \raisebox{0.45\prelim}{$=$} &

\includegraphics[width=45mm]{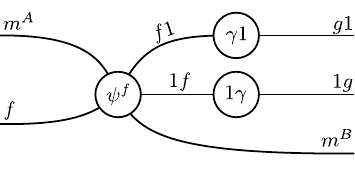},

\end{tabular}
\end{equation}
\end{center}

\begin{enumerate}
\item [b.] We have:
\end{enumerate}

\settoheight{\prelim}{\includegraphics[width=45mm]{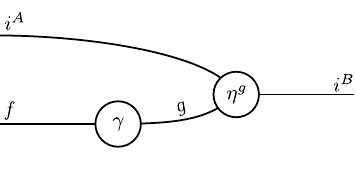}}

\begin{center}
\begin{equation} \label{eqn:algebra2morphismcoh2}
\begin{tabular}{@{}ccc@{}}

\includegraphics[width=45mm]{Pictures/Preliminaries/Algebra/algebra2morcoh2left.pdf} & \raisebox{0.45\prelim}{$=$} &

\includegraphics[width=30mm]{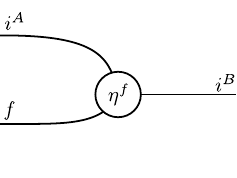}.

\end{tabular}
\end{equation}
\end{center}
\end{Definition}

\begin{Lemma}
    Algebras, algebra 1-morphisms and algebra 2-morphisms in the monoidal 2-category $\mathfrak{C}$ together form a 2-category.
\end{Lemma}

Now, let $\mathfrak{B}$ be a semi-strict braided monoidal 2-category. We recall the definition of braided algebra from \cite{DY}.

\begin{Definition}\label{def:braidedalgebra}
A braided algebra in $\mathfrak{B}$ consists of:
\begin{enumerate}
    \item An algebra $B$ in $\mathfrak{B}$;
    \item A 2-isomorphisms
\end{enumerate}

$$\begin{tikzcd}[sep=small]
                                  &  & BB \arrow[rrdd, "m"] \arrow[dd, Rightarrow, "\beta", shorten > = 1ex, shorten < = 2ex] &  &   \\
                                  &  &                                             &  &   \\
BB \arrow[rruu, "b"] \arrow[rrrr, "m"'] &  & {}                                          &  & B,
\end{tikzcd}$$

satisfying:

\begin{enumerate}
\item [a.] We have:
\end{enumerate}

\newlength{\braid}

\settoheight{\braid}{\includegraphics[width=52.5mm]{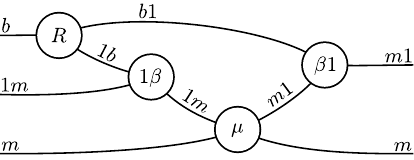}}

\begin{equation}\label{eqn:braidedalgebra1}
\begin{tabular}{@{}ccc@{}}

\includegraphics[width=45mm]{Pictures/Preliminaries/Algebra/braidedalgebra1.pdf} & \raisebox{0.45\braid}{$=$} &
\includegraphics[width=45mm]{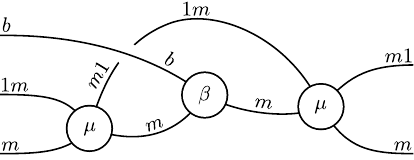},

\end{tabular}
\end{equation}

\begin{enumerate}
\item [b.] We have:
\end{enumerate}

\settoheight{\braid}{\includegraphics[width=52.5mm]{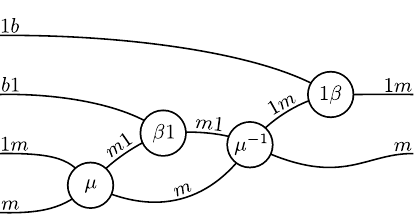}}

\begin{equation}\label{eqn:braidedalgebra2}
\begin{tabular}{@{}ccc@{}}

\includegraphics[width=45mm]{Pictures/Preliminaries/Algebra/braidedalgebra3.pdf} & \raisebox{0.45\braid}{$=$} &

\includegraphics[width=45mm]{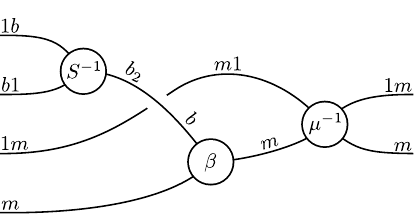},

\end{tabular}
\end{equation}

\begin{enumerate}
\item [c.] We have:
\end{enumerate}

\settoheight{\braid}{\includegraphics[width=30mm]{Pictures/Preliminaries/Algebra/braidedalgebra4.pdf}}

\begin{equation}\label{eqn:braidedalgebra3}
\begin{tabular}{@{}ccc@{}}

\includegraphics[width=30mm]{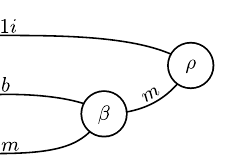} & \raisebox{0.45\braid}{$=$} &

\includegraphics[width=30mm]{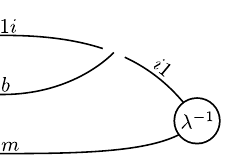}.

\end{tabular}
\end{equation}
\end{Definition}

\begin{Definition}\label{def:braidedalgebrahomomorphism}
Let $A$ and $B$ be two braided algebras in $\mathfrak{B}$. A braided algebra 1-homomorphism $f:A\rightarrow B$ is an algebra 1-homomorphism $f:A\rightarrow B$ that satisfies:

\settoheight{\braid}{\includegraphics[width=45mm]{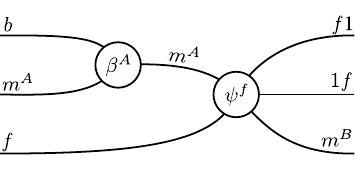}}

\begin{equation} \label{eqn:braidedalgebra1morphism}
\begin{tabular}{@{}ccc@{}}

\includegraphics[width=42mm]{Pictures/Preliminaries/Algebra/braidedalgebra1morcohleft.pdf} & \raisebox{0.45\braid}{$=$} &

\includegraphics[width=56mm]{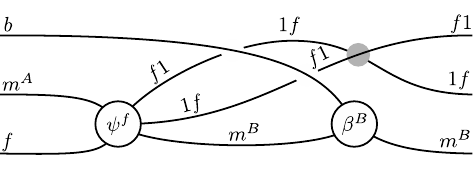}.

\end{tabular}
\end{equation}
\end{Definition}

\begin{Definition} \label{def:braidedalgebra2morphism}
    Let $f,g$ be two braided algebra 1-morphisms between braided algebras $A$ and $B$ in $\mathfrak{B}$. Then a braided algebra 2-morphism between them is just an ordinary algebra 1-morphism $\gamma:f \to g$, see Definition \ref{def:algebra2morphism}.
\end{Definition}

\begin{Lemma}
    Braided algebras, braided algebra 1-morphisms and 2-morphisms in the braided monoidal 2-category $\mathfrak{B}$ form a 2-category.
\end{Lemma}

Let $\mathfrak{S}$ be a semi-strict sylleptic monoidal 2-category. We recall the definition of symmetric algebra from \cite{DY}.

\begin{Definition} \label{def:symmetricalgebra}
    A symmetric algebra in $\mathfrak{S}$ consists of a braided algebra $(S,m,i,\lambda,\mu,\rho,\beta)$ satisfying the additional coherence condition:
    \begin{center}
    \settoheight{\diagramwidth}{\includegraphics[width=45mm]{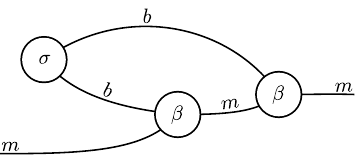}}  
    \begin{equation} \label{eqn:symmetricalgebra}
    \begin{tabular}{@{}ccc@{}}
        \includegraphics[width=45mm]{Pictures/Preliminaries/Algebra/symmetricalgebra.pdf}&\raisebox{0.45\diagramwidth}{$=$} &
        \raisebox{0.45\diagramwidth}{$1_m$.} 
        \end{tabular}
    \end{equation}
    \end{center}
    1-morphisms and 2-morphisms between symmetric algebras are just braided algebra 1-morphisms and 2-morphisms.
\end{Definition}

\begin{Lemma}
    Symmetric algebras, with 1-morphisms and 2-morphisms in the sylleptic monoidal 2-category $\mathfrak{S}$ form a 2-category.
\end{Lemma}

The definition of a rigid algebra is due to \cite{Gai}. Its relevance in the context of fusion 2-categories was first observed in \cite{JFR}. The definition of a separable algebra has its origin in the \cite{GJF}. These notions were studied extensively in \cite{D7}, to which we refer the reader for further discussion as well as examples.

\begin{Definition} \label{def:RigidAlgebra}
A rigid algebra in $\mathfrak{C}$ is an algebra $A$ such that the 1-morphism $m:A \, \Box \, A \rightarrow A$ admits a right adjoint $m^*$ as an $A$-$A$-bimodule 1-morphism.
\end{Definition}

\begin{Definition} \label{def:SeparableAlgebra}
A separable algebra in $\mathfrak{C}$ is a rigid algebra $A$ in $\mathfrak{C}$ such that the counit $\epsilon^m:m\circ m^* \to 1_A$ splits as an $A$-$A$-bimodule 2-morphism.
\end{Definition}

\begin{Definition} \label{def:ConnectedAlgebra}
A separable algebra $A$ in $\mathfrak{C}$ is connected if its unit 1-morphism $i:I \to A$ is simple in the finite semisimple 1-category $Hom_\mathfrak{C}(I,A)$.
\end{Definition}

\subsection{Modules, Local Modules}

Let us now recall the definition of a right $A$-module in $\mathfrak{C}$ from Definition 1.2.3 of \cite{D7}. We invite the reader to consult Definition 3.2.1 of \cite{D4} for the definition in a general monoidal 2-category.

\begin{Definition}\label{def:module}
A right $A$-module in $\mathfrak{C}$ consists of:
\begin{enumerate}
    \item An object $M$ of $\mathfrak{C}$;
    \item A 1-morphism $n^M:M \, \Box \, A \rightarrow M$;
    \item Two 2-isomorphisms
\end{enumerate}
\begin{center}
\begin{tabular}{@{}c c@{}}
$\begin{tikzcd}[sep=small]
MAA \arrow[dd, "1m"'] \arrow[rr, "n^M1"]    &  & MA \arrow[dd, "n^M"] \\
                                            &  &                      \\
MA \arrow[rr, "n^M"'] \arrow[rruu, Rightarrow, "\nu^M", shorten > = 2.5ex, shorten < = 2.5ex] &  & M,                   
\end{tikzcd}$

&

$\begin{tikzcd}[sep=small]
                                  &  & MA \arrow[rrdd, "n^M"] \arrow[dd, Rightarrow, "\rho^M", shorten > = 1ex, shorten < = 2ex] &  &   \\
                                  &  &                                             &  &   \\
M \arrow[rruu, "1i"] \arrow[rrrr,equal] &  & {}                                          &  & M,
\end{tikzcd}$
\end{tabular}
\end{center}

satisfying:

\begin{enumerate}
\item [a.] We have:
\end{enumerate}

\settoheight{\prelim}{\includegraphics[width=52.5mm]{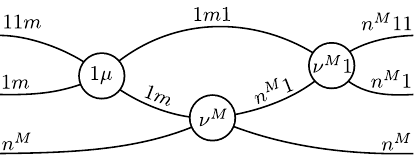}}

\begin{equation}\label{eqn:moduleassociativity}
\begin{tabular}{@{}ccc@{}}

\includegraphics[width=52.5mm]{Pictures/Preliminaries/Module/associativity1.pdf} & \raisebox{0.45\prelim}{$=$} &
\includegraphics[width=45mm]{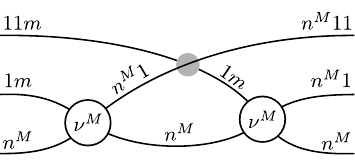},

\end{tabular}
\end{equation}

\begin{enumerate}
\item [b.] We have:
\end{enumerate}

\settoheight{\prelim}{\includegraphics[width=22.5mm]{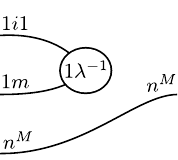}}

\begin{equation}\label{eqn:moduleunitality}
\begin{tabular}{@{}ccc@{}}

\includegraphics[width=22.5mm]{Pictures/Preliminaries/Module/unitality1.pdf} & \raisebox{0.45\prelim}{$=$} &

\includegraphics[width=37.5mm]{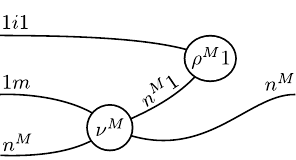}.

\end{tabular}
\end{equation}
\end{Definition}

Finally, let us recall definitions 3.2.6 and 3.2.7 of \cite{D4}.

\begin{Definition}\label{def:modulemap}
Let $M$ and $N$ be two right $A$-modules. A right $A$-module 1-morphism consists of a 1-morphism $f:M\rightarrow N$ in $\mathfrak{C}$ together with an invertible 2-morphism

$$\begin{tikzcd}[sep=small]
MA \arrow[dd, "f1"'] \arrow[rr, "n^M"]    &  & M \arrow[dd, "f"] \\
                                            &  &                      \\
NA \arrow[rr, "n^N"'] \arrow[rruu, Rightarrow, "\psi^f", shorten > = 2.5ex, shorten < = 2.5ex] &  & N,                   
\end{tikzcd}$$

subject to the coherence relations:

\begin{enumerate}
\item [a.] We have:
\end{enumerate}

\settoheight{\prelim}{\includegraphics[width=48mm]{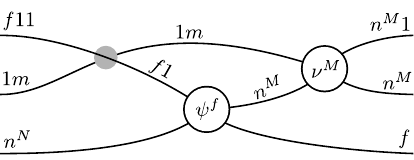}}

\begin{equation}\label{eqn:modulemapassociativity}
\begin{tabular}{@{}ccc@{}}

\includegraphics[width=48mm]{Pictures/Preliminaries/Module/map1.pdf} & \raisebox{0.45\prelim}{$=$} &

\includegraphics[width=48mm]{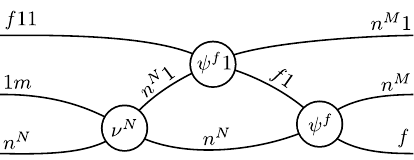},

\end{tabular}
\end{equation}

\begin{enumerate}
\item [b.] We have:
\end{enumerate}

\settoheight{\prelim}{\includegraphics[width=30mm]{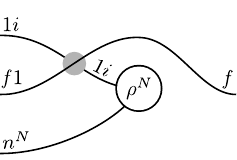}}

\begin{equation}\label{eqn:modulemapunitality}
\begin{tabular}{@{}ccc@{}}

\includegraphics[width=30mm]{Pictures/Preliminaries/Module/map3.pdf} & \raisebox{0.45\prelim}{$=$} &

\includegraphics[width=30mm]{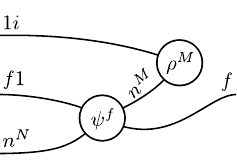}.

\end{tabular}
\end{equation}
\end{Definition}

\begin{Definition}\label{def:moduleintertwiner}
Let $M$ and $N$ be two right $A$-modules, and $f,g:M \to N$ two right $A$-module 1-morphisms. A right $A$-module 2-morphism $f\to g$ is a 2-morphism $\gamma:f\to g$ in $\mathfrak{C}$ that satisfies the following equality:

\settoheight{\prelim}{\includegraphics[width=30mm]{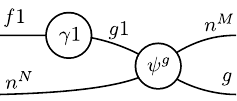}}

$$\label{eqn:module2map}\begin{tabular}{@{}ccc@{}}

\includegraphics[width=30mm]{Pictures/Preliminaries/Module/2morphism1.pdf} & \raisebox{0.45\prelim}{$=$} &

\includegraphics[width=30mm]{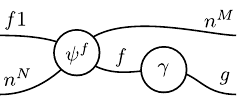}.

\end{tabular}$$
\end{Definition}

These object assemble to form a 2-category as was proven in \cite{D4}.

\begin{Lemma}
Let $A$ be an algebra in a monoidal 2-category $\mathfrak{C}$. Right $A$-modules, right $A$-module 1-morphisms, and right $A$-module 2-morphisms form a 2-category, which we denote by $\mathbf{Mod}_{\mathfrak{C}}(A)$.
\end{Lemma}

\begin{Theorem}
Let $A$ be a separable algebra in a multifusion 2-category $\mathfrak{C}$. Then $\mathbf{Mod}_{\mathfrak{C}}(A)$ is a finite semisimple 2-category.
\end{Theorem}

\begin{proof}
    See \cite[Proposition 3.1.2]{D4}.
\end{proof}

One can similarly define the notion of left modules, 1-morphisms and 2-morphisms between them. A detailed list of axioms can be found in \cite[Section 2.2]{D8} if necessary. Let $A$ be an algebra in a monoidal 2-category $\mathfrak{C}$. We denote the 2-category of left $A$-modules, left $A$-module 1-morphisms, and left $A$-module 2-morphisms by $\mathbf{Lmod}_{\mathfrak{C}}(A)$.

Let $A$ and $B$ be two algebras in monoidal 2-category $\mathfrak{C}$.

\begin{Definition} \label{def:bimodule}
    An $(A,B)$-bimodule in $\mathfrak{C}$ consists of:
    \begin{enumerate}
        \item An object $M$ in $\mathfrak{C}$;
        
        \item A left $A$-module $(M,l^M,\lambda^M,\kappa^M)$;

        \item A right $B$-module $(M,n^M,\nu^M,\rho^M)$;

        \item A balancing structure, i.e. a 2-isomorphism
        $$\begin{tikzcd}[sep=small]
        AMB \arrow[dd, "l^M1"'] \arrow[rr, "1n^M"]    &  & AM \arrow[dd, "l^M"] \arrow[ddll, Rightarrow, "\mu^M"', shorten > = 2.5ex, shorten < = 2.5ex] \\
                                                    &  &                      \\
        NB \arrow[rr, "n^M"']  &  & N,                   
        \end{tikzcd}$$ satisfying:
    \end{enumerate}
    
    \begin{itemize}
        \item [a.] We have:
    \end{itemize}

    \settoheight{\prelim}{\includegraphics[width=48mm]{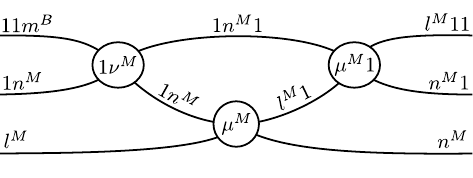}}

    \begin{equation}\label{eqn:bimoduleaxiom1}
    \begin{tabular}{@{}ccc@{}}
    
    \includegraphics[width=48mm]{Pictures/Preliminaries/Bimodule/bimodule1left.pdf} & \raisebox{0.45\prelim}{$=$} &
    
    \includegraphics[width=48mm]{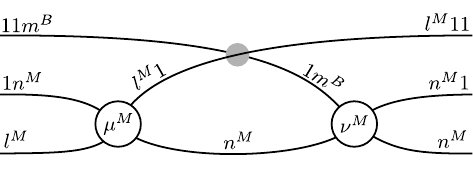},
    
    \end{tabular}
    \end{equation}

    \begin{itemize}
        \item [b.] We have:
    \end{itemize}

    \settoheight{\prelim}{\includegraphics[width=48mm]{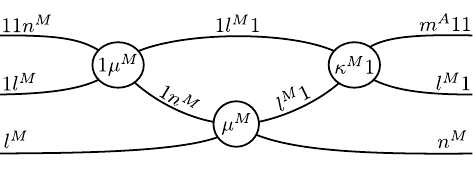}}

    \begin{equation}\label{eqn:bimoduleaxiom2}
    \begin{tabular}{@{}ccc@{}}
    
    \includegraphics[width=48mm]{Pictures/Preliminaries/Bimodule/bimodule2left.pdf} & \raisebox{0.45\prelim}{$=$} &
    
    \includegraphics[width=48mm]{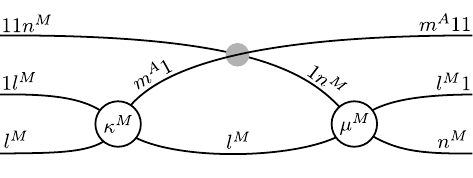},
    
    \end{tabular}
    \end{equation}

    \begin{itemize}
        \item [c.] We have:
    \end{itemize}

    \settoheight{\prelim}{\includegraphics[width=24mm]{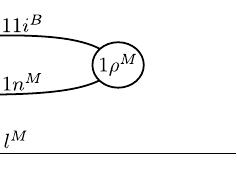}}

    \begin{equation}\label{eqn:bimoduleaxiom3}
    \begin{tabular}{@{}ccc@{}}
    
    \includegraphics[width=24mm]{Pictures/Preliminaries/Bimodule/bimodule3left.pdf} & \raisebox{0.45\prelim}{$=$} &
    
    \includegraphics[width=48mm]{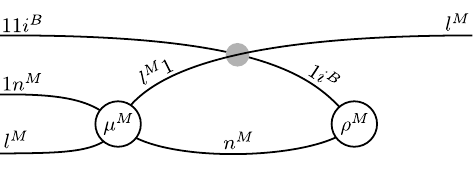},
    
    \end{tabular}
    \end{equation}

    \begin{itemize}
        \item [a.] We have:
    \end{itemize}

    \settoheight{\prelim}{\includegraphics[width=36mm]{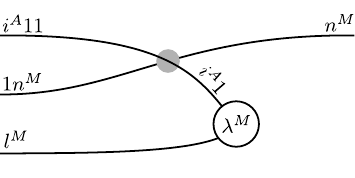}}

    \begin{equation}\label{eqn:bimoduleaxiom4}
    \begin{tabular}{@{}ccc@{}}
    
    \includegraphics[width=36mm]{Pictures/Preliminaries/Bimodule/bimodule4left.pdf} & \raisebox{0.45\prelim}{$=$} &
    
    \includegraphics[width=36mm]{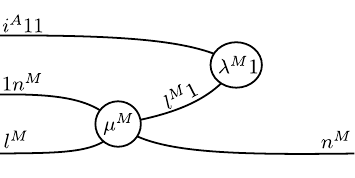}.
    
    \end{tabular}
    \end{equation}
\end{Definition}

\begin{Definition} \label{def:bimodule1morphism}
    Let $M$ and $N$ be two $(A,B)$-bimodules in $\mathfrak{C}$. An $(A,B)$-bimodule 1-morphism consists of:
    \begin{enumerate}
        \item A 1-morphism $f:M \to N$ in $\mathfrak{C}$;

        \item A 2-isomorphism $\chi^f$ such that $(f,\chi^f)$ is a left $A$-module 1-morphism;
        
        \item A 2-isomorphism $\psi^f$ such that $(f,\psi^f)$ is a right $B$-module 1-morphism;

        \item $\chi^f$ and $\psi^f$ satisfy the following condition:
    \end{enumerate}

    \settoheight{\prelim}{\includegraphics[width=48mm]{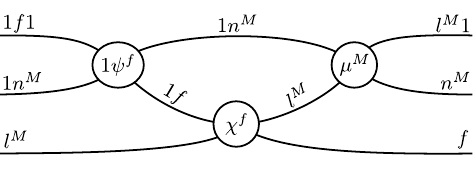}}

    \begin{equation}\label{eqn:bimodule1morphism}
    \begin{tabular}{@{}ccc@{}}
    
    \includegraphics[width=48mm]{Pictures/Preliminaries/Bimodule/bimodule1morphismleft.pdf} & \raisebox{0.45\prelim}{$=$} &
    
    \includegraphics[width=48mm]{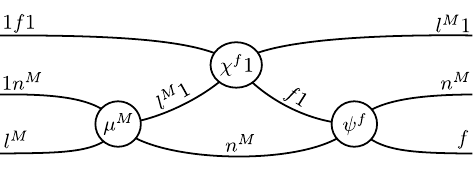}.
    
    \end{tabular}
    \end{equation}
\end{Definition}

\begin{Definition} \label{def:bimodule2morphism}
    Let $M$ and $N$ be two $(A,B)$-bimodules in $\mathfrak{C}$, and let $f$ and $g$ be two $(A,B)$-bimodule 1-morphisms from $M$ to $N$. An $(A,B)$-bimodule 2-morphism from $f$ to $g$ is a 2-morphism $\gamma:f \to g$ in $\mathfrak{C}$ such that $\gamma$ is both a left $A$-module 2-morphism and a right $B$-module 2-morphism.
\end{Definition}

\begin{Lemma}
    $(A,B)$-bimodules, $(A,B)$-bimodule 1-morphisms and $(A,B)$-bimodule 2-morphisms in $\mathfrak{C}$ form a 2-category, denoted by $\mathbf{Bimod}_\mathfrak{C}(A,B)$.
\end{Lemma}

We will mainly focus on the cases when two algebras are identical from now on. We abbreviate the 2-category of $(A,A)$-bimodules in $\mathfrak{C}$ as $\mathbf{Bimod}_\mathfrak{C}(A)$. In \cite[Section 3]{D8}, one can define a relative tensor product $\Box_A$ under some additional assumption, and then $\mathbf{Bimod}_\mathfrak{C}(A)$ can be endowed with a monoidal 2-category structure, with product given by relative tensor product $\Box_A$ and unit given by $A$ as an $(A,A)$-bimodule.

From now on, we will assume relative tensor products always exist, and they commute with monoidal product. This assumption always holds when $A$ is a separable algebra in a multifusion 2-category $\mathfrak{C}$.

\begin{Definition} \label{def:E1localmodule}
    Let us rename the monoidal 2-category $\mathbf{Bimod}_\mathfrak{C}(A)$ as $\mathbf{Mod}^{\mathbb{E}_1}_\mathfrak{C}(A)$, that is the 2-category of monoidal local  modules (or $\mathbb{E}_1$ local modules) over algebra $A$ in $\mathfrak{C}$.
\end{Definition}

\begin{Theorem}
    Suppose $A$ is a separable algebra in a multifusion 2-category $\mathfrak{C}$, then $\mathbf{Mod}^{\mathbb{E}_1}_\mathfrak{C}(A)$ is a multifusion 2-category. If $A$ is connected, then $\mathbf{Mod}^{\mathbb{E}_1}_\mathfrak{C}(A)$ is fusion.
\end{Theorem}

\begin{proof}
    See \cite[Theorem 3.2.8]{D8}.
\end{proof}

This change of terminology is by analogy with the following notions of $\mathbb{E}_n$ local modules.

Let $\mathfrak{B}$ be a semi-strict braided monoidal 2-category, and $B$ be a braided algebra in $\mathfrak{B}$.

\begin{Definition} \label{def:E2localmodule}
    Recall from \cite{DX} that a braided module (or $\mathbb{E}_2$ local module) over $B$ in $\mathfrak{B}$ consists of:
    \begin{enumerate}
    \item A right $B$-module $(M,n^M,\nu^M,\rho^M)$ in $\mathfrak{B}$,
    
    \item A $2$-isomorphism, called a holonomy,
    $$\begin{tikzcd}[sep=small]
    BM \arrow[rr, "b"]                    & {} \arrow[dd, Rightarrow, "h^{M}", shorten > = 2ex, shorten < = 1ex, near start] & MB \arrow[dd, "n^M"] \\
                                      &                        &                      \\
    MB \arrow[uu, "b"] \arrow[rr, "n^M"'] & {}                     & M,                  
    \end{tikzcd}$$
    \end{enumerate}

    satisfying:

    \begin{enumerate}
        \item [a.] We have:
    \end{enumerate}

        \settoheight{\diagramwidth}{\includegraphics[width=75mm]{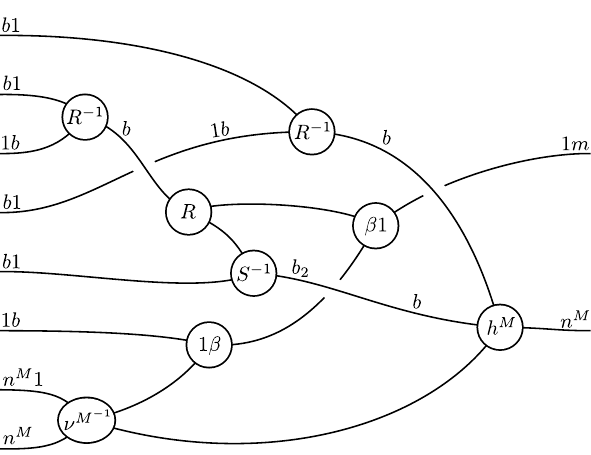}}
        
        \begin{landscape}
        \vspace*{2.2cm}
        \begin{equation}\label{eqn:holonomyassociativity}
        \begin{tabular}{@{}ccc@{}}
        \includegraphics[width=75mm]{Pictures/Preliminaries/E2LocalModule/holonomycoh1left.pdf}&\raisebox{0.48\diagramwidth}{$=$} &
        \includegraphics[width=60mm]{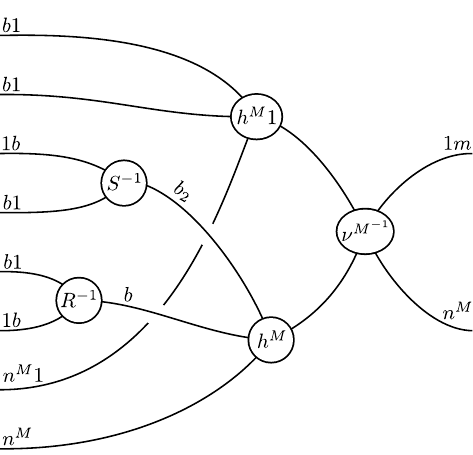},
        \end{tabular}
        \end{equation}
        \end{landscape}

    \begin{enumerate}
        \item [b.] We have:
     \end{enumerate}   
        \settoheight{\diagramwidth}{\includegraphics[width=60mm]{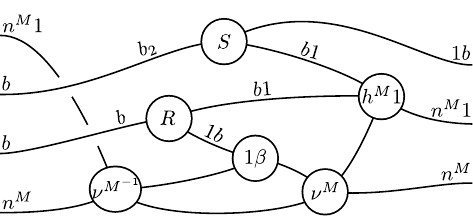}}

        \begin{equation}\label{eqn:holonomyassociativityprime}
        \begin{tabular}{@{}ccc@{}}
        \includegraphics[width=60mm]{Pictures/Preliminaries/E2LocalModule/holonomycoh1primeleft.pdf}&\raisebox{0.45\diagramwidth}{$=$} &
        \includegraphics[width=52.5mm]{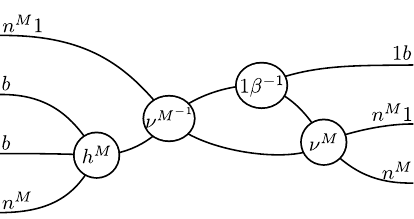},
        \end{tabular}
        \end{equation}
        
    \begin{enumerate}
        \item [c.] We have:
     \end{enumerate}   
        \settoheight{\diagramwidth}{\includegraphics[width=30mm]{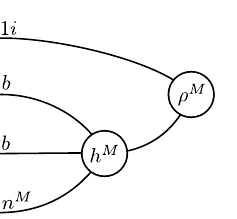}}

        \begin{equation}\label{eqn:holonomyunitality}
        \begin{tabular}{@{}ccc@{}}
        \includegraphics[width=30mm]{Pictures/Preliminaries/E2LocalModule/holonomycoh2left.pdf}&\raisebox{0.45\diagramwidth}{$=$} &
        \includegraphics[width=30mm]{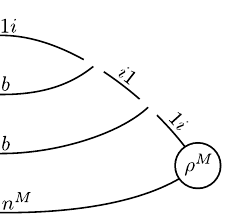}.
        \end{tabular}
        \end{equation}
\end{Definition}

\begin{Definition} \label{def:E2localmorphism}
Given $M$ and $N$ two braided local $B$-modules. A 1-morphism of braided local $B$-modules is a right $B$-module 1-morphism $(f,\psi^f)$ in $\mathfrak{B}$ satisfying the following equation
    
\settoheight{\diagramwidth}{\includegraphics[width=40mm]{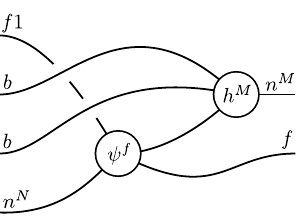}}

\begin{equation}\label{eqn:E2localmorphism}
    \begin{tabular}{@{}ccc@{}}
        \includegraphics[width=37.5mm]{Pictures/Preliminaries/E2LocalModule/localmodule1morcohleft.pdf}&\raisebox{0.45\diagramwidth}{$=$} &
        \includegraphics[width=37.5mm]{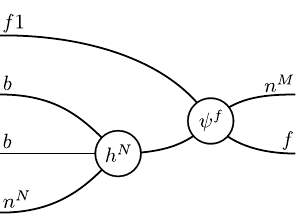}.
    \end{tabular}
\end{equation}
        
\end{Definition}

\begin{Definition} \label{def:E2local2morphism}
Given $f,g:M\rightarrow N$ two 1-morphisms of braided local $B$-modules. A 2-morphisms of braided local $B$-modules is a right $B$-module 2-morphism $f \to g$.
\end{Definition}

\begin{Lemma}
Braided local $B$-modules in $\mathfrak{B}$, 1-morphisms of braided local $B$-modules, and 2-morphisms of braided local $B$-modules form a 2-category, which we denote by $\mathbf{Mod}^{\mathbb{E}_2}_{\mathfrak{B}}(B)$.
\end{Lemma}

\begin{Theorem}
    Suppose $B$ is a separable braided algebra in a braided multifusion 2-category $\mathfrak{B}$, then $\mathbf{Mod}^{\mathbb{E}_2}_\mathfrak{B}(B)$ is a braided multifusion 2-category. If $B$ is connected, then $\mathbf{Mod}^{\mathbb{E}_2}_\mathfrak{B}(B)$ is fusion.
\end{Theorem}

\begin{proof}
    See \cite[Theorem 2.3.4]{DX}.
\end{proof}

\begin{Remark} \label{rmk:EmbeddingE2LocalModulesInDrinfeldCenter}
    By \cite[Remark 2.2.5]{DX}, one can embed $\mathbb{E}_2$ local $B$-modules into the Drinfeld center (see Definition \ref{def:DrinfeldCenter}) of the 2-category of right $B$-modules via the braided 2-functor \[\mathbf{Mod}^{\mathbb{E}_2}_\mathfrak{B}(B) \to \mathscr{Z}_1(\mathbf{Mod}_\mathfrak{B}(B)); \quad M \mapsto (M,\widetilde{b}_{M,-},\widetilde{R}_{M,-,-},\widetilde{S}_{M,-,-}).\]
\end{Remark}

Finally, let $\mathfrak{S}$ be a semi-strict sylleptic monoidal 2-category, and $S$ be a symmetric algebra in $\mathfrak{S}$.

\begin{Definition} \label{def:E3localmodule}
    A symmetric local module (or $\mathbb{E}_3$ local module) over $S$ in $\mathfrak{S}$ consists of an braided local $S$-module $(M,n^M,\nu^M,\rho^M,h^M)$ such that its holonomy satisfies with an extra coherence condition:
    \begin{center}
    \settoheight{\diagramwidth}{\includegraphics[width=30mm]{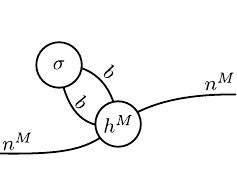}}  
    \begin{equation} \label{eqn:E3localmodule}
    \begin{tabular}{@{}cccccc@{}}
        \raisebox{0.45\diagramwidth}{$1_{n^M}$} & \raisebox{0.45\diagramwidth}{$=$} & \includegraphics[width=30mm]{Pictures/Preliminaries/E3LocalModule/E3localmodule1.pdf}&\raisebox{0.45\diagramwidth}{$=$} & \includegraphics[width=40mm]{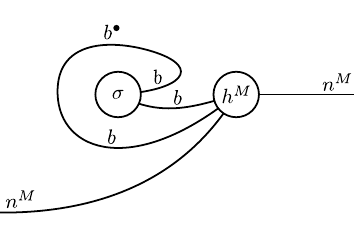}&
        \raisebox{0.45\diagramwidth}{.} 
        \end{tabular}
    \end{equation}
    \end{center}
\end{Definition}

\begin{Definition}
    Symmetric local module 1-morphisms are just braided local module 1-morphisms; symmetric local module 2-morphisms are just braided local module 2-morphisms. Thus, symmetric local modules over $S$ form a 2-subcategory of $\mathbf{Mod}^{\mathbb{E}_2}_\mathfrak{S}(S)$. Let us denote it as $\mathbf{Mod}^{\mathbb{E}_3}_\mathfrak{S}(S)$.
\end{Definition}

\begin{Theorem}
    Suppose $S$ is a separable symmetric algebra in a sylleptic multifusion 2-category $\mathfrak{S}$, then $\mathbf{Mod}^{\mathbb{E}_3}_\mathfrak{S}(S)$ is a sylleptic multifusion 2-category. If $S$ is connected, then $\mathbf{Mod}^{\mathbb{E}_3}_\mathfrak{S}(S)$ is fusion.
\end{Theorem}

\begin{proof}
    First, recall that $\mathbf{Mod}^{\mathbb{E}_2}_\mathfrak{S}(S)$ is finite semisimple by \cite[Proposition 2.3.1]{DX}. By definition, $\mathbf{Mod}^{\mathbb{E}_3}_\mathfrak{S}(S)$ is a 2-subcategory of $\mathbf{Mod}^{\mathbb{E}_2}_\mathfrak{S}(S)$, hence it is also finite semisimple. Then recall that $\mathbf{Mod}^{\mathbb{E}_2}_\mathfrak{S}(S)$ is rigid by \cite[Proposition 2.3.2]{DX}, since the monoidal structure of $\mathbf{Mod}^{\mathbb{E}_3}_\mathfrak{S}(S)$ is inherited from that on $\mathbf{Mod}^{\mathbb{E}_2}_\mathfrak{S}(S)$, by the monoidality of dual 2-functor \cite[Appendix A]{DX}, $\mathbf{Mod}^{\mathbb{E}_3}_\mathfrak{S}(S)$ is also rigid. If $S$ is connected, then $\mathbf{Mod}^{\mathbb{E}_2}_\mathfrak{S}(S)$ is a fusion 2-category and $\mathbf{Mod}^{\mathbb{E}_3}_\mathfrak{S}(S)$ becomes a fusion 2-subcategory.

    The braiding on $\mathbf{Mod}^{\mathbb{E}_2}_\mathfrak{S}(S)$ was described explicitly in \cite[Theorem 2.2.3]{DX}. It is immediate that the braiding preserves $\mathbf{Mod}^{\mathbb{E}_3}_\mathfrak{S}(S)$ as a monoidal 2-subcategory. Lastly, we can construct the syllepsis on $\mathbf{Mod}^{\mathbb{E}_3}_\mathfrak{S}(S)$ as follows.

    Let $M$ and $N$ be two $\mathbb{E}_3$ local $S$-modules, and write $$ t_{M,N}: M \, \Box \, N \to M \, \Box_S \, N, \quad t_{N,M}: N \, \Box \, M \to N \, \Box_S \, M, $$ for the universal right $S$-balanced 1-morphisms from relative tensor product \cite[Definition 3.1.3]{D8}. Recall from the proof of \cite[Theorem 2.2.3]{DX} that the braiding $(\widetilde{b},\widetilde{R},\widetilde{S})$ on $\mathbf{Mod}^{\mathbb{E}_3}_\mathfrak{S}(S)$ is induced from braiding $(b,R,S)$ in $\mathfrak{S}$ using the universal property of relative tensor product $\Box_S$:
    \[\begin{tikzcd}[sep=40pt]
        {MN}
            \arrow[dd,equal]
            \arrow[r,"t_{M,N}"]
            \arrow[rd,"b_{M,N}"']
        & {M \, \Box_S \, N}
            \arrow[rd,"\widetilde{b}_{M,N}"]
        & {}
        \\ {}
            \arrow[r,Rightarrow,shorten <=15pt,shorten >=15pt,"\sigma_{M,N}"]
        & {NM}
            \arrow[r,"t_{N,M}"]
            \arrow[ld,"b_{N,M}"']
            \arrow[u,Rightarrow,shorten <=10pt,shorten >=10pt,"\xi_{M,N}"]
        & {N \, \Box_S \, M}
            \arrow[ld,"\widetilde{b}_{N,M}"]
        \\ {MN}
            \arrow[r,"t_{M,N}"']
            \arrow[rru,Rightarrow,shorten <=60pt,shorten >=60pt,"\xi_{N,M}"]
        & {M \, \Box_S \, N}
        & {}
    \end{tikzcd}\]
    where $\xi_{M,N}$ and $\xi_{N,M}$ are $S$-balanced 2-isomorphisms created via 2-universal properties. Then we can consider the composition of the above pasting diagram, it gives a 2-isomorphism $1_{M \Box_S N} \circ t_{M,N} \to \widetilde{b}_{N,M} \circ \widetilde{b}_{M,N} \circ t_{M,N} $. Then one can adapt the proof of \cite[Proposition 3.11]{DY} to check that it is $S$-balanced. Hence, by the 2-universal property of $t_{M,N}$ again, there exists a unique 2-isomorphism $\widetilde{\sigma}_{M,N}: 1_{M \Box_S N} \to \widetilde{b}_{N,M} \circ \widetilde{b}_{M,N}$. This is the syllepsis on $\mathbf{Mod}^{\mathbb{E}_3}_\mathfrak{S}(S)$. The naturality of the relative tensor product $\Box_S$ ensures $(M,N) \mapsto \widetilde{\sigma}_{M,N}$ extends to an invertible modification.

    Finally, with some careful elaboration one can check that $\widetilde{\sigma}$ and braiding $(\widetilde{b},\widetilde{R},\widetilde{S})$ indeed satisfies the conditions for sylleptic monoidal 2-categories: (\ref{eqn:syllepsisaxiom1}), (\ref{eqn:syllepsisaxiom2}) and unitality.
\end{proof}

\begin{Remark} \label{rmk:EmbeddingE3LocalModulesInSyllepticCenter}
    One might notice in the proof of the above theorem, the existence of the half-syllepsis $\widetilde{\sigma}_{M,-}$ is solely determined by the $\mathbb{E}_3$ local $S$-module structure on $M$. Thus, one can extend the half-syllepsis $\widetilde{\sigma}_{M,N}$ to any right $S$-module $N$ rather than just $\mathbb{E}_3$ local $S$-module $N$. More generally, this provides an embedding of $\mathbf{Mod}^{\mathbb{E}_3}_\mathfrak{S}(S)$ into the sylleptic center $\mathscr{Z}_2(\mathbf{Mod}_\mathfrak{S}(S))$, where for an $\mathbb{E}_3$ local $S$-module $M$ we assign an object $(M,\widetilde{\sigma}_{M,-})$, see Definition \ref{def:SyllepticCenter}.
\end{Remark}

\subsection{Induction 2-Functors}

    In this section, we would like to recall the main results from \cite{DY}.

    \begin{Theorem} \label{thm:ModulesOfE2algebraIsE1}
        Let $B$ be a separable braided algebra in a braided multifusion 2-category $\mathfrak{B}$, then $\mathbf{Mod}_\mathfrak{B}(B)$ is equipped with a canonical multifusion 2-category structure.
    \end{Theorem}

    \begin{Theorem} \label{thm:ModulesOfE3algebraIsE2}
        Let $S$ be a separable symmetric algebra in a sylleptic multifusion 2-category $\mathfrak{S}$, then $\mathbf{Mod}_\mathfrak{S}(S)$ is equipped with a canonical braided multifusion 2-category structure.
    \end{Theorem}

    The above canonical monoidal structures on the 2-category of modules is defined by the so-called induction 2-functors.

    \begin{Lemma}
        Let $B$ be a braided algebra in a braided multifusion 2-category $\mathfrak{B}$, then there is a 2-functor $Ind^+:\mathbf{Mod}_\mathfrak{B}(B) \to \mathbf{Bimod}_\mathfrak{B}(B)$ which is fully faithful on 2-morphisms.
    \end{Lemma}

    \begin{proof}
        Given any right $B$-module $(M,n^M,\nu^M,\rho^M)$, we can induce a left $B$-action $l^M:B \, \Box \, M \xrightarrow{b_{B,M}} M \, \Box \, B \xrightarrow{n^M} M $ with 2-isomorphisms
        
        \settoheight{\braid}{\includegraphics[width=30mm]{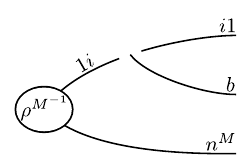}}
        
        $$\raisebox{15pt}{\raisebox{0.45\braid}{$\lambda^M:=\ $}
        \includegraphics[width=30mm]{Pictures/Preliminaries/Induction/lambdaM.pdf},}\ \ \ \settoheight{\braid}{\includegraphics[width=52.5mm]{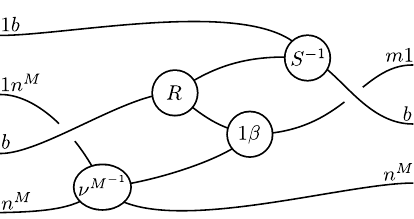}}
        \raisebox{0.45\braid}{$\kappa^M:=\ $}
        \includegraphics[width=52.5mm]{Pictures/Preliminaries/Induction/kappaM.pdf}.$$

        The left and right $B$-actions on $M$ are compatible with the balancing 2-isomorphism
        \settoheight{\braid}{\includegraphics[width=52.5mm]{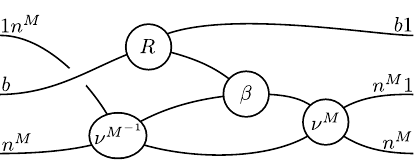}}
        
        $$\raisebox{0.45\braid}{$\beta^M:=\ $}
        \includegraphics[width=52.5mm]{Pictures/Preliminaries/Induction/betaM.pdf}.$$

        Given a right $B$-module 1-morphism $f:M \to N$, we induce a left $B$-module 1-morphism structure
        \settoheight{\braid}{\includegraphics[width=37.5mm]{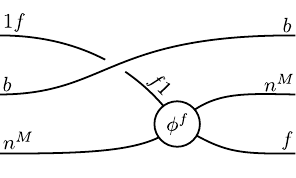}}
        
        $$\raisebox{0.45\braid}{$\chi^f:=\ $}
        \includegraphics[width=37.5mm]{Pictures/Preliminaries/Induction/chif.pdf}.$$
        Finally, it is routine to check that the above structures satisfies the corresponding coherence conditions. This gives us a 2-functor $Ind^+:\mathbf{Mod}_\mathfrak{B}(B) \to \mathbf{Bimod}_\mathfrak{B}(B)$. Then it is straight forward to see this 2-functor is fully faithful on 2-morphisms.
    \end{proof}

    \begin{Remark}
        One can replace the induced $B$-action on $M$ by $B \, \Box \, M \xrightarrow{b^\bullet_{M,B}} M \, \Box \, B \xrightarrow{n^M} M$ and define the other coherence data by replacing braiding $(b,R,S)$ by its reverse $(b^\bullet,S^\bullet,R^\bullet)$. The resulting 2-functor $Ind^-:\mathbf{Mod}_\mathfrak{B}(B) \to \mathbf{Bimod}_\mathfrak{B}(B)$ is again fully faithful on 2-morphisms.
    \end{Remark}

    \begin{Corollary} \label{cor:MonoidalStructureOnModulesByInduction2Functor}
        Moreover, suppose $B$ is a separable braided algebra in $\mathfrak{B}$, then $\mathbf{Bimod}_\mathfrak{B}(B)$ is equipped with a monoidal product given by the relative tensor product $\Box_B$. Hence, the induction 2-functors $Ind^\pm:\mathbf{Mod}_\mathfrak{B}(B) \to \mathbf{Bimod}_\mathfrak{B}(B)$ induce (non-equivalent) monoidal structures on 2-category $\mathbf{Mod}_\mathfrak{B}(B)$.
    \end{Corollary}

    We will denote the monoidal structures on $\mathbf{Mod}_\mathfrak{B}(B)$ induced by $Ind^\pm$ as $\Box^\pm_B$, respectively. For consistency, when we refer to the monoidal structure on $\mathbf{Mod}_\mathfrak{B}(B)$, unless otherwise specified, we will always refer to the monoidal structure $\Box^+_B$.

    An analogous construction exists for separable symmetric algebras within a sylleptic multifusion 2-category.

    \begin{Lemma}
        Let $S$ be a separable symmetric algebra in a sylleptic multifusion 2-category $\mathfrak{S}$, then there is a monoidal 2-functor $Ind^+:\mathbf{Mod}_\mathfrak{S}(S) \to \mathbf{Mod}^{\mathbb{E}_2}_\mathfrak{S}(S)$ which is an embedding (i.e. it induces equivalences on the level of hom 1-categories).
    \end{Lemma}

    \begin{proof}
        Given any right $S$-module $(M,n^M,\nu^M,\rho^M)$, we can induce a holonomy
        \settoheight{\braid}{\includegraphics[width=30mm]{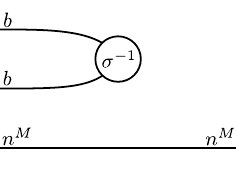}}
        
        $$\raisebox{15pt}{\raisebox{0.45\braid}{$h^M:=\ $}
        \includegraphics[width=30mm]{Pictures/Preliminaries/Induction/holonomyM.pdf}.}$$ Applying properties of syllepsis: (\ref{eqn:syllepsisaxiom1}), (\ref{eqn:syllepsisaxiom2}) and unitality we can show that coherence conditions (\ref{eqn:holonomyassociativity}), (\ref{eqn:holonomyassociativityprime}) and (\ref{eqn:holonomyunitality}) are satisfied. Coherence conditions on 1- and 2-morphisms are automatic, hence we obtain a 2-functor $Ind^+:\mathbf{Mod}_\mathfrak{S}(S) \to \mathbf{Mod}^{\mathbb{E}_2}_\mathfrak{S}(S)$. Then it is straight-forward to check this is an embedding.

        Finally, the monoidal structure on $\mathbf{Mod}^{\mathbb{E}_2}_\mathfrak{S}(S)$ is by definition (c.f. \cite[Proposition 2.2.1]{DX}) induced from the monoidal structure $\Box^+_S$ on the underlying right $S$-modules, hence the above 2-functor is equipped with a canonical monoidal 2-functor structure.
    \end{proof}

    \begin{Remark}
        We can also replace the syllepsis $\sigma$ in the above lemma by its reverse $\overline{\sigma}$. This produces another induction 2-functor $Ind^-:\mathbf{Mod}_\mathfrak{S}(S) \to \mathbf{Mod}^{\mathbb{E}_2}_\mathfrak{S}(S)$, which is again an embedding of monoidal 2-categories.
    \end{Remark}

    \begin{Corollary} \label{cor:BraidingOnModulesByInduction2Functor}
        Induction 2-functors $Ind^\pm:\mathbf{Mod}_\mathfrak{S}(S) \to \mathbf{Mod}^{\mathbb{E}_2}_\mathfrak{S}(S)$ endow monoidal 2-category $\mathbf{Mod}_\mathfrak{S}(S)$ with (non-equivalent) braidings.
    \end{Corollary}

    Compared with \cite[Theorem 3.8]{DY}, the braiding there was exactly induced by induction 2-functor $Ind^+:\mathbf{Mod}_\mathfrak{S}(S) \to \mathbf{Mod}^{\mathbb{E}_2}_\mathfrak{S}(S)$. Let us denote the braidings induced by induction 2-functor $Ind^\pm$ as $\widetilde{b}^\pm$. Unless otherwise specified, we will always use braiding $\widetilde{b}^+$ on $\mathbf{Mod}_\mathfrak{S}(S)$.

    Lastly, one can obtain an equivalent characterizations of $\mathbb{E}_n$ local module structures via induction 2-functors.

    \begin{Lemma}
        Suppose $B$ is a braided algebra in braided monoidal 2-category $\mathfrak{B}$ and $M$ is a right $B$-module. Then the $\mathbb{E}_2$ local $B$-module structure on $M$, i.e. a holonomy on $M$, corresponds to the data of upgrading the right $B$-module 1-morphism $1_M: M \to M$ to an $\mathbb{E}_1$ local $B$-module 1-morphism $Ind^+(M) \to Ind^-(M)$, for induction 2-functors $Ind^\pm:\mathbf{Mod}_\mathfrak{B}(B) \to \mathbf{Mod}^{\mathbb{E}_1}_\mathfrak{B}(B)$.
    \end{Lemma}

    \begin{proof}
        See \cite[Lemma 2.1.5]{DX}.
    \end{proof}
    
    \begin{Lemma}
        Suppose $S$ is a symmetric algebra in sylleptic monoidal 2-category $\mathfrak{S}$ and $M$ is a right $S$-module. Then the $\mathbb{E}_3$ local $S$-module structure on $M$ corresponds to the data of upgrading the right $S$-module 1-morphism $1_M: M \to M$ to an $\mathbb{E}_2$ local $S$-module 1-morphism $Ind^+(M) \to Ind^-(M)$, for induction 2-functors $Ind^\pm:\mathbf{Mod}_\mathfrak{S}(S) \to \mathbf{Mod}^{\mathbb{E}_2}_\mathfrak{S}(S)$.
    \end{Lemma}

    \begin{proof}
        This follows immediately from the Definition \ref{def:E3localmodule}.
    \end{proof}

\subsection{Module 2-Categories and Enrichment}

In this section, we would like to recall some basic notions and properties for modules 2-categories and enriched 2-categories. For more details, readers may consult \cite{GS,D4,D8,D9}. Here we do not explicitly assume existence of linear 2-category structures. However, readers can readily reconstruct the correct linear 2-categorical concepts by substituting specific adjectives in appropriate locations.

Let $\mathfrak{C}$ and $\mathfrak{D}$ be two monoidal 2-categories.
\begin{Definition} \label{def:Module2Category}
    A left $\mathfrak{C}$-module 2-category consists of:
    \begin{enumerate}
        \item[a.] An underlying 2-category $\mathfrak{M}$;

        \item[b.] A monoidal 2-functor $F:\mathfrak{C} \to \mathbf{End}(\mathfrak{M})$, where $\mathbf{End}(\mathfrak{M})$ is the 2-category of endo-2-functors on 2-category $\mathfrak{M}$ and the monoidal structure is given by composition of 2-functors.
    \end{enumerate}
    Similarly, a right $\mathfrak{D}$-module 2-category consists of:
    \begin{enumerate}
        \item[a.] An underlying 2-category $\mathfrak{N}$;

        \item[b.] A monoidal 2-functor $G:\mathfrak{D}^{1mp} \to \mathbf{End}(\mathfrak{N})$.
    \end{enumerate}
    A $(\mathfrak{C},\mathfrak{D})$-bimodule 2-category consists of:
    \begin{enumerate}
        \item[a.] An underlying 2-category $\mathfrak{L}$;

        \item[b.] A monoidal 2-functor $H:\mathfrak{C} \times \mathfrak{D}^{1mp} \to \mathbf{End}(\mathfrak{L})$.
    \end{enumerate}
\end{Definition}

\begin{Proposition} \label{prop:Monoidal2CategoryAction}
    Equivalently, a left $\mathfrak{C}$-module 2-category consists of:
    \begin{enumerate}
        \item A 2-functor $\Diamond: \mathfrak{C} \times \mathfrak{M} \to \mathfrak{M}$;

        \item 2-natural equivalences \[\alpha^\mathfrak{M}_{x,y,m}:  (x \, \Box \, y) \, \Diamond \, m \to x \, \Diamond \, (y \, \Diamond \, m),\] \[l^\mathfrak{M}_m: I \, \Diamond \, m \to m,\] given on objects $x,y$ in $\mathfrak{C}$ and $m$ in $\mathfrak{M}$;

        \item Invertible modifications
        \[\begin{tikzcd}[row sep=30pt,column sep=50pt]
             {((x \, \Box \, y) \, \Box \, z) \, \Diamond \, m}
                \arrow[r,"\alpha^\mathfrak{M}_{xy,z,m}"]
                \arrow[d,"\alpha^\mathfrak{C}_{x,y,z} 1"']
            & {(x \, \Box \, y) \, \Diamond \, (z \, \Diamond \, m)}
                \arrow[dd,"\alpha^\mathfrak{M}_{x,y,zm}"]
            \\ {(x \, \Box \, (y \, \Box \, z)) \, \Diamond \, m}
                \arrow[d,"\alpha^\mathfrak{M}_{x,yz,m}"']
            & {}
            \\ {x \, \Diamond \, ((y \, \Box \, z) \, \Diamond \, m)}
                \arrow[r,"1 \alpha^\mathfrak{M}_{y,z,m}"']
                \arrow[uur,Rightarrow,shorten <=40pt, shorten >=50pt,"\pi^\mathfrak{M}_{x,y,z,m}"']
            & {x \, \Diamond \, (y \, \Diamond \, (z \, \Diamond \, m))}
        \end{tikzcd},\]
        \[\begin{tikzcd}[sep=30pt]
            {x \, \Diamond \, m}
                \arrow[d,equal]
                \arrow[dr,Rightarrow,shorten <=20pt, shorten >=20pt,"\lambda^\mathfrak{M}_{x,m}"]
            & {(I \, \Box \, x) \, \Diamond \, m}
                \arrow[l,"l^\mathfrak{C}_x 1"']
                \arrow[d,"\alpha^\mathfrak{M}_{I,x,m}"]
            \\ {x \, \Diamond \, m}
            & {I \, \Diamond \, (x \, \Diamond \, m)}
                \arrow[l,"l^\mathfrak{M}_{xm}"]
        \end{tikzcd}, \quad \begin{tikzcd}[sep=30pt]
            {x \, \Diamond \, m}
                \arrow[d,equal,"" {name=A}]
            & {x \, \Diamond \, (I \, \Diamond \, m)}
                \arrow[l,"1 l^\mathfrak{M}_m"']
            \\ {x \, \Diamond \, m}
                \arrow[r,"r^\mathfrak{C}_x 1"']
            & {(x \, \Box \, I) \, \Diamond \, m}
                \arrow[u,"\alpha^\mathfrak{M}_{x,I,m}"' {name=B}]
                \arrow[Rightarrow,from=B,to=A,shorten <=30pt, shorten >=20pt,"\mu^\mathfrak{M}_{x,m}"']
        \end{tikzcd}.\]
    \end{enumerate}
    subject to conditions a), b), c) in \cite[Definition 2.1.3]{D4}.
\end{Proposition}

\begin{Remark}
    Similarly, one can define a right $\mathfrak{D}$-module and $(\mathfrak{C},\mathfrak{D})$-bimodule 2-categories using action 2-functors together with various associativity and unitality coherence data. We omit details here.
\end{Remark}

\begin{Example}
    For any algebra $A$ in monoidal 2-category $\mathfrak{C}$, the 2-category of right $A$-modules $\mathbf{Mod}_\mathfrak{C}(A)$ has a canonical left $\mathfrak{C}$-module 2-category structure via tensoring from the left. Similarly, the 2-category of left $A$-modules $\mathbf{Lmod}_\mathfrak{C}(A)$ has a canonical right $\mathfrak{C}$-module 2-category structure via tensoring from the right.
\end{Example}

\begin{Definition} \label{def:Module2Functor}
    Suppose $\mathfrak{C}$ is a monoidal 2-category. A left $\mathfrak{C}$-module 2-functor between left $\mathfrak{C}$-module 2-categories $\mathfrak{M}$ and $\mathfrak{N}$ consists of:
    \begin{itemize}
        \item An underlying 2-functor $F:\mathfrak{M} \to \mathfrak{N}$;

        \item 2-natural equivalence $\chi^F_{x,m}:x \, \Diamond^\mathfrak{N} \, F(m) \to F(x \, \Diamond^\mathfrak{M} \, m)$, given on object $x$ in $\mathfrak{C}$ and $m$ in $\mathfrak{M}$;

        \item Invertible modifications \[\begin{tikzcd}[row sep=30pt,column sep=25pt]
            {(x \, \Box \, y) \, \Diamond^\mathfrak{N} \, F(m)}
                \arrow[r,"\alpha^\mathfrak{N}_{x,y,F(m)}"]
                \arrow[d,"\chi^F_{xy,m}"']
            & {x \, \Diamond^\mathfrak{N} \, (y \, \Diamond^\mathfrak{N} \, F(m))}
                \arrow[r,"1 \chi^F_{y,m}"]
            & {x \, \Diamond^\mathfrak{N} \, F(y \, \Diamond^\mathfrak{M} \, m)}
                \arrow[d,"\chi^F_{x,F(ym)}"]
                \arrow[dll,Rightarrow,shorten <=40pt, shorten >=30pt,"\omega^F_{x,y,m}"']
            \\ {F((x \, \Box \, y) \, \Diamond^\mathfrak{M} \, m)}
                \arrow[rr,"F(\alpha^\mathfrak{M}_{x,y,m})"']
            & {}
            & {F(x \, \Diamond^\mathfrak{M} \, (y \, \Diamond^\mathfrak{M} \, m))}
        \end{tikzcd},\]
        \[\begin{tikzcd}[sep=30pt]
            {I^\mathfrak{N} \, \Diamond^\mathfrak{N} \, F(m)}
                \arrow[d,"\chi^F_{I,m}"']
                \arrow[r,equal]
            & {I^\mathfrak{N} \, \Diamond^\mathfrak{N} \, F(m)}
                \arrow[d,"l^\mathfrak{N}_{F(m)}"]
            \\ {F(I^\mathfrak{M} \, \Diamond^\mathfrak{M} \, m)}
                \arrow[r,"F(l^\mathfrak{M}_m)"']
                \arrow[ur,Rightarrow,shorten <=10pt, shorten >=20pt,"\gamma^F_m"]
            & {F(m)}
        \end{tikzcd},\] given on objects $x,y$ in $\mathfrak{C}$ and $m$ in $\mathfrak{M}$,
    \end{itemize}
    subject to conditions a), b), c) in \cite[Definition 2.1.7]{D4}.
\end{Definition}

\begin{Definition} \label{def:Module2NatTrans}
    Let $F,G$ be two left $\mathfrak{C}$-module 2-functors between $\mathfrak{M}$ and $\mathfrak{N}$. A left $\mathfrak{C}$-module 2-natural transform from $F$ to $G$ consists of:
    \begin{enumerate}
        \item An underlying 2-natural transform $\theta:F \to G$;

        \item An invertible modification \[\begin{tikzcd}[sep=30pt]
            {x \, \Diamond^\mathfrak{N} \, F(m)}
                \arrow[r,"1 \theta_m"]
                \arrow[d,"\chi^F"']
            & {x \, \Diamond^\mathfrak{N} \, G(m)}
                \arrow[d,"\chi^G"]
                \arrow[ld,Rightarrow,shorten <=20pt, shorten >=20pt,"\Pi^\theta_{x,m}"']
            \\ {F(x \, \Diamond^\mathfrak{M} \, m)}
                \arrow[r,"\theta_{xm}"']
            & {G(x \, \Diamond^\mathfrak{M} \, m)}
        \end{tikzcd},\] given on object $x$ in $\mathfrak{C}$ and $m$ in $\mathfrak{M}$,
    \end{enumerate}
    subject to conditions a), b) in \cite[Definition 2.1.9]{D4}.
\end{Definition}

\begin{Definition} \label{def:ModuleModification}
    Let $\theta,\xi$ be two left $\mathfrak{C}$-module 2-natural transforms between $F,G$, which are two left $\mathfrak{C}$-module 2-functors between $\mathfrak{M}$ and $\mathfrak{N}$. A left $\mathfrak{C}$-module modification from $\theta$ to $\xi$ consists of a modification $\Xi:\theta \to \xi$ satisfying the condition in \cite[Definition 2.1.10]{D4}.
\end{Definition}

\begin{Proposition}
    For a given monoidal 2-category $\mathfrak{C}$, we obtain a 3-category consists of left $\mathfrak{C}$-module 2-categories, left $\mathfrak{C}$-module 2-functors, left $\mathfrak{C}$-module 2-natural transforms and left $\mathfrak{C}$-module modifications, denoted as $\mathbf{Lmod}(\mathfrak{C})$.
\end{Proposition}

\begin{Remark}
    One can similarly define right module and bimodule 2-natural transforms, modifications. They form 3-categories in the same way as stated in the proposition above. 
    
    We denote the 3-category of right $\mathfrak{D}$-module 2-categories as $\mathbf{Mod}(\mathfrak{D})$, and denote the 3-category of $(\mathfrak{C},\mathfrak{D})$-bimodule 2-categories as $\mathbf{Bimod}(\mathfrak{C},\mathfrak{D})$, for any monoidal 2-categories $\mathfrak{C}$ and $\mathfrak{D}$.
\end{Remark}

\begin{Definition} \label{def:Enriched2Category}
    An enriched 2-category over a monoidal 2-category $\mathfrak{C}$ consists of:
    \begin{enumerate}
        \item A set\footnote{We assume smallness in general and ignore the size issue.} of objects $Ob(\mathfrak{M})$;

        \item For objects $x,y$ in $Ob(\mathfrak{M})$, an enriched $Hom$-object $[x,y]^\mathfrak{M}_\mathfrak{C}$ in $\mathfrak{C}$;

        \item For object $x$ in $Ob(\mathfrak{M})$, an enriched unit $j_x:I \to [x,y]^\mathfrak{M}_\mathfrak{C}$;

        \item For objects $x,y,z$ in $Ob(\mathfrak{M})$, an enriched multiplication \[m_{x,y,z}:[y,z]^\mathfrak{M}_\mathfrak{C} \,\Box \, [x,y]^\mathfrak{M}_\mathfrak{C} \to [x,z]^\mathfrak{M}_\mathfrak{C};\]

        \item For objects $x,y$ in $Ob(\mathfrak{M})$, enriched unitors in $\mathfrak{C}$ \[\begin{tikzcd}[sep=30pt]
            {I \, \Box \, [x,y]^\mathfrak{M}_\mathfrak{C}}
                \arrow[r,equal]
                \arrow[d,"j_y 1"']
            & {I \, \Box \, [x,y]^\mathfrak{M}_\mathfrak{C}}
                \arrow[d,"l_{[x,y]}"]
            \\ {[y,y]^\mathfrak{M}_\mathfrak{C} \, \Box \, [x,y]^\mathfrak{M}_\mathfrak{C}}
                \arrow[r,"m_{x,y,y}"']
                \arrow[ur,Rightarrow,shorten <=20pt, shorten >=20pt,"\rho_{x,y}"]
            & {[x,y]^\mathfrak{M}_\mathfrak{C}}
        \end{tikzcd},\] \[\begin{tikzcd}[sep=30pt]
            {[x,y]^\mathfrak{M}_\mathfrak{C} \, \Box \, I}
                \arrow[r,equal]
                \arrow[d,"1 j_x"']
            & {[x,y]^\mathfrak{M}_\mathfrak{C} \, \Box \, I}
                \arrow[d,"r_{[x,y]}"]
            \\ {[x,y]^\mathfrak{M}_\mathfrak{C} \, \Box \, [x,x]^\mathfrak{M}_\mathfrak{C}}
                \arrow[r,"m_{x,x,y}"']
                \arrow[ur,Rightarrow,shorten <=20pt, shorten >=20pt,"\lambda_{x,y}"]
            & {[x,y]^\mathfrak{M}_\mathfrak{C}}
        \end{tikzcd};\]

        \item For objects $x,y,z,w$ in $Ob(\mathfrak{M})$, an enriched associator in $\mathfrak{C}$ \[\begin{tikzcd}[sep=30pt]
            {([z,w]^\mathfrak{M}_\mathfrak{C} \, \Box \, [y,z]^\mathfrak{M}_\mathfrak{C}) \, \Box \, [x,y]^\mathfrak{M}_\mathfrak{C}}
                \arrow[r,"m_{y,z,w} 1"]
                \arrow[d,"\alpha_{[z,w],[y,z],[x,y]}"']
             & {[y,w]^\mathfrak{M}_\mathfrak{C} \, \Box \, [x,y]^\mathfrak{M}_\mathfrak{C}}
                \arrow[dd,"m_{[x,y,w]}"]
             \\{[z,w]^\mathfrak{M}_\mathfrak{C} \, \Box \,( [y,z]^\mathfrak{M}_\mathfrak{C} \, \Box \, [x,y]^\mathfrak{M}_\mathfrak{C})}
                \arrow[d,"1 m_{x,y,z}"']
                \arrow[r,Rightarrow,shorten <=10pt, shorten >=10pt,"\pi_{x,y,z,w}"]
            & {}
             \\ {[z,w]^\mathfrak{M}_\mathfrak{C} \, \Box \, [x,z]^\mathfrak{M}_\mathfrak{C}}
                \arrow[r,"m_{x,z,w}"']
             & {[x,w]^\mathfrak{M}_\mathfrak{C}}
        \end{tikzcd};\]
    \end{enumerate}
    subject to conditions a) and b) in \cite[Definition 4.2.1]{D4}.
\end{Definition}

\begin{Remark}
    Post-composing the enriched $Hom$-object with the forgetful 2-functor $\mathfrak{C} \to \mathbf{Cat}; \, x \mapsto Hom_\mathfrak{C}(I,x)$, we obtain an enriched 2-category structure over $\mathbf{Cat}$. This is nothing else but the underlying 2-category structure on $\mathfrak{M}$, with $Ob(\mathfrak{M})$ as its genuine set of objects now. Hence, we can also view an enriched 2-category over $\mathfrak{C}$ as a 2-category $\mathfrak{M}$ equipped with addition structures described in the above definition.
\end{Remark}

For finite semisimple module 2-categories, \cite{D4} has proven the equivalence between module 2-category structures and enriched 2-category structures.
\begin{Proposition}
    A finite semisimple module 2-category over a multifusion 2-category has a canonical enriched 2-category structure. Moreover, every finite semisimple module 2-category over this multifusion 2-category can be reconstructed from a rigid algebra.
\end{Proposition}

\begin{Remark}
    In general, it is not true that every enriched 2-category $\mathfrak{M}$ over a monoidal 2-category $\mathfrak{C}$ can be realized as $\mathbf{Mod}_\mathfrak{C}(A)$ for some rigid algebra $A$ in $\mathfrak{C}$. Rather, provided that there is a tensor-hom 2-adjunction \[Hom_\mathfrak{M}(x \, \Diamond \, m,n) \simeq Hom_\mathfrak{C}(x,[m,n]_\mathfrak{C})\] given on object $x$ in $\mathfrak{C}$ and objects $m,n$ in $\mathfrak{M}$, the hypothetical left $\mathfrak{C}$-action $x \mapsto x \, \Diamond \, - $ turns out to be only an \textit{oplax monoidal 2-functor} from $\mathfrak{C}$ to $\mathbf{End}(\mathfrak{M})$. We may call the enriched structure on $\mathfrak{M}$ \textit{strong} if this induced $\mathfrak{C}$-action is \textit{strong}. 
    
    Thence, this implies that a finite semisimple 2-category $\mathfrak{M}$ is strongly enriched over multifusion 2-category $\mathfrak{C}$ if and only if we can find some rigid algebra $A$ in $\mathfrak{C}$ such that $\mathfrak{M} \simeq \mathbf{Mod}_\mathfrak{C}(A)$.
\end{Remark}

\begin{Remark}
    Using the above equivalence, one can get an impression of what the \textit{correct definition} of enriched 2-functors, enriched 2-natural transforms and enriched modifications should look like. In other word, we expect an independent set of definitions for enriched 2-functors, enriched 2-natural transforms and enriched modifications between 2-categories strongly enriched over a fixed monoidal 2-category $\mathfrak{C}$, and they should together form a sub-3-category among all $\mathfrak{C}$-enriched 2-categories, $\mathbf{Enrich}(\mathfrak{C})$. Then the above proposition can be lifted to an equivalence between $\mathbf{Lmod}(\mathfrak{C})$ and its image in $\mathbf{Enrich}(\mathfrak{C})$.
\end{Remark}

\begin{Definition}
    Let $\mathfrak{B}$ be a braided 2-category. A monoidal $\mathfrak{B}$-module 2-category consists of:
    \begin{enumerate}
        \item[a.] A monoidal 2-category $\mathfrak{C}$;

        \item[b.] A braided 2-functor $F:\mathfrak{B} \to \mathscr{Z}_1(\mathfrak{C})$, where $\mathscr{Z}_1(\mathfrak{C})$ is the Drinfeld center of $\mathfrak{C}$ defined in Definition \ref{def:DrinfeldCenter}.
    \end{enumerate}
    Let $\mathfrak{T}$ be a sylleptic 2-category. A braided $\mathfrak{T}$-module 2-category consists of:
    \begin{enumerate}
        \item[a.] A braided 2-category $\mathfrak{A}$;

        \item[b.] A sylleptic 2-functor $G:\mathfrak{T} \to \mathscr{Z}_2(\mathfrak{A})$, where $\mathscr{Z}_2(\mathfrak{A})$ is the sylleptic center of $\mathfrak{A}$ defined in Definition \ref{def:SyllepticCenter}.
    \end{enumerate}
    Let $\mathfrak{V}$ be a symmetric 2-category. A sylleptic $\mathfrak{V}$-module 2-category consists of:
    \begin{enumerate}
        \item[a.] A sylleptic 2-category $\mathfrak{S}$;

        \item[b.] A sylleptic 2-functor $H:\mathfrak{V} \to \mathscr{Z}_3(\mathfrak{S})$, where $\mathscr{Z}_3(\mathfrak{S})$ is the symmetric center of $\mathfrak{S}$ defined in Definition \ref{def:SymmetricCenter}.
    \end{enumerate}
\end{Definition}

\begin{Remark}
    Post-composed with the canonical forgetful 2-functors $\mathscr{Z}_1(\mathfrak{C}) \to \mathbf{End}(\mathfrak{C})$, $\mathscr{Z}_2(\mathfrak{A}) \to \mathbf{End}(\mathfrak{A})$ and $\mathscr{Z}_3(\mathfrak{S}) \to \mathbf{End}(\mathfrak{S})$, one can obtain genuine \textit{underlying module 2-category} structures on $\mathfrak{C}$, $\mathfrak{A}$ and $\mathfrak{S}$, respectively. One should expect that:
    \begin{enumerate}
        \item A monoidal $\mathfrak{B}$-module 2-category $\mathfrak{C}$ is the same as an algebra in the monoidal 3-category $\mathbf{Lmod}(\mathfrak{B})$;

        \item A braided $\mathfrak{T}$-module 2-category $\mathfrak{A}$ is the same as a braided algebra (aka. $\mathbb{E}_2$ algebra) in the braided 3-category $\mathbf{Lmod}(\mathfrak{T})$;

        \item A sylleptic $\mathfrak{V}$-module 2-category $\mathfrak{S}$ is the same as a sylleptic algebra (aka. $\mathbb{E}_3$ algebra) in the sylleptic (in fact symmetric) 3-category $\mathbf{Lmod}(\mathfrak{V})$.
    \end{enumerate}

    Therefore, there are two natural types of higher morphisms we can impose upon them:
    \begin{enumerate}
        \item Algebra higher morphisms, i.e.
        \begin{itemize}
            \item Algebra 1-, 2- and 3-morphisms in $\mathbf{Lmod}(\mathfrak{B})$;

            \item Braided algebra 1-, 2- and 3-morphisms in $\mathbf{Lmod}(\mathfrak{T})$;

            \item Sylleptic algebra 1-, 2- and 3-morphisms in $\mathbf{Lmod}(\mathfrak{V})$;
        \end{itemize}
        These all provide 3-categories;
        \item Higher bimodules between them, i.e.
            \begin{itemize}
            \item Bimodule 1-, 2- and 3-morphisms in $\mathbf{Lmod}(\mathfrak{B})$;

            \item Monoidal bimodules, bimodule between monoidal bimodules, and then bimodule 1-, 2- and 3-morphisms in $\mathbf{Lmod}(\mathfrak{T})$;

            \item Braided bimodules, and then monoidal bimodules between them, then again iterated bimodules between these monoidal bimodules, and finally bimodule 1-, 2- and 3-morphisms in $\mathbf{Lmod}(\mathfrak{V})$;
        \end{itemize}
        They provide higher Morita 4-, 5-, and 6-categories, respectively.
    \end{enumerate}
\end{Remark}

\begin{Remark} \label{rmk:E1E2E3Enriched2Categories}
    Following \cite{Kel,For,KYZZ}, one should expect that: \begin{enumerate}
        \item A monoidal 2-category can be enriched over a braided 2-category $\mathfrak{B}$;

        \item A braided 2-category can be enriched over a sylleptic 2-category $\mathfrak{T}$;

        \item A sylleptic 2-category can be enriched over a symmetric 2-category $\mathfrak{V}$.
    \end{enumerate} Moreover, \begin{enumerate}
        \item The 3-category of $\mathfrak{B}$-enriched 2-categories is equipped with a monoidal structure, and an algebra in $\mathbf{Enrich}(\mathfrak{B})$ is the same as a monoidal 2-category enriched over $\mathfrak{B}$;

        \item The 3-category of $\mathfrak{T}$-enriched 2-categories is equipped with a braided monoidal structure, and a braided algebra in $\mathbf{Enrich}(\mathfrak{T})$ is the same as a braided 2-category enriched over $\mathfrak{T}$;

        \item The 3-category of $\mathfrak{V}$-enriched 2-categories is equipped with a sylleptic (in fact symmetric) monoidal structure, and a sylleptic algebra in $\mathbf{Enrich}(\mathfrak{V})$ is the same as a sylleptic 2-category enriched over $\mathfrak{V}$.
    \end{enumerate}
\end{Remark}

\begin{Remark}
    By analogy with the main result of \cite{D4}, we conjecture equivalences:
    \begin{enumerate}
        \item Between multifusion $\mathfrak{B}$-module 2-categories and multifusion 2-categories strongly enriched over $\mathfrak{B}$, where $\mathfrak{B}$ is a braided multifusion 2-category;

        \item Between braided multifusion $\mathfrak{T}$-module 2-categories and braided multifusion 2-categories strongly enriched over $\mathfrak{T}$, where $\mathfrak{T}$ is a sylleptic multifusion 2-category;

        \item Between sylleptic multifusion $\mathfrak{V}$-module 2-categories and sylleptic multifusion 2-categories strongly enriched over $\mathfrak{V}$ where $\mathfrak{V}$ is a symmetric multifusion 2-category.
    \end{enumerate}
\end{Remark}

\section{Monoidal Centralizers and Braided Centers}

In this section, we will continue using notations and conventions defined in the preliminary.

\subsection{Definitions and Basic Properties}

The following definitions generalize the notion of center for monoidal 2-categories \cite[Section 3.1]{Cr}.

\begin{Definition} \label{def:DrinfeldCentralizer}
    Given a monoidal 2-functor $F:\mathfrak{C} \to \mathfrak{D}$ between monoidal 2-categories $\mathfrak{C}$ and $\mathfrak{D}$, its monoidal centralizer (or Drinfeld centralizer, $\mathbb{E}_1$ centralizer) is a monoidal 2-category $\mathscr{Z}_1(F)$ where
    \begin{itemize}
        \item [a.] An Object is a triple $(x,b^x_{-},R^x_{-,-})$, where 
        \begin{itemize}
            \item $x$ is an object in $\mathfrak{D}$;

            \item $b^x_{-}$ is a pseudo-natural isomorphism
            $$\begin{tikzcd}[sep=small]
            {x F(y_0)} 
                \arrow[rrr,"b^x_{y_0}"] 
                \arrow[ddd,"1_{x} F(f)"']
            & {} 
            & {} 
            & {F(y_0) x}
                \arrow[ddd,"F(f) 1_{x}"]
                \arrow[dddlll,Rightarrow,"b^x_{f}"',shorten >=3ex, shorten <= 3ex] 
            \\ {} & {} & {} & {}
            \\ {} & {} & {} & {}
            \\ {x F(y_1)}
                \arrow[rrr,"b^x_{y_1}"']
            & {} 
            & {} 
            & {F(y_1) x}
            \end{tikzcd}$$ which is natural in 1-morphism $f:y_0 \to y_1$ in $\mathfrak{C}$;

            \item $R^x_{-,-}$ is an invertible modification 
            $$\begin{tikzcd}[sep=25pt]
            xF(yz) 
                \arrow[rr, "b^x_{yz}"]
            & {} \arrow[d, Rightarrow, "R^x_{y,z}"]          & F(yz)x 
            \\ {x F(y)F(z)}
                \arrow[r, "b^x_{y}1_{F(z)}"'] 
                \arrow[u,"1_xF_{y,z}"]
            & {F(y)xF(z) }
                \arrow[r, "1_{F(y)}b^x_{z}"'] 
            & {F(y)F(z)x}    
                \arrow[u,"F_{y,z}1_x"']
            \end{tikzcd}$$ satisfying the coherence condition
            \end{itemize}
            
            \settoheight{\calculus}{\includegraphics[width=45mm]{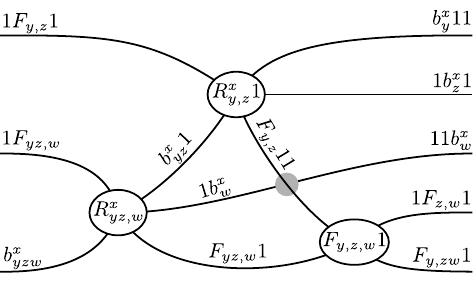}}

            \begin{equation}\label{eqn:DrinfeldCenterObject}
            \begin{tabular}{@{}ccc@{}}
            
            \includegraphics[width=45mm]{Pictures/E1centralizer/DrinfeldCenterObjectLeft.pdf} & \raisebox{0.45\calculus}{$=$} &
            \includegraphics[width=45mm]{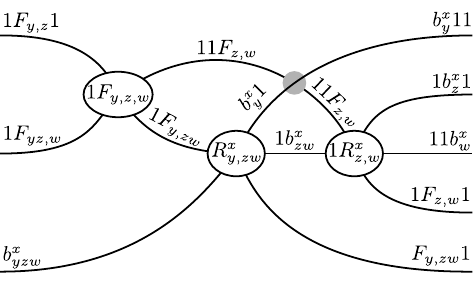},
            
            \end{tabular}
            \end{equation} 
            
        \begin{itemize}
            \item [] which is natural for objects $y,z,w$ in $\mathfrak{C}$.
        \end{itemize}

        \item [b.] A 1-morphisms between objects $(x,b^x_{-},R^x_{-,-})$ and $(y,b^y_{-},R^y_{-,-})$ is a pair $(g,b^g_{-})$ where
        \begin{itemize}
            \item $g:x \to y$ is a 1-morphism in $\mathfrak{D}$;

            \item $b^g_{-}$ is an invertible modification
            $$\begin{tikzcd}[sep=small]
            {x F(z)} 
                \arrow[rrr,"b^x_{z}"] 
                \arrow[ddd,"g 1_{F(z)}"']
            & {} 
            & {} 
            & {F(z) x}
                \arrow[ddd,"1_{F(z)} g"]
                \arrow[dddlll,Rightarrow,"b^g_{z}"',shorten >=3ex, shorten <= 3ex] 
            \\ {} & {} & {} & {}
            \\ {} & {} & {} & {}
            \\ {y F(z)}
                \arrow[rrr,"b^y_{z}"']
            & {} 
            & {} 
            & {F(z) y}
            \end{tikzcd}$$ which is natural for object $z$ in $\mathfrak{C}$, satisfying the coherence condition
            \end{itemize}
            \settoheight{\calculus}{\includegraphics[width=50mm]{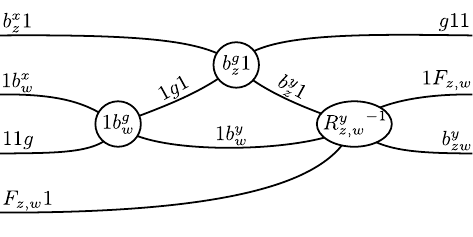}}

            \begin{equation}\label{eqn:DrinfeldCenter1Morphism}
            \begin{tabular}{@{}ccc@{}}
            
            \includegraphics[width=50mm]{Pictures/E1centralizer/DrinfeldCenter1MorphismLeft.pdf} & \raisebox{0.45\calculus}{$=$} &
            \includegraphics[width=37.5mm]{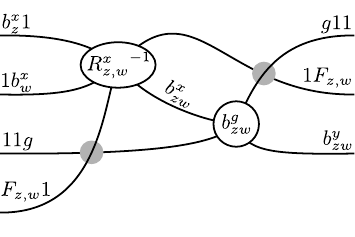},
            
            \end{tabular}
            \end{equation} 
            \begin{itemize}
                \item [] which is natural for objects $z,w$ in $\mathfrak{C}$.
            \end{itemize}
        
        \item [c.] A 2-morphism between $(g,b^g_{-})$ and $(h,b^h_{-})$, which are both 1-morphisms between $(x,b^x_{-},R^x_{-,-})$ and $(y,b^y_{-},R^y_{-,-})$, consists of 2-morphism $\varphi:g \to h$ in $\mathfrak{D}$ satisfying the coherence condition
            \settoheight{\calculus}{\includegraphics[width=45mm]{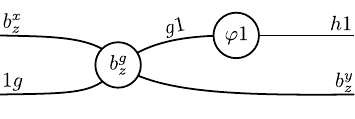}}

            \begin{equation}\label{eqn:DrinfeldCenter2Morphism}
            \begin{tabular}{@{}ccc@{}}
            
            \includegraphics[width=45mm]{Pictures/E1centralizer/DrinfeldCenter2MorphismLeft.pdf} & \raisebox{0.45\calculus}{$=$} &
            \includegraphics[width=45mm]{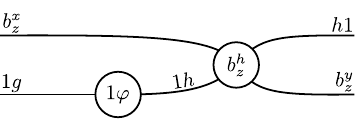},
            
            \end{tabular}
            \end{equation} 
            which is natural for object $z$ in $\mathfrak{C}$.
            
        \item [d.] Monoidal unit is given by the monoidal unit $I$ in $\mathfrak{D}$ together with $b^I_{-}$ given by the identity pseudo-natural isomorphism on $F$, and $R^I_{-,-}$ given by the identity modification.

        \item [e.] Monoidal product of $(x,b^x_{-},R^x_{-,-})$ and $(y,b^y_{-},R^y_{-,-})$ is given by the monoidal product $x \, \Box \, y$ in $\mathfrak{D}$ together with $b^{xy}_{-} :=(b^x_{-} 1_y) \circ (1_x b^y_{-})$ and invertible modification $R^{xy}_{-,-}$ given by $$\begin{tikzcd}[column sep=35pt]
            {xyF(zw)}
                \arrow[rr, "b^{xy}_{zw}"]
                \arrow[rd, "1 b^y_{zw}"] 
            & {} %\arrow[d, Rightarrow, "R_{x,y,z}"]          
            & {F(zw)xy}
            \\ {xyF(z)F(w)}
                \arrow[dd,"1 b^y_{z}1"']
                \arrow[u,"11F_{z,w}"]
            & {xF(zw)y}
                \arrow[ru, "b^x_{zw}1"] 
            & {F(z)F(w)xy}  
                \arrow[u,"F_{z,w}11"']
            \\ {}
                \arrow[r, Rightarrow, "1 R^y_{z,w}",shorten >= 20pt,shorten <= 20pt]
            & {xF(z)F(w)y} 
                \arrow[rd,"b^x_{z}11"']
                \arrow[r, Rightarrow, "R^x_{z,w} 1",shorten >= 20pt,shorten <= 20pt]
                \arrow[u,"1F_{z,w}1"]
            & {}
            \\ {xF(z)yF(w)}
                \arrow[ru,"11b^y_{w}"']
            & {}
            & {F(z)xF(w)y}
                \arrow[uu,"1b^x_{w}1"']
            \end{tikzcd}.$$

        \item [f.] Monoidal product of 1-morphism $(g,b^g_{-})$, which is from $(x_0,b^{x_0}_{-},R^{x_0}_{-,-})$ to $(y_0,b^{y_0}_{-},R^{y_0}_{-,-})$, and 1-morphism $(h,b^h_{-})$, which is from $(x_1,b^{x_1}_{-},R^{x_1}_{-,-})$ to $(y_1,b^{y_1}_{-},R^{y_1}_{-,-})$, consists of 1-morphism $g\, \Box \, h:x_0 \, \Box \, x_1 \to y_0 \, \Box \, y_1$ in $\mathfrak{D}$ together with an invertible modification
        \settoheight{\calculus}{\includegraphics[width=60mm]{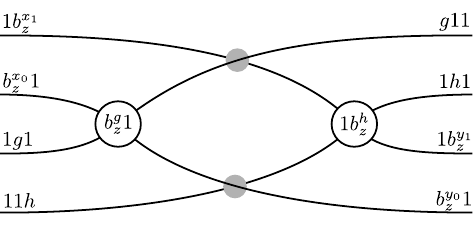}}
        \begin{equation}\label{eqn:DrinfeldCenterProduct1Morphism}
        \begin{tabular}{@{}cc@{}}
           \raisebox{0.45\calculus}{$b^{gh}_{z}$ :=} & \includegraphics[width=60mm]{Pictures/E1centralizer/DrinfeldCenterProduct1Morphism.pdf}.
        \end{tabular}
        \end{equation} 
        % $$\begin{tikzcd}
        %     {F(x_0x_1)z}
        %         \arrow[d,equal]
        %         \arrow[r,"F(g1)1"]
        %     & {F(y_0x_1)z}
        %         \arrow[d,equal]
        %         \arrow[r,"F(1h)1"]
        %     & {F(y_0y_1)z}
        %         \arrow[d,equal]
        %     \\ {F(x_0)F(x_1)z}
        %         \arrow[r,"F(g)11"]
        %     & {F(y_0)F(x_1)z}
        %         \arrow[r,"1F(h)1"]
        %     & {F(y_0)F(y_1)z}
        %     \\ {F(x_0)zF(x_1)}
        %         \arrow[r,"F(g)11"]
        %     & {F(y_0)zF(x_1)}
        %         \arrow[r,"11F(h)"]
        %     & {F(y_0)zF(y_1)}
        %     \\ {zF(x_0)F(x_1)}
        %         \arrow[d,equal]
        %         \arrow[r,"1F(g)1"]
        %     & {zF(y_0)F(x_1)}
        %         \arrow[d,equal]
        %         \arrow[r,"11F(h)"]
        %     & {zF(y_0)F(y_1)}
        %         \arrow[d,equal]
        %     \\ {z F(x_0x_1)}
        %         \arrow[r,"1F(g1)"']
        %     & {zF(y_0x_1)}
        %         \arrow[r,"1F(1h)"']
        %     & {z F(y_0y_1)}
        %     %\arrow[r,"\beta_{x_0x_1,z}"]
        %     %\arrow[r,"\beta_{y_0x_1,z}"]
        %     %\arrow[r,"\beta_{y_0y_1,z}"']
        % \end{tikzcd}$$

        \item [g.] Interchangers of 1-morphsms are induced by those in $\mathfrak{D}$.
        
        \item [h.] Monoidal product of 2-morphisms is just the monoidal product of the corresponding 2-morphisms in $\mathfrak{D}$.
    \end{itemize}
\end{Definition}

\begin{Lemma}
    By the above construction, there is a canonical forgetful functor $\mathscr{Z}_1(F) \to \mathfrak{D}$, preserving the monoidal structures.
\end{Lemma}

\begin{Proposition}
    Suppose $\mathfrak{C}$ is a monoidal 2-category and $\mathfrak{M}$ is a left $\mathfrak{C}$-module 2-category, then the Drinfeld centralizer of the left action $\mathfrak{C} \to \mathbf{End}(\mathfrak{M})$ is equivalent to $\mathbf{End}_{\mathfrak{C}}(\mathfrak{M})$.
\end{Proposition}

\begin{proof}
    The proposition can be easily proven by unfolding the definitions of the Drinfeld centralizer and the 2-category of left $\mathfrak{C}$-module functors.
\end{proof}

\begin{Definition}
    Given monoidal 2-functor $F:\mathfrak{C} \to \mathfrak{D}$, and monoidal 2-functor $G:\mathfrak{E}^{1mp} \to \mathfrak{D}$, we can construct a $(\mathfrak{C},\mathfrak{E})$-bimodule 2-category whose underlying 2-category is $\mathfrak{D}$, with left $\mathfrak{C}$-action given by $$ \mathfrak{C} \to  \mathbf{End}(\mathfrak{D}); \quad x \mapsto F(x) \, \Box \, -, $$ and right $\mathfrak{E}$-action given by $$ \mathfrak{E}^{1mp} \to \mathbf{End}(\mathfrak{D}); \quad y \mapsto - \, \Box \, G(y). $$ We denote this bimodule 2-category as ${}_{\langle F \rangle}\mathfrak{D}_{\langle G \rangle}$.
\end{Definition}

\begin{Corollary} \label{cor:DrinfeldCenterAsBimoduleFunctors}
    $\mathscr{Z}_1(F) \simeq \mathbf{End}_{\mathfrak{C} |\mathfrak{D}}({}_{\langle F \rangle}\mathfrak{D})$.
\end{Corollary}

\begin{proof}
    Notice that $\mathfrak{D} \simeq \mathbf{End}_{\mathfrak{D}^{1mp}}(\mathfrak{D}) \simeq \mathscr{Z}_1(\mathfrak{D}^{1mp} \to \mathbf{End}(\mathfrak{D}))$ is the centralizer of the right $\mathfrak{D}$-action on ${}_{\langle F \rangle}\mathfrak{D}$. Hence, we can rewrite the centralizer of $F:\mathfrak{C} \to \mathfrak{D}$ as the centralizer of the left $\mathfrak{C}$-action on ${}_{\langle F \rangle}\mathfrak{D}$ preserving the right $\mathfrak{D}$-action, i.e. $\mathscr{Z}_1(\mathfrak{C} \xrightarrow{F} \mathbf{End}_{\mathfrak{D}^{1mp}}(\mathfrak{D})) \simeq \mathscr{Z}_1(\mathfrak{C} \times \mathfrak{D}^{1mp} \xrightarrow{F \, \Box \, \mathrm{Id}_\mathfrak{D}} \mathbf{End}(\mathfrak{D})) \simeq \mathbf{End}_{\mathfrak{C} |\mathfrak{D}}({}_{\langle F \rangle}\mathfrak{D})$.
\end{proof}

\begin{Definition} \label{def:DrinfeldCenter}
    The braided center (or Drinfeld center, $\mathbb{E}_1$ center), see \cite[Section 2.1]{D9}, of a monoidal 2-category $\mathfrak{C}$ is defined to be the monoidal centralizer $\mathscr{Z}_1(\mathrm{Id}_\mathfrak{C})$. We denote it by $\mathscr{Z}_1(\mathfrak{C})$.
\end{Definition}

\begin{Lemma}
    Drinfeld center $\mathscr{Z}_1(\mathfrak{C})$ is equipped with a canonical braided monoidal structure.
\end{Lemma}

\begin{proof}
    
We define the braiding on objects $(x, b^x_{-}, R^x_{-,-})$ and $(y, b^y_{-}, R^y_{-,-})$ in $\mathscr{Z}_1(\mathfrak{C})$ as $b_{x,y} := b^x_{y}$. This braiding can be extended to an adjoint 2-natural equivalence. Additionally, for $(x, b^x_{-}, R^x_{-,-})$, $(y, b^y_{-}, R^y_{-,-})$ and $(z, b^z_{-}, R^z_{-,-})$ in $\mathscr{Z}_1(\mathfrak{C})$, we introduce invertible modifications $R_{x,y,z}:= R^x_{y,z}$ and $S^x_{y,z}:=1_{xyz}$. These assignments establish a braiding on Drinfeld center $\mathscr{Z}_1(\mathfrak{C})$ according to \cite[Theorem 3.1]{Cr}.
\end{proof}

\begin{Remark} \label{rmk:Z_1(C)=ΩEnd(∑C)}
    One can also realize the Drinfeld center as $\mathbf{Fun}_{\mathfrak{C} \boxtimes \mathfrak{C}^{1mp}}(\mathfrak{C},\mathfrak{C})$, see \cite[Lemma 2.2.1]{D9}. Under some finiteness conditions on $\mathfrak{C}$, the 3-category $\mathbf{Bimod}(\mathfrak{C})$ of $(\mathfrak{C},\mathfrak{C})$-bimodule 2-categories becomes a monoidal 3-category under the relative Deligne 2-tensor product $\boxtimes_\mathfrak{C}$, see \cite{D10}. $\mathfrak{C}$ is the monoidal unit in this monoidal 3-category. Thus $\mathscr{Z}_1(\mathfrak{C}) \simeq \mathbf{Fun}_{\mathfrak{C} \boxtimes \mathfrak{C}^{1mp}}(\mathfrak{C},\mathfrak{C})$ is the endo-hom 2-category on the monoidal unit, hence it is equipped with a canonical braided monoidal 2-category structure. Finally, these two braided monoidal 2-category structures on $\mathscr{Z}_1(\mathfrak{C})$ will agree with each other.\footnote{Depending on the exact definitions of Drinfeld center and relative 2-tensor product, these two braided monoidal 2-category structure might agree up to flipping the direction of the braidings.}
\end{Remark}

\begin{Remark} \label{rmk:CentralActionOfDrinfeldCenterOnMoritaDual}
    Drinfeld center is a Morita invariant, i.e. when $\mathfrak{C}$ is a monoidal 2-category, $\mathfrak{M}$ is a left $\mathfrak{C}$-module 2-category, we have $\mathscr{Z}_1(\mathfrak{C}) \simeq \mathscr{Z}_1(\mathbf{End}_\mathfrak{C}(\mathfrak{M})^{1mp})$ as braided monoidal 2-categories. Hence, it induces monoidal 2-functor $\mathscr{Z}_1(\mathfrak{C}) \to \mathbf{End}_\mathfrak{C}(\mathfrak{M})^{1mp}$ factoring through its Drinfeld center, i.e. a central $\mathscr{Z}_1(\mathfrak{C})$-action on $\mathbf{End}_\mathfrak{C}(\mathfrak{M})^{1mp}$. 
    
    More explicitly, following the notations from \cite{D4}, denote the left $\mathfrak{C}$-action on $\mathfrak{M}$ by $(\Diamond,\alpha^\mathfrak{M},l^\mathfrak{M},\mu^\mathfrak{M},\lambda^\mathfrak{M},\pi^\mathfrak{M})$, then for each object $(x,b^x_{-},R^x_{-,-})$ in $\mathscr{Z}_1(\mathfrak{C})$, and object $y$ in $\mathfrak{C}$, the interchanger $\chi^{x \, \Diamond \, -}_{y \, \Diamond \, -,-}$ is given by $$ x \, \Diamond \, (y \, \Diamond \, -) \xrightarrow{(\alpha^\mathfrak{M}_{x,y,-})^{-1}} (x \, \Box \, y) \, \Diamond \, - \xrightarrow{b^x_{y} \, \Diamond \, -} (y \, \Box \, x) \, \Diamond \, - \xrightarrow{\alpha^\mathfrak{M}_{y,x,-}} y \, \Diamond \, (x \, \Diamond \, -).$$ Higher coherence data $(\omega^{x \, \Diamond \, -},\gamma^{x \, \Diamond \, -})$ are induced by $R^x_{-,-}$ and $\mu^\mathfrak{M},\lambda^\mathfrak{M},\pi^\mathfrak{M}$.
\end{Remark}

\begin{Remark}
    A braiding on monoidal 2-category $\mathfrak{C}$ is equivalent to a monoidal section of the canonical forgetful 2-functor $\mathscr{Z}_1(\mathfrak{C}) \to \mathfrak{C}$, i.e. an embedding $\mathfrak{C} \to \mathscr{Z}_1(\mathfrak{C})$ preserving monoidal structures such that its composition with the forgetful functor is equivalent to the identity functor $\mathrm{Id}_\mathfrak{C}$. 
    
    Moreover, a braiding on $\mathfrak{C}$ induces a monoidal equivalence between $\mathfrak{C}$ and $\mathfrak{C}^{1mp}$, which can either be proven explicitly, or can be implied by observing that $\mathscr{Z}_1(\mathfrak{C})^{2mp} \simeq \mathscr{Z}_1(\mathfrak{C}^{1mp})$ as braided monoidal 2-categories.
\end{Remark}

\subsection{Factorization of Endofunctors into Free Modules and Monoidal Local Modules}

\begin{Proposition} \label{prop:MonoidalCentralizerIsMultifusion}
    For any tensor 2-functor $F:\mathfrak{C} \to \mathfrak{D}$ between multifusion 2-categories $\mathfrak{C}$ and $\mathfrak{D}$, its monoidal centralizer is multifusion.
\end{Proposition}

\begin{proof}
    This is basically the same as the proof for \cite[Corollary 2.2.2, Lemma 2.2.4]{D9}. First, to see $\mathscr{Z}_1(F)$ is finite semisimple, we combine Corollary \ref{cor:DrinfeldCenterAsBimoduleFunctors} with \cite[Proposition 5.1.3, Theorem 5.2.7]{D4} and \cite[Theorem 3.1.6]{D7}.

    Next, we need to show that $\mathscr{Z}_1(F)$ is rigid. For any object $(x,b^x_{-},R^x_{-,-})$ in the centralizer $\mathscr{Z}_1(F)$, we will construct its left dual. The right dual can be construct in a similar way.

    Since $\mathfrak{D}$ is rigid, we can find a left dual ${}^\lor x$ for object $x$, together with unit $i_x:I \to x \, \Box \, {}^\lor x$ and counit $e_x:{}^\lor x \, \Box \, x \to I$ and 2-isomorphisms witnessing zigzag conditions:
    $$\Xi_x:(1_x \, \Box \, e_x) \circ (i_x \, \Box \, 1_x) \to 1_x, $$
    $$ \Phi_x:(e_x \, \Box \, 1_{{}^\lor x}) \circ (1_{{}^\lor x} \, \Box \, i_x) \to 1_{{}^\lor x}.$$
    Furthermore, we can assume this left dual is coherent \cite{Pstr}, i.e. $\Xi_x$ and $\Phi_x$ satisfy the coherence conditions (\ref{eqn:coherentleftdual1}) and (\ref{eqn:coherentleftdual2}).

    We promote the object ${}^\lor x$ to an object in the centralizer $\mathscr{Z}_1(F)$ by
    $$ b^{{}^\lor x}_{y} := (e_x \, \Box \, 1_{F(y)} \, \Box \, 1_{{}^\lor x}) \circ (1_{{}^\lor x} \, \Box \, {b^x_{y}}^\bullet \, \Box \, 1_{{}^\lor x}) \circ (1_{{}^\lor x} \, \Box \, 1_{F(y)} \, \Box \, i_x), $$
    \settoheight{\calculus}{\includegraphics[width=60mm]{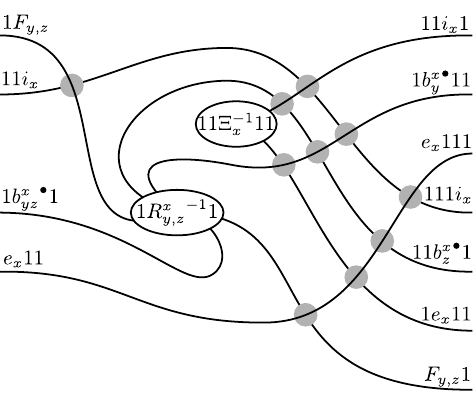}}
        \begin{center}
        \begin{tabular}{@{}cc@{}}
           \raisebox{0.45\calculus}{$R^{{}^\lor x}_{y,z} :=$} & \includegraphics[width=60mm]{Pictures/E1centralizer/DrinfeldCentralizerCondensationDualObject.pdf},
        \end{tabular}
        \end{center} because we assume the inverse of half-braiding $b^x_{-}$ is coherent.

    Then we promote the unit $i_x:I \to x \, \Box \, {}^\lor x$ and counit $e_x:{}^\lor x \, \Box \, x \to I$ to 1-morphisms in $\mathscr{Z}_1(F)$ by
    \settoheight{\calculus}{\includegraphics[width=45mm]{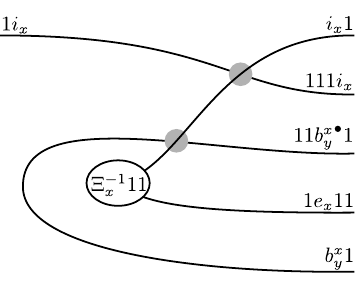}}
        \begin{center}
        \begin{tabular}{@{}cc@{}}
           \raisebox{0.45\calculus}{$b^{i_x}_{y} :=$} & \includegraphics[width=45mm]{Pictures/E1centralizer/DrinfeldCentralizerCondensationDualUnit.pdf},
        \end{tabular}
        \end{center}

    \settoheight{\calculus}{\includegraphics[width=45mm]{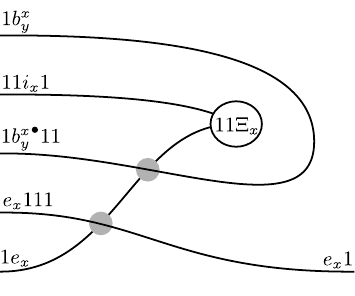}}
        \begin{center}
        \begin{tabular}{@{}cc@{}}
           \raisebox{0.45\calculus}{$b^{e_x}_{y} :=$} & \includegraphics[width=45mm]{Pictures/E1centralizer/DrinfeldCentralizerCondensationDualCounit.pdf}.
        \end{tabular}
        \end{center}
    Finally, we can check $\Xi_x$ and $\Phi_x$ satisfy condition (\ref{eqn:DrinfeldCenter2Morphism}), so they can be promoted to 2-morphisms in $\mathscr{Z}_1(F)$. In summary, we have constructed a left dual $({}^\lor x, b^{{}^\lor x}_{-},R^{{}^\lor x}_{-,-})$ for object $(x,b^x_{-},R^x_{-,-})$ in the Drinfeld centralizer $\mathscr{Z}_1(F)$.
\end{proof}

\begin{Remark}
    Direct sum of two objects $(x,b^x_{-},R^x_{-,-})$ and $(y,b^y_{-},R^y_{-,-})$ is given by $(x \boxplus y,b^x_{-} \boxplus b^y_{-}, R^x_{-,-} \boxplus R^y_{-,-})$.
\end{Remark}

\begin{Remark}
    We can explicitly demonstrate that $\mathscr{Z}_1(F)$ is 2-condensation complete. Given an object $(x,b^x_{-},R^x_{-,-})$ in $\mathscr{Z}_1(F)$, any 2-condensation monad $(x,b^x_{-},R^x_{-,-};e,b^e_{-};\xi,\delta)$ induces a 2-condensation monad $(x,e,\xi,\delta)$ in $\mathfrak{D}$. Since $\mathfrak{D}$ is assumed to be 2-condensation complete, this 2-condensation monad has a splitting witnessed by a 2-condensation $(x,y,f,g,\phi,\gamma)$ in $\mathfrak{D}$. Furthermore, by \cite[Proposition 3.1.5]{GJF}, we can assume this 2-condensation is unital without loss of generality, i.e. $g$ is the right adjoint of $f$ with unit $\eta: 1_x \to g \circ f$ and counit $\phi: f \circ g \to 1_y$.

    We can promote object $y$ to an object $(y,b^y_{-},R^y_{-,-})$ in $\mathscr{Z}_1(F)$, where the half-braiding is defined as follows:
        $$ b^y_{z} := {yF(z)} \xrightarrow{g 1} {xF(z)} \xrightarrow{b^x_{z}} {F(z)x} \xrightarrow{1 f} F(z)y.$$
        \settoheight{\calculus}{\includegraphics[width=60mm]{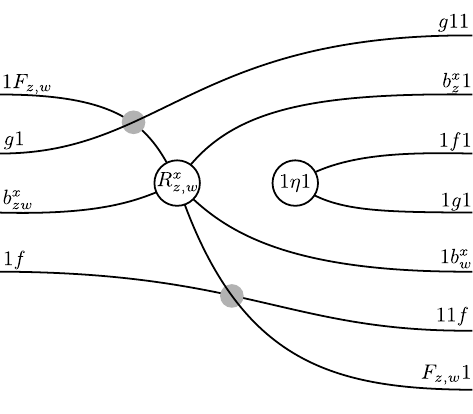}}
        \begin{center}
        \begin{tabular}{@{}cc@{}}
           \raisebox{0.45\calculus}{$R^y_{z,w} :=$} & \includegraphics[width=60mm]{Pictures/E1centralizer/DrinfeldCentralizerCondensationComplete.pdf}.
        \end{tabular}
        \end{center}

        Next, we can promote 1-morphisms $f:X \to Y$ and $g: Y \to X$ to 1-morphisms in $\mathscr{Z}_1(F)$, where
        \settoheight{\calculus}{\includegraphics[width=30mm]{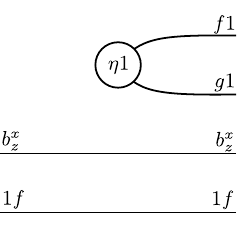}}
        \begin{center}
        \begin{tabular}{@{}cc@{}}
           \raisebox{0.45\calculus}{$b^f_{z} :=$} & \includegraphics[width=30mm]{Pictures/E1centralizer/DrinfeldCentralizerCondensationf.pdf},
        \end{tabular}
        \end{center}

        \settoheight{\calculus}{\includegraphics[width=30mm]{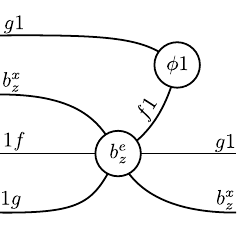}}
        \begin{center}
        \begin{tabular}{@{}cc@{}}
           \raisebox{0.45\calculus}{$b^g_{z} :=$} & \includegraphics[width=30mm]{Pictures/E1centralizer/DrinfeldCentralizerCondensationg.pdf}.
        \end{tabular}
        \end{center}
        
        Lastly, $\phi,\gamma$ both satisfy condition (\ref{eqn:DrinfeldCenter2Morphism}), so the 2-condensation $(x,y,f,g,\phi,\gamma)$ can be upgraded to a 2-condensation $(x,b^x_{-},R^x_{-,-};y,b^y_{-},R^y_{-,-};f,b^f_{-};g,b^g_{-};\phi,\gamma)$ splitting the given 2-condensation monad.
\end{Remark}

\begin{Remark} We can explicitly construct adjoints of 1-morphisms in $\mathscr{Z}_1(F)$. For any 1-morphism $f$ from object $(x,b^x_{-},R^x_{-,-})$ to $(y,b^y_{-},R^y_{-,-})$, we can promote the left adjoint ${}^*f: y \to x$ of the underlying 1-morphism in $\mathfrak{D}$ to $\mathscr{Z}_1(F)$ with
        \settoheight{\calculus}{\includegraphics[width=30mm]{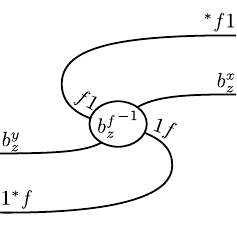}}
        \begin{center}
        \begin{tabular}{@{}cc@{}}
           \raisebox{0.45\calculus}{$b^{{}^*f}_{z} :=$} & \includegraphics[width=30mm]{Pictures/E1centralizer/DrinfeldCentralizerCondensationLeftAdjoint.pdf}.
        \end{tabular}
        \end{center}
        Similarly, we can also promote the right adjoint $f^*:y \to x$ to a 1-morphism in $\mathscr{Z}_1(F)$ with
        \settoheight{\calculus}{\includegraphics[width=30mm]{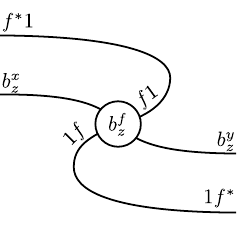}}
        \begin{center}
        \begin{tabular}{@{}cc@{}}
           \raisebox{0.45\calculus}{$(b^{{f}^*}_{z})^{-1} :=$} & \includegraphics[width=30mm]{Pictures/E1centralizer/DrinfeldCentralizerCondensationRightAdjoint.pdf}.
        \end{tabular}
        \end{center}
\end{Remark}

\begin{Proposition} \label{prop:braidedcenterpreservessumandproduct}
    Let $\mathfrak{C}$ and $\mathfrak{D}$ be multifusion 2-categories, then we have the follows:
    \begin{enumerate}
        \item $\mathscr{Z}_1(\mathfrak{C} \boxplus \mathfrak{D}) \simeq \mathscr{Z}_1(\mathfrak{C}) \boxplus \mathscr{Z}_1(\mathfrak{D})$,

        \item $\mathscr{Z}_1(\mathfrak{C} \boxtimes \mathfrak{D}) \simeq \mathscr{Z}_1(\mathfrak{C}) \boxtimes \mathscr{Z}_1(\mathfrak{D})$.
    \end{enumerate}
\end{Proposition}

\begin{proof}
    The first equation is by direct inspection. The second equation follows from \cite[Corollary 5.3.2]{D9}.
\end{proof}

\begin{Definition} \label{def:NonDegenerateFusion2Category}
    A multifusion 2-category $\mathfrak{C}$ is non-degenerate if $\mathscr{Z}_1(\mathfrak{C}) \simeq \mathbf{2Vect}$.
\end{Definition}

\begin{Lemma} \label{lem:SendingFreeModulesToEndohomOfModules}
    Let $A$ be a separable algebra in a multifusion 2-category $\mathfrak{C}$. Then we have a tensor 2-functor \[ \mathfrak{C}^{1mp} \to \mathbf{End}(\mathbf{Mod}_\mathfrak{C}(A)), \] which sends an object $x$ in $\mathfrak{C}$ to the endo-2-functor on $\mathbf{Mod}_\mathfrak{C}(A)$ determined by $A \mapsto x \, \Box \, A$.
\end{Lemma}

\begin{Theorem} \label{thm:E1LocalModulesAsE1Centralizer}
    Let $A$ be a separable algebra in a multifusion 2-category $\mathfrak{C}$. Then we have an equivalence of multifusion 2-categories $$ \mathbf{Mod}^{\mathbb{E}_1}_\mathfrak{C}(A) \simeq \mathscr{Z}_1(\mathfrak{C}^{1mp} \to \mathbf{End}(\mathbf{Mod}_\mathfrak{C}(A)) ). $$
\end{Theorem}

\begin{proof}
    This is essentially a reformulation of \cite[Theorem 5.3.2]{D8}.
\end{proof}

\begin{Theorem} \label{thm:DrinfeldCenterAndE1LocalModules}
    Let $A$ be a separable algebra in an indecomposable multifusion 2-category $\mathfrak{C}$. Then we have an equivalence of braided fusion 2-categories $$ \mathscr{Z}_1(\mathfrak{C}) \simeq \mathscr{Z}_1(\mathfrak{C}^{1mp} \boxtimes \mathbf{Mod}^{\mathbb{E}_1}_\mathfrak{C}(A) \to \mathbf{End}(\mathbf{Mod}_\mathfrak{C}(A)) ). $$
\end{Theorem}

\begin{proof}
    See \cite[Proposition 2.3.1]{D9}.
\end{proof}

\begin{Corollary} \label{cor:E1LocalModulesInNondegenerateFusion2Cateory}
    Suppose $\mathfrak{C}$ is a non-degenerate multifusion 2-category and $A$ is a separable algebra in $\mathfrak{C}$, then we have an equivalence of multifusion 2-categories $$ \mathbf{End}(\mathbf{Mod}_\mathfrak{C}(A)) \simeq \mathfrak{C}^{1mp} \boxtimes \mathbf{Mod}^{\mathbb{E}_1}_\mathfrak{C}(A). $$
\end{Corollary}

\begin{Corollary} \label{cor:FactorizationOfNondegenerateFusion2Cateory}
    In particular, for any non-degenerate multifusion 2-category $\mathfrak{C}$, we have $ \mathbf{End}(\mathfrak{C}) \simeq \mathfrak{C}^{1mp} \boxtimes \mathfrak{C}. $
\end{Corollary}

\begin{Corollary} \label{cor:E0CenterOfModuleCategoryIsE1LocalModuleInE0Center}
    Any separable algebra $A$ in non-degenerate multifusion 2-category $\mathfrak{C}$ can be viewed as a separable algebra in $\mathbf{End}(\mathfrak{C})$, with \[\mathbf{Mod}^{\mathbb{E}_1}_{\mathbf{End}(\mathfrak{C})}(A) \simeq \mathbf{End}(\mathbf{Mod}_\mathfrak{C}(A)).\]
\end{Corollary}

\begin{Remark} \label{rmk:nondegeneratemultifusion2categoryisindecomposable}
    A non-degenerate multifusion 2-category $\mathfrak{C}$ is indecomposable. Moreover, one has $\mathfrak{C} \simeq \boxplus_{i,j} \, \mathfrak{C}_{ij}$ where indices $i,j$ run across the set of isomorphism classes of simple objects $\{e_i\}$ in $\Omega \mathfrak{C} := Hom_\mathfrak{C}(I,I)$, and $\mathfrak{C}_{ij} := e_i \mathfrak{C} e_j$. Then component $\mathfrak{C}_{ii}$ is a non-degenerate fusion 2-category, while $\mathfrak{C}_{ij}$ is an invertible bimodule 2-category between $\mathfrak{C}_{ii}$ and $\mathfrak{C}_{jj}$.
\end{Remark}

\begin{Remark}
    By \cite[Lemma 4.1.3]{D9}, a non-degenerate fusion 2-category is always equivalent to $\mathbf{Mod}(\mathcal{B})$ for some non-degenerate braided fusion 1-category $\mathcal{B}$, i.e. a braided fusion 1-category $\mathcal{B}$ with a trivial Müger center: $\mathcal{Z}_2(\mathcal{B}) \simeq \mathbf{Vect}$.

    By \cite[Theorem 3.1.4]{D9}, non-degenerate fusion 2-categories $\mathbf{Mod}(\mathcal{B}_0)$ and $\mathbf{Mod}(\mathcal{B}_1)$ are Morita equivalent if and only if two non-degenerate braided fusion 1-categories $\mathcal{B}_0$ and $\mathcal{B}_1$ are Witt equivalent in the sense of \cite[Definition 5.1]{DMNO}.
\end{Remark}

\begin{Remark} \label{rmk:E1centerofE0centeristrivial}
    The 2-category of endofunctors on a finite semisimple 2-category $\mathfrak{M}$ is always non-degenerate, i.e. $\mathscr{Z}_1(\mathbf{End}(\mathfrak{M})) \simeq \mathbf{2Vect}$. Conversely, a non-degenerate fusion 2-category $\mathbf{Mod}(\mathcal{B})$ is equivalent to $\mathbf{End}(\mathfrak{M})$ for some finite semisimple 2-category $\mathfrak{M}$ if and only if $\mathcal{{B}}$ is Witt trivial, i.e. $\mathcal{B} \simeq \mathcal{Z}_1(\mathcal{M})$ for some fusion 1-category $\mathcal{M}$, in which case we can take $\mathfrak{M} \simeq \mathbf{Mod}(\mathcal{M})$.
\end{Remark}

\section{Braided Centralizers and Sylleptic Centers}

In this section, we will continue using notations and conventions defined in the preliminary.

\subsection{Definitions} \label{def:BraidedCentralizer}

We first generalize the notion of 2-center of braided monoidal 2-categories \cite[Section 5.1]{Cr}.
    \begin{Definition}
    Braided centralizer (or $\mathbb{E}_2$ centralizer) of a braided 2-functor $F:\mathfrak{A} \to \mathfrak{B}$ between braided monoidal 2-categories $\mathfrak{A}, \mathfrak{B}$ is a braided monoidal 2-category $\mathscr{Z}_{2}(F)$ where:
    
    \begin{enumerate}
        \item [a.] An object is a pair $(x,\sigma^x_{-})$ where $x$ is an object in $\mathfrak{B}$, $\sigma^x_{-}$ is an invertible modification $$\begin{tikzcd}[sep=small]
        {x F(y)} \arrow[rr,equal] \arrow[dddr,"b_{x,F(y)}"'] & {} \arrow[ddd,Rightarrow,"\sigma^x_{y}",shorten >=2ex, shorten <= 2ex] & {x F(y)}
        \\ {} & {} & {} & {}
        \\ {} & {} & {} & {}
        \\ {} & {F(y) x} \arrow[uuur,"b_{F(y),x}"'] & {} & {}
        \end{tikzcd}$$
    \end{enumerate}
    
    \begin{enumerate}
    \item[] given on object $y$ in $\mathfrak{A}$, satisfying:
    \end{enumerate}
    
    \settoheight{\diagramwidth}{\includegraphics[width=60mm]{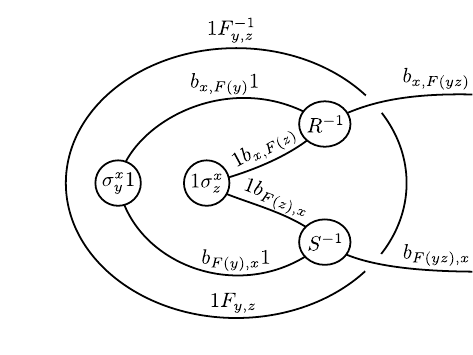}}

    \begin{equation}\label{eqn:braidedcentralizerobject}
    \begin{tabular}{@{}cccc@{}}
    \raisebox{0.15\diagramwidth}{\includegraphics[width=30mm]{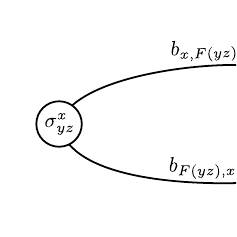}}&
    \raisebox{0.45\diagramwidth}{$=$} &
    \includegraphics[width=60mm]{Pictures/E2centralizer/syllepseleft.pdf} & \raisebox{0.45\diagramwidth}{.}
    \end{tabular}
    \end{equation}
    
    \begin{enumerate}
    \item [b.] A 1-morphism between objects $(x,\sigma^x_{-})$ and $(x',\sigma^{x'}_{-})$ is a 1-morphism $f:x \to x'$ in $\mathfrak{B}$ such that
    \end{enumerate}
    
    \settoheight{\diagramwidth}{\includegraphics[height=35mm]{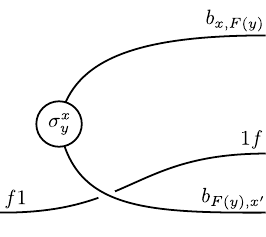}}

    \begin{equation}\label{eqn:braidedcentralizer1morphism}
    \begin{tabular}{@{}cccc@{}}
    \includegraphics[height=35mm]{Pictures/E2centralizer/morsyllepseleft.pdf}&
    \raisebox{0.45\diagramwidth}{$=$} &
    \includegraphics[height=35mm]{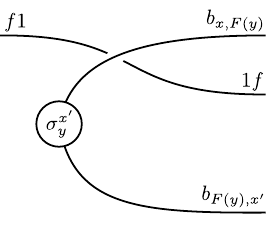} & \raisebox{0.45\diagramwidth}{.}
    \end{tabular}
    \end{equation}
    
    \begin{enumerate}
    \item [c.] A 2-morphism between 1-morphisms $f,g:x \to x'$ is just a 2-morphism $f \to g$ in $\mathfrak{B}$.
    \end{enumerate}

    \begin{enumerate}
    \item [d.] Monoidal product of two objects $(x,\sigma^x_{-})$ and $(y,\sigma^y_{-})$ consists of an object $x \, \Box \, y$ and an invertible modification
    \end{enumerate}

    \settoheight{\diagramwidth}{\includegraphics[width=60mm]{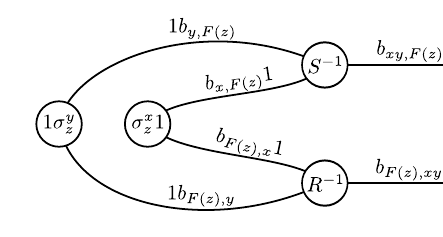}}

    \begin{equation}\label{eqn:braidedcentralizerproduct}
    \begin{tabular}{@{}cccc@{}}
    \raisebox{0.45\diagramwidth}{$\sigma^{xy}_{z}$} &
    \raisebox{0.45\diagramwidth}{$:=$} &
    \includegraphics[width=60mm]{Pictures/E2centralizer/syllepseprod.pdf} & \raisebox{0.45\diagramwidth}{.}
    \end{tabular}
    \end{equation}

    \begin{enumerate}
    \item [e.] Monoidal unit in $\mathscr{Z}_2(F)$ is given by the monoidal unit $I$ in $\mathfrak{B}$ with the identity modification $\sigma^I_{z} := 1_{F(z)}$.
    \end{enumerate}
    
    \begin{enumerate}
    \item [f.] Monoidal product of two 1-morphisms $f:x \to y$ and $g:x' \to y'$ is just 1-morphism $f \, \Box \, g: x \, \Box \, x' \to y \, \Box \, y'$ in $\mathfrak{B}$, which always turns out to satisfy equation (\ref{eqn:braidedcentralizer1morphism}).
    \end{enumerate}

    \begin{enumerate}
    \item [g.] Interchangers of 1-morphisms are induced from those in $\mathfrak{B}$.
    \end{enumerate}

    \begin{enumerate}
    \item [h.] Monoidal product of two 2-morphisms is just the monoidal product of the underlying 2-morphisms in $\mathfrak{B}$.
    \end{enumerate}

    \begin{enumerate}
    \item [i.] Braiding is induced from the braiding $b_{x,y}:x \, \Box \, y \to y \, \Box \, x$ in $\mathfrak{B}$, together with the invertible modifications $R$ and $S$.
    \end{enumerate}
    
    \end{Definition}
    
    \begin{Lemma}
    By definition, there is a canonical braiding-preserving forgetful functor $\mathscr{Z}_{2}(F) \to \mathfrak{B}$.
    \end{Lemma}

\subsection{Enveloping Algebras}
    \begin{Lemma} \label{lem:CentralActionInducedByBraidedFunctor}
    Braided 2-functor $F:\mathfrak{A} \to \mathfrak{B}$ induces a monoidal structure on the $\mathfrak{A} \boxtimes \mathfrak{A}^{1mp}$-module 2-category ${}_{\langle F \rangle} \mathfrak{B}_{\langle F \rangle}$, i.e. braided monoidal 2-category $\mathfrak{A} \boxtimes \mathfrak{A}^{2mp}$ acts centrally on monoidal 2-category ${}_{\langle F \rangle} \mathfrak{B}_{\langle F \rangle}$.
    \end{Lemma}

    \begin{proof}
    The left $\mathfrak{A} \boxtimes \mathfrak{A}^{1mp}$-action on $\mathfrak{B}$ is defined via $$\mathfrak{A} \boxtimes \mathfrak{A}^{1mp} \xrightarrow{F \, \boxtimes \, F} \mathfrak{B} \boxtimes \mathfrak{B}^{1mp} \xrightarrow{x \, \boxtimes \, y \, \longmapsto \, x \, \Box \, - \, \Box \, y} \mathbf{End}(\mathfrak{B}),$$ with monoidal 2-functor structure induced from that on $F$.
    
    In comparison, the central left $\mathfrak{A} \boxtimes \mathfrak{A}^{2mp}$-action on $\mathfrak{B}$ is defined via $$ \mathfrak{A} \boxtimes \mathfrak{A}^{2mp} \xrightarrow{F \, \boxtimes \, F} \mathfrak{B} \boxtimes \mathfrak{B}^{2mp} \xrightarrow{x \, \boxtimes \, y \, \longmapsto \, x \, \Box \, - \, \Box \, y} \mathscr{Z}_1(\mathfrak{B}),$$ with braided monoidal 2-functor structure induced from that on $F$. Note that the embeddings $\mathfrak{B} \to \mathscr{Z}_1(\mathfrak{B}); x \mapsto (x,b_{x,-},R_{x,-,-})$ and $\mathfrak{B}^{2mp} \to \mathscr{Z}_1(\mathfrak{B}); y \mapsto (y,b^\bullet_{-,y},S^\bullet_{-,-,y})$ preserves braidings. Using Remark \ref{rmk:Z_1(C)=ΩEnd(∑C)}, we identify $\mathscr{Z}_1(\mathfrak{B})$ with $\mathbf{End}_{\mathfrak{B} \boxtimes \mathfrak{B}^{1mp}}(\mathfrak{B})$, hence the two embeddings combine together as $\mathfrak{B} \boxtimes \mathfrak{B}^{2mp} \to \mathscr{Z}_1(\mathfrak{B})$ where object $x \boxtimes y$ in $\mathfrak{B} \boxtimes \mathfrak{B}^{2mp}$ is send to the 2-endofunctor $x \, \Box - \Box \, y: \mathfrak{B} \to \mathfrak{B}$ with $(\mathfrak{B},\mathfrak{B})$-bimodule 2-functor structure
    $$ a \, \Box \, x \, \Box - \Box \, y \, \Box \, b \xrightarrow{b_{a,x} \, \Box \, - \, \Box \, b^\bullet_{b,y}} x \, \Box \, a \, \Box - \Box \, b \, \Box \, y$$ plus the higher coherence data induced from the braiding on $\mathfrak{B}$.
    \end{proof}

\begin{Remark}
    $F:\mathfrak{A} \to {}_{\langle F \rangle} \mathfrak{B}_{\langle F \rangle}$ has a canonical $\mathfrak{A} \boxtimes \mathfrak{A}^{2mp}$-module 2-functor structure.
\end{Remark}

    \begin{Lemma}  \label{lem:EnvelopingAlgebraForBraided2Functor}
    $\int_{S^1}{F} := \mathbf{Fun}_{\mathfrak{A} \boxtimes \mathfrak{A}^{2mp}}({}_{\langle F \rangle} \mathfrak{B}_{\langle F \rangle},\mathfrak{A})$ is a monoidal 2-category, provided that $\mathfrak{A}$ has 2-coends, its monoidal product preserves colimits and the $\mathfrak{A} \boxtimes \mathfrak{A}^{2mp}$-action on $\mathfrak{B}$ admits an internal hom $[-,-]_{\mathfrak{A} \boxtimes \mathfrak{A}^{2mp}}: \mathfrak{B}^{1op} \boxtimes \mathfrak{B} \to \mathfrak{A} \boxtimes \mathfrak{A}^{2mp}$. In particular, these conditions are satisfied if both $\mathfrak{A}$ and $\mathfrak{B}$ are multifusion.
    \end{Lemma}

    \begin{proof}
    The central $\mathfrak{A} \boxtimes \mathfrak{A}^{2mp}$-module 2-category structure on ${}_{\langle F \rangle} \mathfrak{B}_{\langle F \rangle}$ is described in Lemma \ref{lem:CentralActionInducedByBraidedFunctor}; on the other hand, $\mathfrak{A}$ is equipped with the canonical central $\mathfrak{A} \boxtimes \mathfrak{A}^{2mp}$-action. Hence, $\int_{S^1}{F}$ is well-defined as a 2-category.

    To define the monoidal structure on $\int_{S^1}{F}$, we will use various notions of enriched 2-categories from Décoppet. Readers unfamiliar with the notations can consult \cite[Section 4]{D4} for details. We will also return to the discussion of enriched 2-categories later in Section 4.
    
    Given two objects $S$ and $T$ in $\int_{S^1}{F}$, viewed as $\mathfrak{A} \boxtimes \mathfrak{A}^{2mp}$-module functors from ${}_{\langle F \rangle} \mathfrak{B}_{\langle F \rangle}$ to $\mathfrak{A}$, their product is the Day convolution $$ (S \boxdot T)(c) := \int^{(a,b):\mathfrak{B} \times \mathfrak{B}} [a \, \Box \, b,c]_{\mathfrak{A} \boxtimes \mathfrak{A}^{2mp}} \, \Diamond \, S(a) \, \Box \, T(b),  $$ together with the canonical $\mathfrak{A} \boxtimes \mathfrak{A}^{2mp}$-module functor structure \vspace{10pt} \newline $ (S \boxdot T)((x \boxtimes z) \, \Diamond \, y) = (S \boxdot T)(F(x) \, \Box \, y \, \Box \, F(z))  $
        \begin{align*}
           := \, & \int^{(a,b):\mathfrak{B} \times \mathfrak{B}} [a \, \Box \, b,F(x) \, \Box \, y \, \Box \, F(z)]_{\mathfrak{A} \boxtimes \mathfrak{A}^{2mp}} \, \Diamond \, S(a) \, \Box \, T(b) \\
           \simeq \, & \, \int^{(a,b):\mathfrak{B} \times \mathfrak{B}} (x \boxtimes z) \, \Diamond \, [a \, \Box \, b,y]_{\mathfrak{A} \boxtimes \mathfrak{A}^{2mp}} \, \Diamond \, S(a) \, \Box \, T(b) \\
           \simeq \, & \, (x \boxtimes z) \, \Diamond \, \int^{(a,b):\mathfrak{B} \times \mathfrak{B}}  [a \, \Box \, b,y]_{\mathfrak{A} \boxtimes \mathfrak{A}^{2mp}} \, \Diamond \, S(a) \, \Box \, T(b) \\
           =: \, & (x \boxtimes z) \, \Diamond \, (S \boxdot T)(y) = s \, \Box \, (S \boxdot T)(y) \, \Box \, z, 
        \end{align*} where the first equivalence follows from \cite[Lemma 4.1.6]{D4}.

    The monoidal unit in $\int_{S^1}F$ is given by $I_{\int_{S^1}F} := [I_{\mathfrak{B}},-]_{\mathfrak{A} \boxtimes \mathfrak{A}^{2mp}} \, \Diamond \,I_\mathfrak{A}$. Functoriality of the monoidal product follows immediately from the definition of 2-coends.

    Lastly, higher coherence data is provided by the 2-universal property of 2-coends.
    \end{proof}

\begin{Remark}
    Consider the monoidal 2-functor $$ \mathfrak{A} \boxtimes \mathfrak{B}^{1mp} \to \int_{S^1}F; \quad x \boxtimes y \mapsto [{y}^\lor ,-]_{\mathfrak{A} \boxtimes \mathfrak{A}^{2mp}} \, \Diamond \, x. $$ It has restrictions to $\mathfrak{A}$ and $\mathfrak{B}^{1mp}$:
    $$ \mathfrak{A} \to \int_{S^1}F; \quad x \boxtimes I_\mathfrak{B} \mapsto [I_\mathfrak{B} ,-]_{\mathfrak{A} \boxtimes \mathfrak{A}^{2mp}} \, \Diamond \, x, $$
    $$ \mathfrak{B}^{1mp} \to \int_{S^1}F; \quad I_\mathfrak{A} \boxtimes y \mapsto [y^\lor,-]_{\mathfrak{A} \boxtimes \mathfrak{A}^{2mp}} \, \Diamond \, I_\mathfrak{A}. $$
    
    Meanwhile, this 2-functor is dominant since it has a left adjoint $$ \int_{S^1}F \to \mathfrak{A} \boxtimes \mathfrak{B}^{1mp}; \quad S \mapsto \int^{y : \mathfrak{B}} S(y) \boxtimes {}^\lor y $$ with the unit of the adjunction given on $\int_{S^1}F$-module 2-functor $S:\mathfrak{B} \to \mathfrak{A}$ as an equivalence $$ S \to \int^{y:\mathfrak{B}} [y,-]_{\mathfrak{A} \boxtimes \mathfrak{A}^{2mp}} \, \Diamond \, S(y). $$
\end{Remark}

\begin{Remark}
    For object $x$ in $\mathfrak{A}$ and $y$ in $\mathfrak{B}$, the image of $x \boxtimes y$ in $\int_{S^1}F$,  $[y^\lor,-]_{\mathfrak{A} \boxtimes \mathfrak{A}^{2mp}} \, \Diamond \, x$, is equipped with a canonical automorphism in $\int_{S^1}F$: $$ [y^\lor,-]_{\mathfrak{A} \boxtimes \mathfrak{A}^{2mp}} \, \Diamond \, x \simeq [y^\lor,-]_{\mathfrak{A} \boxtimes \mathfrak{A}^{2mp}} \, \Diamond \, (x \boxtimes I_\mathfrak{A}) \, \Diamond \, I_\mathfrak{A} $$
    $$ \simeq [y^\lor,F(x) \, \Box \, -]_{\mathfrak{A} \boxtimes \mathfrak{A}^{2mp}} \, \Diamond \, I_\mathfrak{A} \xrightarrow{b_{F(x)},-} [y^\lor,- \, \Box \, F(x)]_{\mathfrak{A} \boxtimes \mathfrak{A}^{2mp}} \, \Diamond \, I_\mathfrak{A} $$
    $$ \simeq [y^\lor,-]_{\mathfrak{A} \boxtimes \mathfrak{A}^{2mp}} \, \Diamond \, (I_\mathfrak{A} \boxtimes x) \, \Diamond \, I_\mathfrak{A} \simeq [y^\lor,-]_{\mathfrak{A} \boxtimes \mathfrak{A}^{2mp}} \, \Diamond \, x.$$ Moreover, this induces an automorphism on an arbitrary object $S$ in $\int_{S^1}F$, which we will denote by $aut^{+}(S)$.
    
    Replacing braiding $b_{F(x),-}$ by the reversed braiding $b^\bullet_{-,F(x)}$ in $\mathfrak{B}$, we can obtain another canonical automorphism on an arbitrary object $S$ in $\int_{S^1}F$, denoted by $aut^{-}(S)$.
\end{Remark}

\begin{Remark}
    By Lemma \ref{lem:CentralActionInducedByBraidedFunctor}, we see that $F$ induces a central $\mathfrak{A} \boxtimes \mathfrak{A}^{2mp}$-action on $\mathfrak{B}^{1mp}$ via the braided functor $\mathfrak{A} \boxtimes \mathfrak{A}^{2mp} \to \mathscr{Z}_1(\mathfrak{B})$, $$ x_0 \boxtimes I_\mathfrak{A} \mapsto (F(x_0),b_{F(x_0),-},R_{F(x_0),-,-})$$
    $$ I_\mathfrak{A} \boxtimes x_1 \mapsto (F(x_1),b^\bullet_{-,F(x_1)},S^\bullet_{-,-,F(x_1)}).$$ Furthermore, the embedding $\mathfrak{B}^{1mp} \to \int_{S^{1}}F$ induces a canonical forgetful functor $\mathfrak{A} \boxtimes \mathfrak{A}^{2mp} \to \int_{S^{1}}F$ factoring through the Drinfeld center $\mathscr{Z}_1(\mathfrak{B})$, i.e.
    $$\begin{tikzcd}
        {\mathfrak{A} \boxtimes \mathfrak{A}^{2mp}}
            \arrow[d]
        & {\mathfrak{B}^{1mp}}
            \arrow[d]
        \\ {\mathscr{Z}_1(\mathfrak{B})}
            \arrow[ur]
            \arrow[r,dashed]
        & {\int_{S^{1}}F}.
    \end{tikzcd}$$
    Fix an object $x$ in $\mathfrak{A}$, we observe that $S = [I_\mathfrak{B},-]_{\mathfrak{A} \boxtimes \mathfrak{A}^{2mp}} \, \Diamond \, x$ is isomorphic to the image of $x \boxtimes I_\mathfrak{A}$ and $I_\mathfrak{A} \boxtimes x$ in $\int_{S^1}F$. Meanwhile, $aut^+(S)$ is just the image of the half-braiding $b_{F(x),-}$ while $aut^-(S)$ is the image of the half-braiding $b^\bullet_{-,F(x)}$.
\end{Remark}

\begin{Remark} \label{rmk:BalancedActionOnEnvelopingAlgebra}
    The central actions of $\mathfrak{A} \boxtimes \mathfrak{A}^{2mp}$ on $\mathfrak{A}$ and $\mathfrak{B}^{1mp}$ are balanced in the sense that there exists a monoidal 2-natural adjoint equivalence filling the following diagram:
    $$\begin{tikzcd}
        {\mathfrak{A} \boxtimes \mathfrak{A}^{2mp}}
            \arrow[rr,"F \boxtimes F"]
            \arrow[d]
        & {}
        & {\mathfrak{B} \boxtimes \mathfrak{B}^{2mp}}
            \arrow[d]
        \\ {\mathscr{Z}_1(\mathfrak{A})}
            \arrow[d]
        & {}
        & {\mathscr{Z}_1(\mathfrak{B})}
            \arrow[d]
        \\ {\mathfrak{A}}
            \arrow[r]
        & {\int_{S^1}F}
        & {\mathfrak{B}^{1mp}}
            \arrow[l]
    \end{tikzcd}$$ Indeed, starting with object $x \boxtimes I_\mathfrak{A} $ in $\mathfrak{A} \boxtimes \mathfrak{A}^{2mp}$ and following the two paths, we obtain two objects in $\int_{S^1}F$ which are canonically equivalent as follows $$[I_\mathfrak{B},-] \, \Diamond \, x \simeq [I_\mathfrak{B},F(x) \, \Box \, - ] \, \Diamond \, I \xrightarrow{b_{F(x),-}} [I_\mathfrak{B},- \, \Box \, F(x)] \, \Diamond \, I \simeq [F(x)^{\lor},-] \, \Diamond \, I_\mathfrak{A}. $$
    Similarly, starting with object $ I_\mathfrak{A} \boxtimes x$ in $\mathfrak{A} \boxtimes \mathfrak{A}^{2mp}$ and following the two paths, there is a canonical equivalence $$ [I_\mathfrak{B},-] \, \Diamond \, x \simeq [I_\mathfrak{B},F(x) \, \Box \, - ] \, \Diamond \, I \xrightarrow{b^\bullet_{-,F(x)}} [I_\mathfrak{B},- \, \Box \, F(x)] \, \Diamond \, I \simeq  [F(x)^{\lor},-] \, \Diamond \, I_\mathfrak{A}. $$

    Moreover, these two equivalent monoidal 2-functors $\mathfrak{A} \boxtimes \mathfrak{A}^{2mp} \to \int_{S^1}F$ can be lifted to the braided 2-functor $\mathfrak{A} \boxtimes \mathfrak{A}^{2mp} \to \mathscr{Z}_1(\int_{S^1}F)$ via the following half-braiding on the image of $x \boxtimes I_\mathfrak{A}$
    $$ \begin{tikzcd}[sep=25pt]
        {([I_\mathfrak{B},-] \, \Diamond \, x) \boxdot ([y^\lor,-] \, \Diamond \, z)}
            \arrow[r,"b_{F(x),-}"]
            \arrow[d,"\sim"']
         & {([F(x)^\lor,-] \, \Diamond \, I_\mathfrak{A}) \boxdot ([y^\lor,-] \, \Diamond \, z)}
            \arrow[d,"\sim"]
        \\ {[y^\lor,-] \, \Diamond \, (x \, \Box \, z) }
            \arrow[d,"b_{x,z}"']
        & {[(y \, \Box \, F(x))^\lor,-] \, \Diamond \, z }
            \arrow[d,"b_{F(x),y}^\lor"]
        \\ {[y^\lor,-] \, \Diamond \, (z \, \Box \, x) }
            \arrow[d,"\sim"']
        & {[(F(x) \, \Box \, y)^\lor,-] \, \Diamond \, z }
            \arrow[d,"\sim"]
        \\ {([y^\lor,-] \, \Diamond \, z) \boxdot ([I_\mathfrak{B},-] \, \Diamond \, x)}
            \arrow[r,"b_{F(x),-}"']
        & {([y^\lor,-] \, \Diamond \, z) \boxdot ([F(x)^\lor,-] \, \Diamond \, I_\mathfrak{A}).}
    \end{tikzcd} $$ Replace the braidings on $\mathfrak{A}$ and $\mathfrak{B}$ by their reverses, this gives a half-braiding on the image of $I_\mathfrak{A} \boxtimes x$.
\end{Remark}

\begin{Remark}
    Given objects $x_0 \boxtimes y_0$ and $x_1 \boxtimes y_1$ in $\mathfrak{A} \boxtimes \mathfrak{B}^{1mp}$, their images in $\int_{S^1}F$ are equivalent, i.e. $$ [y_0^\lor,-]_{\mathfrak{A} \boxtimes \mathfrak{A}^{2mp}} \, \Diamond \, x_0 \simeq [y_1^\lor,-]_{\mathfrak{A} \boxtimes \mathfrak{A}^{2mp}} \, \Diamond \, x_1, $$ if and only if $F(x_0) \, \Box \, y_0 \simeq F(x_1) \, \Box \, y_1$ as objects in $\mathfrak{B}$.
\end{Remark}

\begin{Remark}
    If braided monoidal 2-categories $\mathfrak{A}$ and $\mathfrak{B}$ are multifusion, then following \cite[Section 2]{D10}, we expect that the relative Deligne 2-tensor product $\boxtimes_{\mathfrak{A} \boxtimes \mathfrak{A}^{2mp}}$ exists, and it preserves monoidal structures. One can rewrite $\int_{S^1}{F} \simeq \mathfrak{A} \boxtimes_{\mathfrak{A} \boxtimes \mathfrak{A}^{2mp}} \mathfrak{B}^{1mp}$ as multifusion 2-categories, and then the canonical $\mathfrak{A} \boxtimes \mathfrak{A}^{2mp}$-balanced 2-functor $$\mathfrak{A} \boxtimes \mathfrak{B}^{1mp} \to \mathfrak{A} \boxtimes_{\mathfrak{A} \boxtimes \mathfrak{A}^{2mp}} \mathfrak{B}^{1mp}$$ is witnessed by the monoidal 2-functor $$\mathfrak{A} \boxtimes \mathfrak{B}^{1mp} \to \int_{S^1}F; \quad x \boxtimes y \mapsto [{y}^\lor ,-]_{\mathfrak{A} \boxtimes \mathfrak{A}^{2mp}} \, \Diamond \, x.$$
\end{Remark}

\subsection{Finite Semisimple Braided Centralizers}
    \begin{Lemma} \label{lem:BraidedCentralizerAsMoritaDual}
    The $(\mathfrak{A},\mathfrak{B})$-bimodule 2-category structure on ${}_{\langle F \rangle} \mathfrak{B}$ is promoted to an $\int_{S^1}{F}$-module 2-category structure. Then we have $\mathscr{Z}_2(F) \simeq \mathbf{End}_{\int_{S^1}{F}}(\mathfrak{B})$ as monoidal 2-categories.
    \end{Lemma}

    \begin{proof}
    The $\int_{S^1}{F}$-action on $\mathfrak{B}$ is defined via $ S \, \Diamond \, z := FS(z),$ for object $S$ in $\int_{S^1}{F}$ and $z$ in $\mathfrak{B}$.  In particular, the monoidal unit $I_{\int_{S^1}F} = [I_\mathfrak{B},-]_{\mathfrak{A} \boxtimes \mathfrak{A}^{2mp}} \, \Diamond \, I_\mathfrak{A}$ acts on $\mathfrak{B}$ via $$ I_{\int_{S^1}F} \, \Diamond \, z := F([I_\mathfrak{B},z]_{\mathfrak{A} \boxtimes \mathfrak{A}^{2mp}} \, \Diamond \, I_\mathfrak{A}) \simeq [I_\mathfrak{B},z]_{\mathfrak{A} \boxtimes \mathfrak{A}^{2mp}} \, \Diamond \, F(I_\mathfrak{A}) \simeq z,  $$ where the equivalences are induced by the monoidal 2-functor structure on $F$.

    The original $(\mathfrak{A},\mathfrak{B})$-action on ${}_{\langle F \rangle} \mathfrak{B}$ can be recovered from the monoidal 2-functor $$ \mathfrak{A} \boxtimes \mathfrak{B}^{1mp} \to \int_{S^1}F; \quad x \boxtimes y \mapsto [{y}^\lor ,-]_{\mathfrak{A} \boxtimes \mathfrak{A}^{2mp}} \, \Diamond \, x. $$ This induces a forgetful 2-functor $\mathbf{End}_{\int_{S^1}{F}}(\mathfrak{B}) \to \mathbf{End}_{\mathfrak{A}|\mathfrak{B}}({}_{\langle F \rangle} \mathfrak{B})$ preserving monoidal structures. By Corollary \ref{cor:DrinfeldCenterAsBimoduleFunctors}, $\mathbf{End}_{\mathfrak{A}|\mathfrak{B}}({}_{\langle F \rangle} \mathfrak{B})$ is equivalent to $\mathscr{Z}_1(F)$, the monoidal centralizer of $F$. Hence, every $\int_{S^1}F$-module 2-functor $T:\mathfrak{B} \to \mathfrak{B}$ can be forgotten to a triple $(\pmb{y},b^{\pmb{y}}_{-},R^{\pmb{y}}_{-,-})$ by the definition of monoidal centralizer.

    More explicitly, one has the interchanger $\chi^T$:
    $$ \begin{tikzcd}
        {S \, \Diamond \, T(z)}
            \arrow[r,"\chi^T_{S,z}"]
            \arrow[d,equal]
        & {T(S \, \Diamond \, z)}
            \arrow[d,equal]
        \\ {[y^\lor,\pmb{y} \, \Box \, z]_{\mathfrak{A} \boxtimes \mathfrak{A}^{2mp}} \, \Diamond \, F(x)}
            \arrow[d,"\sim"']
        & {\pmb{y} \, \Box \, [y^\lor,z]_{\mathfrak{A} \boxtimes \mathfrak{A}^{2mp}} \, \Diamond \, F(x)}
            \arrow[d,"\sim"]
        \\ {F(x) \, \Box \, \pmb{y} \, \Box \, z \, \Box \, y}
        & {\pmb{y} \, \Box \, F(x) \, \Box \, z \, \Box \, y}
            \arrow[l,"b^{\pmb{y}}_{F(x)}"]
    \end{tikzcd} $$ which is natural for object $S = [y^\lor,-]_{\mathfrak{A} \boxtimes \mathfrak{A}^{2mp}} \, \Diamond \, x $ in $\int_{S^1}F$ and object $z$ in $\mathfrak{B}$.

    For any object $\pmb{x}$ in $\mathfrak{A}$, we notice that under the central $\mathfrak{A} \boxtimes \mathfrak{A}^{2mp}$-action on $\int_{S^1}F$ (c.f. Remark \ref{rmk:BalancedActionOnEnvelopingAlgebra}), the image of $\pmb{x} \boxtimes I_\mathfrak{A}$ in $\int_{S^1}F$ is equipped with a half-braiding, hence by Remark \ref{rmk:CentralActionOfDrinfeldCenterOnMoritaDual}, it also acts centrally on $\mathbf{End}_{\int_{S^1}F}(\mathfrak{B})$. Forgetting to $\mathcal{Z}_1(F)$, we obtain a triple $(F(\pmb{x}),b_{F(\pmb{x}),-},R_{F(\pmb{x}),-,-})$. On the other hand, the image of $I_\mathfrak{A} \boxtimes \,\pmb{x}$ in $\int_{S^1}F$ is equivalent to that of $\pmb{x} \,\boxtimes I_\mathfrak{A}$, while its image in $\mathscr{Z}_1(F)$ is the triple $(F(\pmb{x}),b^\bullet_{-,F(\pmb{x})},S^\bullet_{-,-,F(\pmb{x})})$. Thus there must be an invertible modification between these two half-braidings $b_{F(\pmb{x}),\pmb{y}} \simeq b^\bullet_{\pmb{y},F(\pmb{x})}$. Now, if we fix $\pmb{y}$ and let $\pmb{x}$ vary in $\mathfrak{A}$, the data we obtain is equivalent to a half-syllepsis on $\pmb{y}$ by Currying, satisfying coherence condition (\ref{eqn:braidedcentralizerobject}).

    Conversely, given an object $(\pmb{y},\sigma^{\pmb{y}}_{-})$ in $\mathscr{Z}_2(F)$, one can realize it as a 2-functor $\pmb{y} \, \Box \, - : \mathfrak{B} \to \mathfrak{B}$ with $\int_{S^1}F$-module 2-functor structure:
    $$ \begin{tikzcd}
        {S \, \Diamond \, T(z)}
            \arrow[r,"\chi^T_{S,z}"]
            \arrow[d,equal]
        & {T(S \, \Diamond \, z)}
            \arrow[d,equal]
        \\ {[y^\lor,\pmb{y} \, \Box \, z]_{\mathfrak{A} \boxtimes \mathfrak{A}^{2mp}} \, \Diamond \, F(x)}
            \arrow[d,"\sim"']
        & {\pmb{y} \, \Box \, [y^\lor,z]_{\mathfrak{A} \boxtimes \mathfrak{A}^{2mp}} \, \Diamond \, F(x)}
            \arrow[d,"\sim"]
        \\ {F(x) \, \Box \, \pmb{y} \, \Box \, z \, \Box \, y}
            \arrow[r,"b^\bullet_{\pmb{y},F(x)}"']
        & {\pmb{y} \, \Box \, F(x) \, \Box \, z \, \Box \, y}
    \end{tikzcd} $$ which is natural for object $S = [y^\lor,-]_{\mathfrak{A} \boxtimes \mathfrak{A}^{2mp}} \, \Diamond \, x $ in $\int_{S^1}F$ and object $z$ in $\mathfrak{B}$. The half-syllepsis on $\pmb{y}$ witnesses the compatibility data for the central $\mathfrak{A} \boxtimes \mathfrak{A}^{2mp}$-action on $\int_{S^1}F$. Lastly, the correspondence between objects can be extended monoidally and functorially to produce an equivalence of monoidal 2-categories $\mathscr{Z}_2(F) \simeq \mathbf{End}_{\int_{S^1}F}(\mathfrak{B})$.
    \end{proof}

\begin{Remark}
    Fix object $(\pmb{y},b^{\pmb{y}}_{-},R^{\pmb{y}}_{-,-})$ in the monoidal centralizer $\mathscr{Z}_1(F)$. The fiber of $\mathbf{End}_{\int_{S^1}{F}}(\mathfrak{B}) \to \mathscr{Z}_1(F) $ over this object consists of all possible $\mathfrak{A} \boxtimes \mathfrak{A}^{2mp}$-balancing structures on the $(\mathfrak{A},\mathfrak{B})$-bimodule 2-functor $\pmb{y} \, \Box \, - :\mathfrak{B} \to \mathfrak{B}$. For a general treatment of balanced 2-functors and relative Deligne 2-tensor product, readers may consult \cite[Section 2]{D10}.
\end{Remark}
    
    \begin{Definition} \label{def:SyllepticCenter}
    Sylleptic center (or $\mathbb{E}_2$ center) of a braided monoidal 2-category $\mathfrak{B}$ is defined to be the braided centralizer $\mathscr{Z}_{2}(\mathrm{Id}_{\mathfrak{B}})$. We denote the sylleptic center of $\mathfrak{B}$ as $\mathscr{Z}_{2}(\mathfrak{B})$.
    \end{Definition}
    
    \begin{Lemma}
    Sylleptic center $\mathscr{Z}_{2}(\mathfrak{B})$ of a braided monoidal 2-category $\mathfrak{B}$ is a sylleptic monoidal 2-category.
    \end{Lemma}    

    \begin{proof}
        The syllepsis of $\mathscr{Z}_{2}(\mathfrak{B})$ is given by $\sigma^x_{y}$ for objects $(x,\sigma^x_{-})$ and $(y,\sigma^y_{-})$. Using functoriality of $\sigma^x_{-}$ and compatibility condition (\ref{eqn:braidedcentralizer1morphism}) for 1-morphisms, we can extend $(x,\sigma^x_{-};y,\sigma^y_{-}) \mapsto \sigma_{x,y}$ to an invertible modification. Finally, this syllepsis and braiding $(b,R,S)$ in $\mathfrak{B}$ satisfy the compatibility conditions (\ref{eqn:syllepsisaxiom1}), (\ref{eqn:syllepsisaxiom2}) and unitality, as checked in \cite[Theorem 5.1]{Cr}.
    \end{proof}

\begin{Remark} \label{rmk:Z_2(C)=ΩΩEnd(∑∑C)}
    By analogy with Remark \ref{rmk:Z_1(C)=ΩEnd(∑C)}, we can understand the sylleptic monoidal structure on $\mathscr{Z}_2(\mathfrak{B})$ with a point of view from higher Morita theory. Since any braided monoidal 2-category $\mathfrak{B}$ can be viewed as a braided algebra in the 3-category $\mathbf{2Cat}$, we expect to construct a 3-category $\mathbf{Mod}^{\mathbb{E}_2}(\mathfrak{B})$, generalizing the notion of braided module 1-category by Davydov and Nikshych \cite[Definition 4.1]{DN}, consisting of braided local module 2-categories over braided monoidal 2-category $\mathfrak{B}$. With some nice assumptions, we expect $\mathbf{Mod}^{\mathbb{E}_2}(\mathfrak{B})$ to equip with a braided monoidal 3-category structure, with relative Deligne 2-tensor product $\boxtimes_\mathfrak{B}$ \cite{D10} as its monoidal product, and $\mathfrak{B}$ as its monoidal unit. Then we should be able to recover the sylleptic center $\mathscr{Z}_2(\mathfrak{B})$ as the endo-hom of the monoidal unit in $\mathbf{Mod}^{\mathbb{E}_2}(\mathfrak{B})$, which is equipped with a canonical sylleptic monoidal 2-category structure.
\end{Remark}

    \begin{Remark}
    A syllepsis on a braided monoidal 2-category $\mathfrak{B}$ is equivalent to a braided section of the canonical forgetful 2-functor $\mathscr{Z}_2(\mathfrak{B}) \to \mathfrak{B}$, i.e. an embedding $\mathfrak{B} \to \mathscr{Z}_2(\mathfrak{B})$ preserving braidings such that its composition with the forgetful functor is equivalent to the identity functor $\mathrm{Id}_\mathfrak{B}$. 
    
    Moreover, a syllepsis on $\mathfrak{B}$ induces a braided equivalence between $\mathfrak{B}$ and $\mathfrak{B}^{2mp}$, which can either be proven explicitly, or can be implied by observing that $\mathscr{Z}_2(\mathfrak{B})^{3mp} \simeq \mathscr{Z}_2(\mathfrak{B}^{2mp})$ as sylleptic monoidal 2-categories.
    \end{Remark}
    
    \begin{Proposition} \label{prop:BraidedCentralizerIsBraidedMultifusion}
        When $F:\mathfrak{A} \to \mathfrak{B}$ is a braided 2-functor between braided multifusion 2-categories, its braided centralizer $\mathscr{Z}_{2}(F)$ is a braided multifusion 2-category.
    \end{Proposition}

    \begin{proof}
    First, $\mathscr{Z}_{2}(F)$ is finite semisimple. Using \cite[Proposition 5.1.3, Theorem 5.2.7]{D4} and \cite[Theorem 3.1.6]{D7}, we see that $\int_{S^1}F$ is a multifusion 2-category when both $\mathfrak{A}$ and $\mathfrak{B}$ are braided multifusion 2-categories. By Lemma \ref{lem:BraidedCentralizerAsMoritaDual}, the claim is justified by applying \cite[Proposition 5.1.3, Theorem 5.2.7]{D4} and \cite[Theorem 3.1.6]{D7} again to $\mathscr{Z}_2(F) \simeq \mathbf{End}_{\int_{S^1}{F}}(\mathfrak{B})$.
    
    We still need to show that braided centralizer $\mathscr{Z}_2(F)$ is rigid. For any object $(x,\sigma^x_{-})$ in the centralizer $\mathscr{Z}_2(F)$, we will construct its left dual. The right dual can be construct in a similar way.

    Since $\mathfrak{B}$ is rigid, we can find a left dual ${}^\lor x$ for object $x$, together with unit $i_x:I \to x \, \Box \, {}^\lor x$ and counit $e_x:{}^\lor x \, \Box \, x \to I$ and 2-isomorphisms witnessing zigzag conditions:
    $$\Xi_x:(1_x \, \Box \, e_x) \circ (i_x \, \Box \, 1_x) \to 1_x, $$
    $$ \Phi_x:(e_x \, \Box \, 1_{{}^\lor x}) \circ (1_{{}^\lor x} \, \Box \, i_x) \to 1_{{}^\lor x}.$$
    Furthermore, we can assume this left dual is coherent \cite{Pstr}, i.e. $\Xi_x$ and $\Phi_x$ satisfy the coherence conditions (\ref{eqn:coherentleftdual1}) and (\ref{eqn:coherentleftdual2}).

    We also assume that the braiding $b$ is an adjoint 2-natural equivalence, and denote its coherent inverse by $b^\bullet$. Without loss of generality, we can assume $$b_{{}^\lor x,y} = (e_x \, \Box \, 1_y \, \Box \, 1_{{}^\lor x}) \circ (1_{{}^\lor x} \, \Box \, b^\bullet_{x,y} \, \Box \, 1_{{}^\lor x}) \circ (1_{{}^\lor x} \, \Box \, 1_y \, \Box \, i_x),$$ for any object $x$ and $y$ in $\mathfrak{B}$

     We promote the object ${}^\lor x$ to an object in the centralizer $\mathscr{Z}_2(F)$ with the half-syllepsis
     \settoheight{\calculus}{\includegraphics[width=60mm]{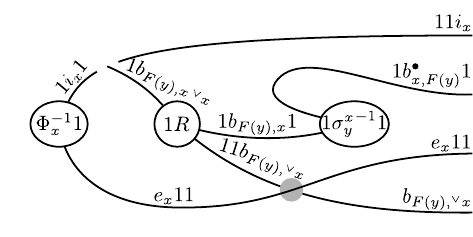}}
        \begin{center}
        \begin{tabular}{@{}ccc@{}}
           \raisebox{0.45\calculus}{$\sigma^{{}^\lor x}_{y} :=$} & \includegraphics[width=60mm]{Pictures/E2centralizer/syllepseDual.pdf} & \raisebox{0.45\calculus}{.}
        \end{tabular}
        \end{center}

    Next, we can promote the 1-morphisms $i_x$ and $e_x$ to 1-morphisms in $\mathscr{Z}_2(F)$ since they satisfy condition (\ref{eqn:braidedcentralizer1morphism}). Finally, we can promote the 2-morphisms $\Phi_x$ and $\Xi_x$ to 2-morphisms in $\mathscr{Z}_2(F)$ for free. In summary, we have lifted the left dual of object $x$ in $\mathfrak{B}$ to the left dual of the object $(x,\sigma^x_{-})$ in the braided centralizer $\mathscr{Z}_2(F)$.
    \end{proof}

    \begin{Remark}
        Direct sum of two objects $(x,\sigma^x_{-})$ and $(y,\sigma^y_{-})$ in $\mathscr{Z}_2(F)$ is given by $(x \boxplus y, \sigma^x_{-} \boxplus \sigma^y_{-})$.
    \end{Remark}

    \begin{Remark}
        $\mathscr{Z}_{2}(F)$ is 2-condensation complete. Given an object $(x,\sigma^x_{-})$ in $\mathscr{Z}_{2}(F)$, any 2-condensation monad $(x,\sigma^x_{-};e,\xi,\delta)$ gives rise to a 2-condensation monad $(x,e,\xi,\delta)$ in $\mathfrak{B}$, hence it splits into a 2-condensation $(x,y,f,g,\phi,\gamma)$ in $\mathfrak{B}$. Given $e$ is a 1-morphism in $\mathfrak{B}$ satisfying condition (\ref{eqn:braidedcentralizer1morphism}), we can then see that $f,g$ are also 1-morphisms in $\mathfrak{B}$ satisfying condition (\ref{eqn:braidedcentralizer1morphism}). 
        
        Thus, $(x,\sigma^x_{-};y,\sigma^y_{-};f,g,\phi,\gamma)$ provides a splitting of 2-condensation monad $(x,\sigma^x_{-};e,\xi,\delta)$ in $\mathscr{Z}_{2}(F)$, where the half-syllepsis on $y$ is defined via
        $$\begin{tikzcd}[sep=30pt]
            {}
            & {}
                \arrow[d,Rightarrow,shorten <=10pt,shorten >=5pt,"\gamma 1" {xshift=2pt,yshift=-5pt}]
            & {}
        \\    {y F(z)}
                \arrow[r,"g 1"]
                \arrow[dd,equal,bend right=40pt]
                \arrow[rr,equal,bend left=50pt]
            & {x F(z)}
                \arrow[r,"f 1"]
                \arrow[d,"b_{x,F(z)}"]
                \arrow[dd,equal,bend right=80pt]
            & {y F(z)}
                \arrow[d,"b_{y,F(z)}"]
        \\    {}
                \arrow[r,Rightarrow,shorten <=25pt,shorten >=5pt,"\sigma^x_{z}" {xshift=10pt,yshift=2pt}]
            & {F(z) x}
                \arrow[r,"1 f"]
                \arrow[d,"b_{F(z),x}"]
            & {F(z) y}
                \arrow[d,"b_{F(z),y}"]
        \\    {y F(z)}
                \arrow[r,"g 1"']
                \arrow[rr,equal,bend right=50pt]
            & {x F(z)}
                \arrow[r,"f 1"']
                \arrow[d,Rightarrow,shorten <=5pt,shorten >=10pt,"\phi 1" {xshift=2pt,yshift=5pt}]
            & {y F(z)}
        \\  {}
            & {}
            & {}
        \end{tikzcd}.$$
    \end{Remark}

    \begin{Remark}
        For any 1-morphism $f$ in $\mathscr{Z}_{2}(F)$ from object $(x,\sigma^x_{-})$ to $(y,\sigma^y_{-})$, the left adjoint ${}^*f$ and right adjoint $f^*$ for the underlying 1-morphism $f$ both exist in $\mathfrak{B}$. Then one can check that ${}^*f$ and $f^*$ both satisfy condition (\ref{eqn:braidedcentralizer1morphism}), hence they become left and right adjoints of $f$ in $\mathscr{Z}_{2}(F)$, respectively.
    \end{Remark}

\subsection{Braided Local Modules form a Braided Centralizer}

    We denote the enveloping algebra $\int_{S^1} \mathrm{Id}_\mathfrak{B}$ as $\int_{S^1} \mathfrak{B}$ for any braided fusion 2-category $\mathfrak{B}$.
    
    \begin{Lemma} \label{lem:TensorProductOfEnvelopingAlgebra}
    Given any two braided fusion 2-categories, one has an equivalence of monoidal 2-categories $\int_{S^1} \mathfrak{A} \boxtimes \int_{S^1} \mathfrak{B} \simeq \int_{S^1}(\mathfrak{A} \boxtimes \mathfrak{B}). $
    \end{Lemma}

    \begin{proof}
        The left hand side is $\mathbf{Fun}_{\mathfrak{A} \boxtimes \mathfrak{A}^{2mp}}(\mathfrak{A},\mathfrak{A}) \boxtimes \mathbf{Fun}_{\mathfrak{B} \boxtimes \mathfrak{B}^{2mp}}(\mathfrak{B},\mathfrak{B})$ while the right hand side is $\mathbf{Fun}_{\mathfrak{A} \boxtimes \mathfrak{B} \boxtimes \mathfrak{B}^{2mp} \boxtimes \mathfrak{A}^{2mp}}(\mathfrak{A} \boxtimes \mathfrak{B},\mathfrak{A} \boxtimes \mathfrak{B})$. By universal properties of Deligne 2-tensor product $\boxtimes$ from \cite{D3}, and matching the module actions carefully on both sides, it is straightforward to see that the underlying 2-categories are equivalent. Fix $\mathfrak{A} \boxtimes \mathfrak{A}^{2mp}$-module 2-functors $S_0,S_1:\mathfrak{A} \to \mathfrak{A}$ and $\mathfrak{B} \boxtimes \mathfrak{B}^{2mp}$-module 2-functors $T_0, T_1:\mathfrak{B} \to \mathfrak{B}$, if we view $S_0 \boxtimes T_0$ and $S_1 \boxtimes T_1$ as objects lying in the left hand side, then their product sends object $x \boxtimes y$ in $\mathfrak{A} \boxtimes \mathfrak{B}$ to $(S_0 \boxtimes T_0) \boxdot (S_1 \boxtimes T_1)(x \boxtimes y)$ \[  := \int^{a_0,a_1:\mathfrak{A}; b_0,b_1: \mathfrak{B}} [(a_0 \boxtimes b_0) \, \Box \, (a_1 \boxtimes b_1), x \boxtimes y]_{\mathfrak{A} \boxtimes \mathfrak{B} \boxtimes \mathfrak{B}^{2mp} \boxtimes \mathfrak{A}^{2mp}} \, \Diamond \, \] \[(S_0 \boxtimes T_0)(a_0 \boxtimes b_0) \, \Box \, (S_1 \boxtimes T_1)(a_1 \boxtimes b_1), \]
        \[= \int^{a_0,a_1:\mathfrak{A}; b_0,b_1: \mathfrak{B}} [(a_0 \, \Box \, a_1) \boxtimes (b_0 \, \Box \, b_1), x \boxtimes y]_{\mathfrak{A} \boxtimes \mathfrak{B} \boxtimes \mathfrak{B}^{2mp} \boxtimes \mathfrak{A}^{2mp}} \, \Diamond \, \] \[(S_0(a_0) \, \Box \, S_1(a_1)) \boxtimes  (T_0(b_0) \, \Box \, T_1(b_1)). \] while when viewed as objects in the right hand side, the product sends object $x$ in $\mathfrak{A}$ and object $y$ in $\mathfrak{B}$ to \[(S_0 \boxdot S_1)(x) \boxtimes (T_0 \boxdot T_1)(y) := \int^{a_0,a_1:\mathfrak{A}}[a_0 \, \Box \, a_1,x]_{\mathfrak{A} \boxtimes \mathfrak{A}^{2mp}} \, \Diamond \, S_0(a_0) \, \Box \, S_1(a_1)\]
        \[ \boxtimes \int^{b_0,b_1:\mathfrak{B}}[b_0 \, \Box \, b_1,y]_{\mathfrak{B} \boxtimes \mathfrak{B}^{2mp}} \, \Diamond \, T_0(b_0) \, \Box \, T_1(b_1).\]
        Again, the monoidal product is canonically equivalent due to the universal properties of Deligne 2-tensor product $\boxtimes$ from \cite{D3}. Notice that we can separate the internal homs into parts because in the $\mathfrak{A} \boxtimes \mathfrak{B} \boxtimes \mathfrak{B}^{2mp} \boxtimes \mathfrak{A}^{2mp}$-enriched structure on $\mathfrak{A} \boxtimes \mathfrak{B}$, essentially only the component $\mathfrak{A}$ is non-trivially enriched in $ \mathfrak{A} \boxtimes \mathfrak{A}^{2mp}$, and only the component $\mathfrak{B}$ is non-trivially enriched in $ \mathfrak{B} \boxtimes \mathfrak{B}^{2mp}$.
    \end{proof}
    
    \begin{Proposition} \label{prop:syllepticcenterpreservessumandproduct}
    Let $\mathfrak{A}$ and $\mathfrak{B}$ be braided multifusion 2-categories, then we have the follows:
    \begin{enumerate}
        \item $\mathscr{Z}_2(\mathfrak{A} \boxplus \mathfrak{B}) \simeq \mathscr{Z}_2(\mathfrak{A}) \boxplus \mathscr{Z}_2(\mathfrak{B})$,

        \item $\mathscr{Z}_2(\mathfrak{A} \boxtimes \mathfrak{B}) \simeq \mathscr{Z}_2(\mathfrak{A}) \boxtimes \mathscr{Z}_2(\mathfrak{B})$.
    \end{enumerate}
    \end{Proposition}

    \begin{proof}
    The first equivalence is by direct inspection. To prove the second equivalence, let's recall from Lemma \ref{lem:BraidedCentralizerAsMoritaDual} that
    \[\mathscr{Z}_2(\mathfrak{A}) \boxtimes \mathscr{Z}_2(\mathfrak{B}) \simeq \mathbf{End}_{\int_{S^1}\mathfrak{A}}(\mathfrak{A}) \boxtimes \mathbf{End}_{\int_{S^1}\mathfrak{B}}(\mathfrak{B}) \simeq \mathbf{End}_{\int_{S^1}\mathfrak{A} \boxtimes \int_{S^1}\mathfrak{B}}(\mathfrak{A} \boxtimes \mathfrak{B}).\]
    Moreover, by Lemma \ref{lem:TensorProductOfEnvelopingAlgebra} one can identify monoidal 2-categories $\int_{S^1}\mathfrak{A} \boxtimes \int_{S^1}\mathfrak{B}$ and $\int_{S^1}(\mathfrak{A} \boxtimes \mathfrak{B})$; also, it is immediate that this equivalence is compatible with their actions on the module 2-category $\mathfrak{A} \boxtimes \mathfrak{B}$. Hence, $\mathscr{Z}_2(\mathfrak{A} \boxtimes \mathfrak{B}) \simeq \mathbf{End}_{\int_{S^1}(\mathfrak{A} \boxtimes \mathfrak{B})}(\mathfrak{A} \boxtimes \mathfrak{B})$ and $\mathscr{Z}_2(\mathfrak{A}) \boxtimes \mathscr{Z}_2(\mathfrak{B}) $ have equivalent underlying 2-categories. Finally, we can explicitly write down the equivalence $\mathscr{Z}_2(\mathfrak{A}) \boxtimes \mathscr{Z}_2(\mathfrak{B}) \to \mathscr{Z}_2(\mathfrak{A} \boxtimes \mathfrak{B})$, \[(x,\sigma^x_{-}) \boxtimes (y,\sigma^y_{-}) \mapsto (x \boxtimes y,\sigma^{x \boxtimes y}_{-})\] where the half-syllepsis $\sigma^{x \boxtimes y}_{-}$ is determined by $\sigma^{x \boxtimes y}_{a \boxtimes b} := \sigma^x_{a} \boxtimes \sigma^y_{b}$ on object $a$ in $\mathfrak{A}$ and object $b$ in $\mathfrak{B}$. Then it is clear that this equivalence of underlying 2-categories also respect the braided monoidal structures on both sides.
    \end{proof}
    
    \begin{Definition} \label{def:NonDegenerateBraidedFusion2Category}
    A braided multifusion 2-category $\mathfrak{B}$ is non-degenerate if $\mathscr{Z}_2(\mathfrak{B})\simeq \mathbf{2Vect}$.
    \end{Definition}

    \begin{Lemma} \label{lem:SendingFreeModulesToDrinfeldCenterOfModules}
        Suppose $\mathfrak{B}$ is a braided multifusion 2-category with a separable braided algebra $B$, then there is a braided 2-functor $\mathfrak{B}^{2mp} \to \mathscr{Z}_1(\mathbf{Mod}_\mathfrak{B}(B))$.
    \end{Lemma}

    \begin{proof}
    We would like to assign a half-braiding $(x \, \Box \, B,\widehat{b}^{x \Box B}_{-},\widehat{R}^{x \Box B}_{-,-})$ for any object $x$ in $\mathfrak{B}^{2mp}$. By Corollary \ref{cor:DrinfeldCenterAsBimoduleFunctors}, this is equivalent to a $\mathbf{Mod}_\mathfrak{B}(B)$-bimodule 2-functor on $\mathbf{Mod}_\mathfrak{B}(B)$, namely,
    \begin{itemize}
        \item 2-functor $\mathbf{Mod}_\mathfrak{B}(B) \to \mathbf{Mod}_\mathfrak{B}(B); $
        \[(M,n^M,\nu^M,\rho^M) \mapsto (x \, \Box \, M,x \, \Box \, n^M,x \, \Box \,\nu^M,x \, \Box \,\rho^M). \]

        \item 2-natural equivalence $\varsigma^x_{M,N}:x \, \Box \, (M \, \Box_B \, N) \simeq (x \, \Box \, M) \, \Box_B \, N$ given on right $B$-modules $M$ and $N$, which is induced by naturality and associativity of monoidal product $\Box$.

        \item 2-natural equivalence $\vartheta^x_{M,N}:x \, \Box \, (M \, \Box_B \, N) \simeq M \, \Box_B \, (x \, \Box \, N) $ given on right $B$-modules $M$ and $N$, which is induced by the 1-morphism \[t^x_{M,N}:M \, \Box \, x \, \Box \, N \xrightarrow{b_{M,x} \, \Box \, N} x \, \Box \, M \, \Box \, N \xrightarrow{x \, \Box \, t_{M,N}} x \, \Box \, (M \, \Box_B \, N) \] with balancing structure
        \settoheight{\braid}{\includegraphics[width=50mm]{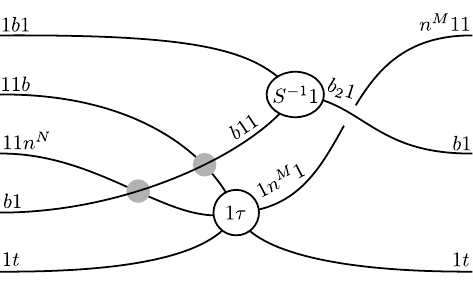}}
        
        $$\raisebox{0.45\braid}{$\tau^x_{M,N}:=\ $}
        \includegraphics[width=50mm]{Pictures/E2centralizer/tauxMN.pdf}.$$

        \item Bimodule 2-functor structure includes higher coherence data for 2-natural equivalences $\varsigma^x_{-,-}$ and $\vartheta^x_{-,-}$. They can also be explicitly constructed from the coherence data for braiding $b$ and monoidal structure $\Box$ on $\mathfrak{B}$; however, to avoid the unnecessarily verbose proof, we choose to omit it here and leave the details to keen readers.
    \end{itemize}
    In particular, viewed as half-braiding $\widehat{b}^{x B}_{M}$ it is given by \[(x \, \Box \, B ) \, \Box_B \, M \xrightarrow{(\varsigma^x_{B,M})^{-1}} x \, \Box \, (B \, \Box_B \, M) \xrightarrow{x \, \Box \, \pmb{l}^{-1}_M} x \, \Box \, M \] \[ \xrightarrow{x \, \Box \, \pmb{r}^{-1}_M} x \, \Box \, (M \, \Box_B \, B) \xrightarrow{\vartheta^x_{M,B}} M \, \Box_B \, (x \, \Box \, B). \]

    Next, using the above explicit construction, it is easy to extend the assignment to a 2-functor $\mathfrak{B}^{2mp} \to \mathscr{Z}_1(\mathbf{Mod}_\mathfrak{B}(B))$ via the naturality of $\Box$.

    Lastly, the braided 2-functor structure is induced as follows. Firstly, the monoidal 2-functor structure (see also \cite[Lemma 3.5]{DY}) is given by adjoint equivalence
    \[x B  \, \Box_B \, y B \xrightarrow{(\varsigma^x_{B,yB})^{-1}} x(B \, \Box_B \, yB) \xrightarrow{x \, \Box \, \pmb{l}^{-1}_{yB}} xyB, \]
    together with invertible modifications \[
    \begin{tikzcd}[column sep=30pt]
        {(xB  \, \Box_B \, y B) \, \Box_B \, z B}
            \arrow[r,"\pmb{\alpha}_{xB,yB,zB}"]
            \arrow[d,"(\varsigma^x)^{-1}\, \Box_B \, 1"']
        & {xB  \, \Box_B \, (y B \, \Box_B \, z B)}
            \arrow[d,"1 \, \Box_B \, (\varsigma^y)^{-1}"]
        \\ {(x(B \, \Box_B \, yB)) \, \Box_B \, z B}
            \arrow[d,"(1 \pmb{l}^{-1}) \, \Box_B \, 1"']
        & {xB  \, \Box_B \, (y(B \, \Box_B \, zB))}
            \arrow[d,"1 \, \Box_B \, (1 \pmb{l}^{-1})"]
        \\ {xyB \, \Box_B \, z B}
            \arrow[d,"(1 \pmb{r}^{-1}) \, \Box_B \, 1"']
            \arrow[r,Rightarrow,shorten >= 20pt,shorten <= 20pt]
        & {xB \, \Box_B \, yz B}
            \arrow[d,"1 \, \Box_B \, (1 \pmb{r}^{-1})"]
        \\ {xy(B \, \Box_B \, zB)}
            \arrow[d,"(\varsigma^{xy})^{-1}"']
        & {x(B \, \Box_B \, yzB)}
            \arrow[d,"(\varsigma^{x})^{-1}"]
        \\ {xyzB}
            \arrow[r,equal]
        & {xyzB,}
    \end{tikzcd}
    \] 
    \[\begin{tikzcd}
        {B \, \Box_B \, xB}
            \arrow[d,"(\varsigma^{I}_{B,xB})^{-1}"']
        & {xB}
            \arrow[d,equal]
            \arrow[l,"\pmb{l}_{xB}"']
        & {x(B \, \Box_B \, B)}
            \arrow[l,"x \, \Box \, \pmb{l}^{-1}_{B}"']
        \\ {B \, \Box_B \, xB}
            \arrow[r,"\pmb{l}^{-1}_{xB}"']
            \arrow[ur,Rightarrow,shorten >= 10pt,shorten <= 10pt]
        & {xB}
            \arrow[r,"\pmb{r}^{-1}_{xB}"']
        & {xB \, \Box_B \, B,}
            \arrow[u,"(\varsigma^{x}_{B,B})^{-1}"']
            \arrow[ul,Rightarrow,shorten >= 10pt,shorten <= 10pt]
    \end{tikzcd}\]
    induced by the associativity and unitality coherence data $\pmb{\alpha},\pmb{l},\pmb{r}$ on the monoidal 2-category $\mathbf{Mod}_\mathfrak{B}(B)$. 
    
    Finally, this 2-functor preserves the braiding with invertible modification
    \[\begin{tikzcd}
        {x B  \, \Box_B \, y B}
            \arrow[d,"(\varsigma^x_{B,yB})^{-1}"']
            \arrow[r,"(\varsigma^x_{B,yB})^{-1}"]
        & {x(B \, \Box_B \, yB)}
            \arrow[r,"x \, \Box \, \pmb{l}^{-1}_{yB}"]
        & {xyB}
            \arrow[dddd,"b^\bullet_{y,x} \, \Box \, B"]
        \\ {x(B \, \Box_B \, yB)}
            \arrow[d,"x \, \Box \, \pmb{l}^{-1}_{yB}"']
        & {}
        & {}
        \\ {xyB}
            \arrow[d,"x \, \Box \, \pmb{r}^{-1}_M"']
            \arrow[rr,Rightarrow,shorten >= 30pt,shorten <= 30pt]
        & {}
        & {}
        \\ {x(yB \, \Box_B \, B)}
            \arrow[d,"\vartheta^x_{yB,B}"']
        & {}
        & {}
        \\ {yB \, \Box_B \, xB}
            \arrow[r,"(\varsigma^y_{B,xB})^{-1}"]
        & {y(B \, \Box_B \, xB)}
            \arrow[r,"y \, \Box \, \pmb{l}^{-1}_{xB}"]
        & {yxB}
    \end{tikzcd}\] induced by the balancing $\tau^x_{yB,B}$ on the 1-morphism $t^x_{yB,B}:yBxB \xrightarrow{b_{yB,x}1} xyBB \xrightarrow{1t_{yB,B}} x(B \, \Box_B \, yB)$.
    \end{proof}
    
    \begin{Theorem} \label{thm:E2LocalModulesAsE2Centralizer}
    Let $B$ be a separable braided algebra in a braided multifusion 2-category $\mathfrak{B}$. Then we have an equivalence of braided multifusion 2-categories $$ \mathbf{Mod}_{\mathfrak{B}}^{\mathbb{E}_2}(B) \simeq  \mathscr{Z}_2(\mathfrak{B}^{2mp} \to \mathscr{Z}_1(\mathbf{Mod}_\mathfrak{B}(B)) ) . $$
    \end{Theorem}

    \begin{proof}
    By \cite[Remark 2.2.5]{DX}, there is a braided 2-functor \[\mathbf{Mod}_{\mathfrak{B}}^{\mathbb{E}_2}(B) \to \mathscr{Z}_1(\mathbf{Mod}_\mathfrak{B}(B))\] sending a fixed local $B$-module $(M,n^M,\nu^M,\rho^M,h^M)$ to the half-braiding $(M,\widetilde{b}_{M,-},\\ \widetilde{R}_{M,-,-})$. Then we can upgrade it to a half-syllepsis \[{\sigma}^M_{x}:1_{M \, \Box_B \, xB} \to \widehat{b}^{xB}_{M} \circ \widetilde{b}_{M,xB}\] given on object $x$ in $\mathfrak{B}$, via the 2-universal property of relative tensor product \[t_{M,xB}:MxB \to M \, \Box_B \, xB \] applied to the $B$-balanced 2-isomorphism:
    \[\begin{tikzcd}[row sep=20pt]
        {}
        & {M \, \Box_B \, xB}
            \arrow[dr,"\widetilde{b}_{M,xB}"]
        & {}
        \\ {MxB}
            \arrow[ur,"t_{M,xB}"]
            \arrow[r,"b_{M,xB}"]
            \arrow[dr,"b_{M,x}1"']
            \arrow[d,"t_{M,xB}"']
        & {xBM}
            \arrow[r,"t_{xB,M}"]
            \arrow[dr,"1t_{B,M}"']
        & {xB \, \Box_B \, M}
        \\ {M \, \Box_B \, xB}
        & {xMB}
            \arrow[d,"1t_{M,B}"]
            \arrow[u,"1b_{M,B}"]
        & {x \, \Box \, (B\, \Box_B \, M)}
            \arrow[u,"\varsigma^x_{B,M}"']
        \\{}
        & {x \, \Box \, (M \, \Box_B \, B)}
            \arrow[ul,"\vartheta^x_{M,B}"]
            \arrow[r,"1\pmb{r}_M"']
            \arrow[ur,"1\widetilde{b}_{M,B}"']
        & {xM}
            \arrow[u,"1\pmb{l}_M"']
    \end{tikzcd}\] where 
    \begin{itemize}
        \item the top 4-gon is induced by the holonomy on $M$, see \cite[Theorem 2.2.3]{DX},

        \item the middle-left triangle is filled by invertible modification $R_{M,x,B}$ for the braiding $b$ in $\mathfrak{B}$,

        \item the middle-right triangle is from the definition of $\varsigma^x$,

        \item the bottom-left 4-gon is from the definition of $\vartheta^x$,

        \item the middle 4-gon is induced by the holonomy on $B$, viewed as the canonical local $B$-module,

        \item the bottom-right triangle is filled by the coherence data for half-braiding $\widetilde{b}_{-,B}$, induced by the holonomy on $B$.
    \end{itemize}

    Conversely, given a right $B$-module $M$ and half-braiding $(M,\widetilde{b}^{M}_{-},\widetilde{R}^M_{-,-})$, with a half-syllepsis ${\sigma}^M_{-}$, viewed as an object in $\mathscr{Z}_2(\mathfrak{B}^{2mp} \to \mathscr{Z}_1(\mathbf{Mod}_\mathfrak{B}(B)) )$, we can reassemble the above commutative diagram:
    \begin{itemize}
        \item replace the upper right edge $\widetilde{b}_{M,xB}$ with $\widetilde{b}^M_{xB}$,
        
        \item the entire region enclosed by the outer boundary is filled by the half-syllepsis ${\sigma}^{M}_{-}$,

        \item all cells except the top 4-gon are still filled as above.
    \end{itemize} Then the above diagram induces an invertible 2-isomorphism \[\begin{tikzcd}
        {MxB}
            \arrow[d,"t_{M,xB}"']
            \arrow[r,"b_{M,xB}"]
        & {xBM}
            \arrow[d,"t_{xB,M}"]
        \\ {M \, \Box \, xB}
            \arrow[r,"\widetilde{b}^M_{xB}"']
        & {xB \, \Box \, M}
    \end{tikzcd},\]which is natural in object $x$ in $\mathfrak{B}$. By the 2-universal property of the relative tensor product $\Box_B$, this identifies the half-braiding $\widetilde{b}^M_{xB}$ with $\widetilde{b}_{M,xB}$.
    
    Recall that any right $B$-module $N$ is the condensate for a condensation monad on the free right $B$-module $N \, \Box \, B$. Karoubi completeness of $\mathbf{Mod}_\mathfrak{B}(B)$ implies that half-braiding $\widetilde{b}^M_{N}$ must be equivalent to $\widetilde{b}_{M,N}$ for any right $B$-module $N$. In other word, by naturality we obtain a 2-natural equivalence filling \[\begin{tikzcd}
        {MN}
            \arrow[d,"t_{M,N}"']
            \arrow[r,"b_{M,N}"]
        & {NM}
            \arrow[d,"t_{N,M}"]
        \\ {M \, \Box \, N}
            \arrow[r,"\widetilde{b}_{M,N}"']
        & {N \, \Box \, M}\end{tikzcd}.\] Then we can transport the canonical $B$-balancing on 1-morphism $\widetilde{b}_{M,N} \circ t_{M,N}$ onto 1-morphism $t_{N,M} \circ b_{M,N}$, which further induces a holonomy $h^M$ on $M$, as discussed in \cite[Theorem 2.2.3]{DX}.

    Lastly, one can check that the above two constructions are mutually inverse, and can be extended naturally to an equivalence of braided multifusion 2-categories.
    \end{proof}

\subsection{Factorization of Drinfeld Center into Free Modules and Braided Local Modules}

    \begin{Theorem} \label{thm:E2CenterAndE2LocalModules}
    Let $B$ be a separable braided algebra in a braided fusion 2-category $\mathfrak{B}$. Then we have an equivalence of sylleptic fusion 2-categories $$ \mathscr{Z}_2(\mathfrak{B}) \simeq  \mathscr{Z}_2(\mathfrak{B}^{2mp} \boxtimes \mathbf{Mod}_{\mathfrak{B}}^{\mathbb{E}_2}(B) \to \mathscr{Z}_1(\mathbf{Mod}_\mathfrak{B}(B)) ) . $$
    \end{Theorem}

    \begin{proof}
    Recall from \cite[Remark 2.2.5]{DX}, there is a canonical braided 2-functor $\mathbf{Mod}^{\mathbb{E}_2}_\mathfrak{B}(B) \to \mathscr{Z}_1(\mathbf{Mod}_\mathfrak{B}(B))$, which by Theorem \ref{thm:E2LocalModulesAsE2Centralizer} manifests itself as the braided centralizer of the braided 2-functor $\mathfrak{B}^{2mp} \to \mathscr{Z}_1(\mathbf{Mod}_\mathfrak{B}(B))$ we defined in Lemma \ref{lem:SendingFreeModulesToDrinfeldCenterOfModules}. Therefore, we can embed $\mathfrak{B}^{2mp}$ into its double centralizer $\mathscr{Z}_2(\mathbf{Mod}^{\mathbb{E}_2}_\mathfrak{B}(B) \to \mathscr{Z}_1(\mathbf{Mod}_\mathfrak{B}(B)))$. 
    
    Assume the embedding is an equivalence, i.e. $\mathfrak{B}^{2mp} \simeq \mathscr{Z}_2(\mathbf{Mod}^{\mathbb{E}_2}_\mathfrak{B}(B) \to \mathscr{Z}_1(\mathbf{Mod}_\mathfrak{B}(B)))$, then we can construct the desired equivalence via \[ \mathscr{Z}_2(\mathfrak{B}^{2mp} \boxtimes \mathbf{Mod}_{\mathfrak{B}}^{\mathbb{E}_2}(B) \to \mathscr{Z}_1(\mathbf{Mod}_\mathfrak{B}(B)) ) \simeq \]
    \[\mathscr{Z}_2(\mathfrak{B}^{2mp} \to \mathscr{Z}_2(\mathbf{Mod}^{\mathbb{E}_2}_\mathfrak{B}(B) \to \mathscr{Z}_1(\mathbf{Mod}_\mathfrak{B}(B)))) \simeq \mathscr{Z}_2(\mathfrak{B}^{2mp}) \simeq \mathscr{Z}_2(\mathfrak{B}). \] Notice that the first equivalence in the row exists since we can either take the centralizer of $\mathfrak{B}^{2mp} \boxtimes \mathbf{Mod}_{\mathfrak{B}}^{\mathbb{E}_2}(B)$ in $\mathscr{Z}_1(\mathbf{Mod}_\mathfrak{B}(B))$, or first take the centralizer of $\mathbf{Mod}_{\mathfrak{B}}^{\mathbb{E}_2}(B)$, then take the centralizer of $\mathfrak{B}^{2mp}$. The results have to agree with each other since $\mathfrak{B}^{2mp}$ and $\mathbf{Mod}_{\mathfrak{B}}^{\mathbb{E}_2}(B)$ centralize each other in $\mathscr{Z}_1(\mathbf{Mod}_\mathfrak{B}(B))$ by Theorem \ref{thm:E2LocalModulesAsE2Centralizer}.

    Therefore, it suffices to show that $\mathfrak{B}^{2mp}$ is equivalent to its double centralizer. We prove this using higher condensation theory in Appendix \ref{sec:DoubleCentralizerTheorem}.
    \end{proof}

    \begin{Corollary} \label{cor:E2LocalModulesInNondegenerateBraidedFusion2Cateory}
    Suppose $\mathfrak{B}$ is a non-degenerate braided fusion 2-category and $B$ is a connected separable braided algebra in $\mathfrak{B}$, then we have an equivalence of braided fusion 2-categories $$ \mathscr{Z}_1(\mathbf{Mod}_\mathfrak{B}(B)) \simeq \mathfrak{B}^{2mp} \boxtimes \mathbf{Mod}_{\mathfrak{B}}^{\mathbb{E}_2}(B). $$
    \end{Corollary}
    
    \begin{Corollary} \label{cor:FactorizationOfNondegenerateBraidedFusion2Cateory}
    In particular, for any non-degenerate braided fusion 2-category $\mathfrak{B}$, we always have $\mathscr{Z}_1(\mathfrak{B}) \simeq \mathfrak{B}^{2mp} \boxtimes \mathfrak{B}.$
    \end{Corollary}
    
    \begin{Corollary} \label{cor:E1CenterOfModuleCategoryIsE2LocalModuleInE1Center}
    Any separable braided algebra $B$ in non-degenerate braided fusion 2-category $\mathfrak{B}$ can be viewed as a separable braided algebra in $\mathscr{Z}_1(\mathfrak{B})$, with \[\mathbf{Mod}^{\mathbb{E}_2}_{\mathscr{Z}_1(\mathfrak{B})}(B) \simeq \mathscr{Z}_1(\mathbf{Mod}_\mathfrak{B}(B)).\]
    \end{Corollary}
    
    \begin{Remark} \label{rmk:nondegeneratebraidedmultifusion2categoryisfusion}
    A non-degenerate braided multifusion 2-category $\mathfrak{B}$ is indecomposable. By a similar argument from Remark \ref{rmk:nondegeneratemultifusion2categoryisindecomposable}, we can decompose it as $\mathfrak{B} = \boxplus_{i,j} \, \mathfrak{B}_{ij}$ where each component $\mathfrak{B}_{ii}$ is a non-degenerate braided fusion 2-category. But the existence of braiding forces all off-diagonal components $\mathfrak{B}_{ij}$ to vanish, so $\mathfrak{B}$ must be fusion.
    \end{Remark}
    
    \begin{Remark} \label{rmk:E2centerofE1centeristrivial}
    Drinfeld center of a fusion 2-category is always non-degenerate \cite[Remark 5.3.8]{D9}, i.e. for any fusion 2-category $\mathfrak{C}$, we have $\mathscr{Z}_2(\mathscr{Z}_1(\mathfrak{C}))\simeq \mathbf{2Vect}$.
    \end{Remark}

    % \begin{Remark} % This is not true, see the main theorem in the last section!
    % In contrast with the $\mathbb{E}_1$ centralizer, where for any separable algebra $A$ in an indecomposable multifusion 2-category $\mathfrak{C}$, one can always find some separable algebra $P$ in $\mathfrak{D}:=\mathbf{Mod}^{\mathbb{E}_1}_\mathfrak{C}(P)$ such that $\mathfrak{C} \simeq \mathbf{Mod}^{\mathbb{E}_1}_\mathfrak{D}(P)$ (see \cite[Theorem 5.4.3]{D7}), for a separable braided algebra $B$ in a braided fusion 2-category $\mathfrak{B}$, in general we cannot find any separable braided algebra $Q$ in $\mathfrak{A}:=\mathbf{Mod}^{\mathbb{E}_2}_\mathfrak{B}(B)$ such that $\mathbf{Mod}^{\mathbb{E}_2}_\mathfrak{A}(Q)$ and $\mathfrak{B}$ becomes braided equivalent.
    % \end{Remark}

\section{Sylleptic Centralizers and Symmetric Centers}

In this section, we will continue using notations and conventions defined in the preliminary.

\subsection{Definitions and Basic Properties}

Following the notion of 3-center of sylleptic monoidal 2-categories \cite[Section 5.3]{Cr}, and \cite[Definition 3.10]{KTZ}, we have the the definition of sylleptic centralizer as follows.
\begin{Definition} \label{def:SyllepticCentralizer}
    Sylleptic centralizer (or $\mathbb{E}_3$ centralizer) of a sylleptic 2-functor $F:\mathfrak{S} \to \mathfrak{T}$ between sylleptic monoidal 2-categories $\mathfrak{S}, \mathfrak{T}$ is a sylleptic monoidal 2-category $\mathscr{Z}_{3}(F)$ where:

    \begin{enumerate}
        \item [a.] An object in $\mathscr{Z}_{3}(F)$ consist of an transparent object $x$ in $\mathfrak{T}$, i.e. $\sigma_{F(y),x} \circ b_{x,F(y)} = b_{F(y),x} \circ \sigma_{x,F(y)}$ for any object $y$ in $\mathfrak{S}$.

        \item [b.] 1-morphisms are just 1-morphisms between the underlying objects in $\mathfrak{T}$.
        
        \item [c.] 2-morphisms are just 2-morphisms in $\mathfrak{T}$.

        \item [d.] The entire sylleptic monoidal structure is inherited from that on $\mathfrak{T}$.
    \end{enumerate}
\end{Definition}

\begin{Lemma}
    By definition, there is a canonical forgetful functor $\mathscr{Z}_{3}(F) \to \mathfrak{T}$ preserving syllepses.
\end{Lemma}

\begin{Remark}
    By analogy with Corollary \ref{cor:DrinfeldCenterAsBimoduleFunctors} and Lemma \ref{lem:BraidedCentralizerAsMoritaDual}, we can also write the sylleptic centralizer $\mathscr{Z}_3(F)$ as the Morita dual of some enveloping algebra $\int_{S^2}F$, provided with some rigidity conditions on $\mathfrak{S}$ and $\mathfrak{T}$. More explicitly, suppose $\mathfrak{S}$ and $\mathfrak{T}$ are multifusion, one should expect that 
    \begin{enumerate}
        \item Recall that $\int_{S^1}F := \mathbf{Fun}_{\mathfrak{S} \boxtimes \mathfrak{S}^{2mp}}({}_{\langle F \rangle} \mathfrak{T}_{\langle F \rangle},\mathfrak{S})$.

        \item There is an $\mathfrak{S} \boxtimes \mathfrak{S}^{3mp}$-action on $\mathfrak{S} \boxtimes \mathfrak{S}^{2mp}$, factoring through its sylleptic center. There is a universal $\mathfrak{S} \boxtimes \mathfrak{S}^{3mp}$-balanced 2-functor $\mathfrak{S} \boxtimes \mathfrak{S}^{2mp} \to \mathfrak{S} \boxtimes_{\mathfrak{S} \boxtimes \mathfrak{S}^{3mp}} \mathfrak{S}^{2mp}$ preserving braided monoidal structures.

        \item There is an induced 2-functor $$ \int_{S^2}F := \mathbf{Fun}_{\mathfrak{S} \boxtimes_{\mathfrak{S} \boxtimes \mathfrak{S}^{3mp}} \mathfrak{S}^{2mp}}({}_{\langle F \rangle} \mathfrak{T}_{\langle F \rangle},\mathfrak{S}) \to \int_{S^1}F.$$

        \item The $\int_{S^1}F$-action on $\mathfrak{T}$ can descend to an $\int_{S^2}F$-action on $\mathfrak{T}$.

        \item $\mathscr{Z}_3(F) \simeq \mathbf{End}_{\int_{S^2}F}(\mathfrak{T})$.
    \end{enumerate}
    Interested readers may also consult \cite{BJS} and \cite{BJSS} for a treatment of general $\mathbb{E}_n$ enveloping algebras and $\mathbb{E}_n$ centers, and their applications to braided tensor 1-categories.
\end{Remark}
    
\begin{Definition} \label{def:SymmetricCenter}
Symmetric center (or $\mathbb{E}_3$ center) of a sylleptic monoidal 2-category $\mathfrak{S}$ is defined to be the sylleptic centralizer $\mathscr{Z}_{3}(\mathrm{Id}_{\mathfrak{S}})$. We denote the symmetric center of $\mathfrak{S}$ as $\mathscr{Z}_{3}(\mathfrak{S})$.
\end{Definition}
    
\begin{Lemma}
Symmetric center $\mathscr{Z}_{3}(\mathfrak{S})$ of a sylleptic monoidal 2-category $\mathfrak{S}$ is a symmetric monoidal 2-category.
\end{Lemma}    

\begin{proof}
The lemma follows easily from the given definition. See also \cite[Theorem 5.2]{Cr}.
\end{proof}

\begin{Remark} \label{rmk:S2EnvelopingAlgebra}
    From the point of view of higher Morita theory, one can understand the symmetric monoidal structure on $\mathscr{Z}_{3}(\mathfrak{S})$ as follows. Every sylleptic monoidal 2-category $\mathfrak{S}$ is an $\mathbb{E}_3$ algebra in the symmetric 3-category $\mathbf{2Cat}$, and we expect to generalize the notion of symmetric local modules in Definition \ref{def:E3localmodule} to $\mathbf{2Cat}$, and obtain a 3-category $\mathbf{Mod}^{\mathbb{E}_3}(\mathfrak{S})$ consisting of symmetric local module 2-categories over $\mathfrak{S}$. With some additional assumptions, $\mathbf{Mod}^{\mathbb{E}_3}(\mathfrak{S})$ can be equipped with a sylleptic monoidal structure, where the relative Deligne 2-tensor product $\boxtimes_{\mathfrak{S}}$ provides the monoidal product \cite{D10}. Finally, one expects to recover $\mathscr{Z}_3(\mathfrak{S})$ as the endo-hom of the monoidal unit $\mathfrak{S}$ in the sylleptic monoidal 3-category $\mathbf{Mod}^{\mathbb{E}_3}(\mathfrak{S})$.
\end{Remark}

    \begin{Remark}
    A sylleptic monoidal 2-category $\mathfrak{S}$ is symmetric if and only if the canonical embedding $\mathscr{Z}_3(\mathfrak{S}) \to \mathfrak{S}$ is a sylleptic equivalence. 
    
    One can also observe that $\mathscr{Z}_3(\mathfrak{S}) \simeq \mathscr{Z}_3(\mathfrak{S}^{3mp})$ as symmetric monoidal 2-categories.
    \end{Remark}

\begin{Proposition} \label{prop:SyllepticCentralizerIsSyllepticMultifusion}
    When $F:\mathfrak{S} \to \mathfrak{T}$ is a sylleptic 2-functor between sylleptic multifusion 2-categories, its sylleptic centralizer $\mathscr{Z}_{3}(F)$ is a sylleptic multifusion 2-category.
\end{Proposition}

\begin{proof}
    By definition, sylleptic center $\mathscr{Z}_3(F)$ is a sylleptic monoidal 2-subcategory of $\mathfrak{T}$. Since $\mathfrak{T}$ is finite semisimple, it follows immediately that the underlying 2-category of $\mathscr{Z}_3(F)$ is finite semisimple. Next, by the functorality of taking dual objects in multifusion 2-categories \cite[Appendix A]{DX}, since $\mathscr{Z}_3(F)$ is closed under the monoidal structure of $\mathfrak{T}$, it is also a rigid monoidal 2-category.
\end{proof}

\begin{Proposition} \label{prop:symmetriccenterpreservessumandproduct}
    Let $\mathfrak{S}$ and $\mathfrak{T}$ be sylleptic multifusion 2-categories, then we have the follows:
    \begin{enumerate}
        \item $\mathscr{Z}_3(\mathfrak{S} \boxplus \mathfrak{T}) \simeq \mathscr{Z}_3(\mathfrak{S}) \boxplus \mathscr{Z}_3(\mathfrak{T})$,

        \item $\mathscr{Z}_3(\mathfrak{S} \boxtimes \mathfrak{T}) \simeq \mathscr{Z}_3(\mathfrak{S}) \boxtimes \mathscr{Z}_3(\mathfrak{T})$.
    \end{enumerate}
\end{Proposition}

\begin{proof}
    The first equivalence is by direct inspection. For the second equivalence, notice that 2-functor in one direction $\mathscr{Z}_3(\mathfrak{S}) \boxtimes \mathscr{Z}_3(\mathfrak{T}) \to \mathscr{Z}_3(\mathfrak{S} \boxtimes \mathfrak{T})$ is straight-forward from the universal properties of Deligne 2-tensor product. It follows that this 2-functor induces equivalences on hom categories and preserves sylleptic monoidal structure. Finally, to see the above 2-functor is an equivalence, pick any object $z$ in $\mathscr{Z}_3(\mathfrak{S} \boxtimes \mathfrak{T})$, i.e. a transparent object in $\mathfrak{S} \boxtimes \mathfrak{T}$. In particular, this object is a condensate of a condensation 2-monad in $\mathfrak{S} \otimes \mathfrak{T}$\footnote{This 2-category has objects given by pairs of objects in $\mathfrak{S}$ and $\mathfrak{T}$; the hom category between pairs $(x_0,y_0)$ and $(x_1,y_1)$ is given by the Deligne tensor product $Hom_\mathfrak{S}(x_0,x_1) \boxtimes Hom_\mathfrak{T}(y_0,y_1)$. The Karoubi completion of $\mathfrak{S} \otimes \mathfrak{T}$ is exactly $\mathfrak{S} \boxtimes \mathfrak{T}$. For details see \cite[Lemma 3.3]{D3}.}. Thus it suffices to show that any 2-condensation $x \boxtimes y$ onto $z$ in $\mathfrak{S} \boxtimes \mathfrak{T}$ factors through another 2-condensation $x^\prime \boxtimes y^\prime$ where $x^\prime$ and $y^\prime$ are now transparent objects in $\mathfrak{S}$ and $\mathfrak{T}$, respectively.

    Lastly, since $\mathfrak{S}$ and $\mathfrak{T}$ are multifusion 2-categories, one can take left adjoint of the canonical embeddings $\mathscr{Z}_3(\mathfrak{S}) \to \mathfrak{S}$ and $\mathscr{Z}_3(\mathfrak{T}) \to \mathfrak{T}$. Let's denote them as $L^\mathfrak{S}:\mathfrak{S} \to \mathscr{Z}_3(\mathfrak{S})$ and $L^\mathfrak{T}:\mathfrak{T} \to \mathscr{Z}_3(\mathfrak{T})$, respectively. Then we observe that the counits $\epsilon^\mathfrak{S}_x:L^\mathfrak{S}(x) \to x$ and $\epsilon^\mathfrak{T}_y:L^\mathfrak{T}(y) \to y$, for object $x$ from $\mathfrak{S}$ and $y$ from $\mathfrak{T}$, can be both extended to be injective parts of 2-condensations. 

    Indeed, $L^\mathfrak{S}(x)$ and $L^\mathfrak{T}(y)$ satisfy all the requirements for transparent objects $x^\prime$ and $y^\prime$ mentioned above, i.e. given any transparent object $z$ and any 2-condensation $x \boxtimes y$ onto $z$ in $\mathfrak{S} \boxtimes \mathfrak{T}$, it factors through 2-condensations from $x$ to $L^\mathfrak{S}(x)$ and from $y$ to $L^\mathfrak{T}(y)$ as the following diagram.
\end{proof}

    \[\begin{tikzcd}
        {x \boxtimes y}
            \arrow[r,two heads,bend right=10pt]
            \arrow[d,two heads,bend right=10pt,shift right=5pt]
        & {x \boxtimes L^\mathfrak{T}(y)}
            \arrow[l,hook,bend right=10pt]
            \arrow[d,two heads,bend right=10pt,shift right=5pt]
        & {}
        \\ {L^\mathfrak{S}(x) \boxtimes y}
            \arrow[r,two heads,bend right=10pt]
            \arrow[u,hook,bend right=10pt,shift right=5pt]
        & {L^\mathfrak{S}(x) \boxtimes L^\mathfrak{T}(y)}
            \arrow[l,hook,bend right=10pt]
            \arrow[u,hook,bend right=10pt,shift right=5pt]
            \arrow[r,two heads,bend right=10pt]
        & {z}
            \arrow[l,hook,bend right=10pt,shift right=2pt]
    \end{tikzcd}\]

\begin{Remark}
    Another way to prove that $\mathscr{Z}_3(\mathfrak{S} \boxtimes \mathfrak{T}) \simeq \mathscr{Z}_3(\mathfrak{S}) \boxtimes \mathscr{Z}_3(\mathfrak{T})$ for sylleptic multifusion 2-categories $\mathfrak{S}$ and $\mathfrak{T}$ utilizes the construction of the $S^2$-enveloping algebra we proposed in Remark \ref{rmk:S2EnvelopingAlgebra}. By analogy with Proposition \ref{prop:syllepticcenterpreservessumandproduct}, we expect that
    \begin{enumerate}
        \item $\int_{S^2}\mathfrak{S} \boxtimes \int_{S^2}\mathfrak{T} \simeq \int_{S^2} (\mathfrak{S} \boxtimes \mathfrak{T})$.

        \item Moreover, the above equivalence should be compatible with the actions on $\mathfrak{S} \boxtimes \mathfrak{T}$ from both sides.

        \item Furthermore, it induces $\mathscr{Z}_3(\mathfrak{S}) \boxtimes \mathscr{Z}_3(\mathfrak{T}) \simeq \mathbf{End}_{\int_{S^2}\mathfrak{S}}(\mathfrak{S}) \boxtimes \mathbf{End}_{\int_{S^2}\mathfrak{T}}(\mathfrak{T}) \simeq \mathbf{End}_{\int_{S^2}\mathfrak{S} \boxtimes \int_{S^2}\mathfrak{T}}(\mathfrak{S} \boxtimes \mathfrak{T}) \simeq \mathbf{End}_{\int_{S^2}(\mathfrak{S} \boxtimes \mathfrak{T})}(\mathfrak{S} \boxtimes \mathfrak{T}) \simeq \mathscr{Z}_3(\mathfrak{S} \boxtimes \mathfrak{T})$.

        \item Finally, the above equivalence should preserve sylleptic structures.
    \end{enumerate}
\end{Remark}

\subsection{Non-Degeneracy of Sylleptic Categories}

\begin{Definition} \label{def:WeakNonDegenerateSyllepticFusion2Category}
    Given a sylleptic multifusion 2-category $\mathfrak{S}$, it is said to be weakly non-degenerate if $\mathscr{Z}_3(\mathfrak{S}) \simeq \mathbf{2Vect}$.
\end{Definition}

\begin{Remark} \label{rmk:nondegeneratesyllepticmultifusion2categoryisfusion}
    By analogy with Remark \ref{rmk:nondegeneratemultifusion2categoryisindecomposable} and Remark \ref{rmk:nondegeneratebraidedmultifusion2categoryisfusion}, any weakly non-degenerate sylleptic multifusion 2-category $\mathfrak{S}$ is fusion.
\end{Remark}

\begin{Remark} \label{rmk:E3centerofE2centerisNOTtrivial}
    In contrast with Remark \ref{rmk:E1centerofE0centeristrivial} and Remark \ref{rmk:E2centerofE1centeristrivial}, for a braided fusion 2-category $\mathfrak{B}$, its sylleptic center $\mathscr{Z}_2(\mathfrak{B})$ is \textit{not} weakly non-degenerate in general. A first counter-example is $\mathfrak{B} = \mathbf{Mod}(\mathbf{Rep}(G))$, when $\mathscr{Z}_2(\mathfrak{B}) \simeq \mathfrak{B}$ and hence $\mathscr{Z}_3(\mathscr{Z}_2(\mathfrak{B})) \simeq \mathfrak{B}$, see \cite[Proposition 3.11]{KTZ}.
\end{Remark}

Weakly non-degeneracy of sylleptic fusion 2-categories is significantly different from the non-degeneracy of fusion 2-categories (Definition \ref{def:NonDegenerateFusion2Category}) and braided fusion 2-categories (Definition \ref{def:NonDegenerateBraidedFusion2Category}). We will compare them under the general picture of condensations of higher fusion categories, \cite{GJF,JF,KZ20,KZ21}.

\begin{Definition} \label{def:NonDegenerateSyllepticFusion2Category}
    A sylleptic fusion 2-category $\mathfrak{S}$ is non-degenerate if its double suspension, $\Sigma^2 \mathfrak{S}$, is a non-degenerate fusion 4-category, i.e. $\mathscr{Z}_1(\Sigma^2 \mathfrak{S}) \simeq \mathbf{4Vect}$.
\end{Definition}

\begin{Corollary}
    By \cite[Corollary IV.3]{JF}, a sylleptic fusion 2-category $\mathfrak{S}$ is non-degenerate if and only if $\Sigma \mathfrak{S}$ is a non-degenerate braided fusion 3-category, i.e. $\mathscr{Z}_2(\Sigma \mathfrak{S}) \simeq \mathbf{3Vect}$.
\end{Corollary}

\begin{Remark}
    Recall that the suspension operation $\Sigma$ is defined for any monoidal $n$-category $\mathfrak{C}$ by first take the one-point delooping $B\mathfrak{C}$ then apply the Karoubi completion \cite{GJF}. Moreover, if $\mathfrak{C}$ is $\mathbb{E}_m$ monoidal then $\Sigma \mathfrak{C}$ is $\mathbb{E}_{m-1}$ monoidal with product given by relative tensor product $\boxtimes_\mathfrak{C}$. Conjecturally, for sylleptic fusion 2-category $\mathfrak{S}$, $\Sigma \mathfrak{S}$ is equivalent as braided fusion 3-category to $\mathbf{Mod}(\mathfrak{S})$, the 3-category of fully dualizable $\mathfrak{S}$-module 2-categories; also, $\Sigma^2 \mathfrak{S}$ is equivalent as fusion 4-category to $\mathbf{Mod}(\mathbf{Mod}(\mathfrak{S}))$, the 4-category of fully dualizable $\mathbf{Mod}(\mathfrak{S})$-module 3-categories.
\end{Remark}

\begin{Conjecture}
    We expect Theorem \ref{thm:E1LocalModulesAsE1Centralizer}, Theorem \ref{thm:DrinfeldCenterAndE1LocalModules}, Corollary \ref{cor:E1LocalModulesInNondegenerateFusion2Cateory}, Corollary \ref{cor:FactorizationOfNondegenerateFusion2Cateory} and Corollary \ref{cor:E0CenterOfModuleCategoryIsE1LocalModuleInE0Center} to generalize to fusion 4-categories and separable algebras within them.

    Also, we expect Theorem \ref{thm:E2LocalModulesAsE2Centralizer}, Theorem \ref{thm:E2CenterAndE2LocalModules}, Corollary \ref{cor:E2LocalModulesInNondegenerateBraidedFusion2Cateory}, Corollary \ref{cor:FactorizationOfNondegenerateBraidedFusion2Cateory} and Corollary \ref{cor:E1CenterOfModuleCategoryIsE2LocalModuleInE1Center} to generalize to braided fusion 3-categories and separable $\mathbb{E}_2$ algebras within them.
\end{Conjecture}

Conversely, \cite[Theorem 5]{JF} states the follows.

\begin{Proposition} \label{prop:NonDegenerateFusion4CategoryAndE3StronglyFusion2Category}
    A non-degenerate fusion 4-category $\mathfrak{C}$ with $\Omega^3\mathfrak{C} \simeq \mathbf{Vect}$ or $\Omega^3\mathfrak{C} \simeq \mathbf{sVect}$ is always equivalent to $\Sigma^2 \mathfrak{S}$ where $\mathfrak{S} = \Omega^2 \mathfrak{C}$ is a strongly fusion sylleptic 2-category in the sense of Definition \ref{def:StronglyFusion2Category}.
\end{Proposition}

\begin{Remark}
    Similarly, a non-degenerate braided fusion 3-category $\mathfrak{B}$ with $\Omega^2\mathfrak{B} \simeq \mathbf{Vect}$ or $\Omega^2\mathfrak{B} \simeq \mathbf{sVect}$ is equivalent to $\Sigma \mathfrak{S}$ where $\mathfrak{S} = \Omega \mathfrak{B}$ is a strongly fusion sylleptic 2-category.
\end{Remark}

Recall \cite[Theorem 3.43]{KZ20}, one can characterize symmetric center of sylleptic fusion 2-category $\mathfrak{S}$ as follows.

\begin{Proposition}
    $\mathscr{Z}_3(\mathfrak{S}) = \Omega \mathscr{Z}_2(\Sigma \mathfrak{S}) = \Omega^2 \mathscr{Z}_1(\Sigma^2 \mathfrak{S})$.
\end{Proposition}

\begin{Corollary}
    For strongly fusion sylleptic 2-categories, the notion of weak non-degeneracy (Definition \ref{def:WeakNonDegenerateSyllepticFusion2Category}) agrees with non-degeneracy (Definition \ref{def:NonDegenerateSyllepticFusion2Category}).
\end{Corollary}

\begin{Remark}
    In general, a weakly non-degenerate sylleptic fusion 2-category $\mathfrak{S}$ is not necessarily non-degenerate. One could have $\mathscr{Z}_3(\mathfrak{S}) \simeq \mathbf{2Vect}$ but $\mathscr{Z}_2(\Sigma \mathfrak{S}) \not\simeq \mathbf{3Vect}$. It would be an interesting question to classify weakly non-degenerate sylleptic fusion 2-categories modulo non-degenerate ones.
\end{Remark}

\begin{Remark}
    A non-degenerate sylleptic fusion 2-category $\mathfrak{S}$ is always weakly non-degenerate. But for a non-degenerate fusion 4-category $\mathfrak{C}$, in general sylleptic 2-category $\Omega^2 \mathfrak{C}$ is not weakly non-degenerate. However, $\mathfrak{C}$ should be Morita equivalent to a non-degenerate fusion 4-category $\mathfrak{D}$ with $\Omega^2 \mathfrak{D}$ weakly non-degenerate via \textit{de-equivariantization} in \cite[Remark V.6]{JF} and \cite{JFY22}. 
    
    More explicitly, $\Omega^3 \mathfrak{C}$ is a symmetric fusion 1-category, which is classified by Deligne to be either $\mathbf{Rep}(G)$ for some finite group $G$ or $\mathbf{Rep}(G,z)$ for some finite super-group $(G,z)$. In the Tannakian case, $\Omega^3 \mathfrak{C}$ admits a fiber functor to $\mathbf{Vect}$ while in general $\Omega^3 \mathfrak{C}$ admits a fiber functor to $\mathbf{sVect}$. This fiber functor induces a monoidal 4-functor $\Sigma^3 \Omega^3 \mathfrak{C} \to \Sigma^3 \mathbf{sVect}$. Hence, we can construct the de-equivariantization via the relative tensor product \[\mathfrak{M} := \mathfrak{C} \boxtimes_{\Sigma^3 \Omega^3 \mathfrak{C}} \Sigma^3 \mathbf{sVect},\] which is a module 4-category over $\mathfrak{C}$. Its Morita dual $\mathfrak{D}:= \mathbf{End}_{\mathfrak{C}}(\mathfrak{M})^{1mp}$ is now a non-degenerate fusion 4-category with $\Omega^3 \mathfrak{D} \simeq \mathbf{sVect}$, hence by Proposition \ref{prop:NonDegenerateFusion4CategoryAndE3StronglyFusion2Category}, $\Omega^2 \mathfrak{D}$ is a weakly non-degenerate sylleptic fusion 2-category.
\end{Remark}

\begin{Definition} \label{def:StronglyFusion2Category}
    A fusion 2-category $\mathfrak{C}$ is said to be bosonic strongly fusion if $\Omega \mathfrak{C} \simeq \mathbf{Vect}$. A fusion 2-category $\mathfrak{C}$ is said to be fermionic strongly fusion if $\Omega \mathfrak{C} \simeq \mathbf{sVect}$.
\end{Definition}

By \cite[Theorem A, Theorem B]{JFY21}, we can characterize strongly fusion 2-categories by 2-categories of (possibly twisted) $G$-graded separable categories or super-categories.

\begin{Proposition}
    Simple objects in a bosonic or fermionic strongly fusion 2-category are invertible.
\end{Proposition}

\begin{Corollary}
    A fusion 2-category $\mathfrak{C}$ is bosonic strongly fusion if and only if it is equivalent to $\mathbf{2Vect}^\pi_G$ for finite group $G = \pi_0(\mathfrak{C})$ with pentagonator twisted by 4-cocycle $\pi$ on $G$.
\end{Corollary}

The classification of fermionic strongly fusion 2-categories is more complicated. Thanks to the theory of relative 2-Deligne tensor product, one can reinterpret this problem into the classification of $G$-graded extensions for fusion 2-category $\mathbf{2sVect}$. Following the general idea of \cite{ENO09,DN}, Décoppet obtained the following proposition in \cite[Proposition 4.3.2]{D10}.

\begin{Proposition}
    A fermionic strongly fusion 2-category $\mathfrak{C}$ is a $G$-graded extension of $\mathbf{2sVect}$ for finite group $G = \pi_0(\mathfrak{C})$, together with a class $\varpi \in \mathrm{H}^2(G;\mathbb{Z}/2)$ and a class $\pi \in \mathrm{SH}^{4+\varpi}(G)$.
\end{Proposition}

\begin{Definition}
    A strongly fusion sylleptic 2-category is a sylleptic fusion 2-category whose underlying fusion 2-category is either bosonic or fermionic strongly fusion.
\end{Definition}

Following the paradigm of classification of braided extensions of braided fusion 1-categories \cite[Section 8.3]{DN}, we conjecture that strongly fusion sylleptic 2-categories can be classified as follows.

\begin{Definition}
    The sylleptic Picard 4-group of a sylleptic fusion 2-category $\mathfrak{S}$, denoted by $\mathbf{Pic}_{\mathbb{E}_3}(\mathfrak{S})$, is defined to be the sub-3-groupoid spanned by invertible objects in the 3-categorical sylleptic center $\mathscr{Z}_2(\mathbf{Mod}(\mathfrak{S}))$.
\end{Definition}

Fix a finite Abelian group $G$. 

\begin{Definition}
    A $G$-graded sylleptic extension of sylleptic fusion 2-category $\mathfrak{S}$ is a sylleptic fusion 2-category $\bigoplus_{g \in G} \mathfrak{T}_g$ with identity component $\mathfrak{T}_e \simeq \mathfrak{S}$.
\end{Definition}

\begin{Conjecture}
    There is an equivalence between $G$-graded sylleptic extensions of sylleptic fusion 2-category $\mathfrak{S}$ and the 3-category of sylleptic 3-functors from $G$ into $\mathbf{Pic}_{\mathbb{E}_3}(\mathfrak{S})$.
\end{Conjecture}

\begin{Conjecture}
    $\mathbf{Pic}_{\mathbb{E}_3}(\mathbf{2Vect}) \simeq B^3 \Bbbk^\times, \mathbf{Pic}_{\mathbb{E}_3}(\mathbf{2sVect}) \simeq \mathbf{Pic}(\mathbf{2sVect})$.
\end{Conjecture}

\begin{Corollary}
    A bosonic strongly fusion sylleptic 2-category is equivalent to $\mathbf{2Vect}^q_A$ where $A$ is a finite Abelian group with a class $q \in \mathrm{H}^3(B^3 A;B^3\Bbbk^\times)$. By \cite[Section 2.2]{JFY22}, if the Pontryagin dual of $A_2 := \mathrm{Hom}(\mathbb{Z}/2,A)$ vanishes, then this class is determined by a skew-symmetric bilinear pairing $A \times A \to \Bbbk^\times$. The pair $(A,q)$ is called a finite pre-symplectic Abelian group. In general, the class $q$ also depends on the Pontryagin dual of $A_2$. From a field-theoretic perspective \cite{JFR}, this extra component measures the partition functions on Klein bottles and only depends on the braiding of the 2-category.
\end{Corollary}

\begin{Corollary}
    A fermionic strongly fusion sylleptic 2-category is a $G$-graded sylleptic extension of $\mathbf{2sVect}$ for a finite Abelian group $G$, which is classified by a class $\pi \in \mathrm{SH}^6(B^3 G)$. By \cite[Proposition 2.3]{JFY22}, this class is determined by a skew-symmetric bilinear pairing $q: A \times A \to \Bbbk^\times$. In other word, fermionic strongly fusion sylleptic 2-categories arise as linearization of finite pre-symplectic Abelian groups.
\end{Corollary}

\begin{Conjecture}[{\cite[Theorem 2.4]{JFY22}}]
    Non-degeneracy condition of a finite pre-symplectic Abelian group is equivalent to the non-degeneracy of its linearization as a sylleptic fusion 2-category.
\end{Conjecture}

\subsection{Factorization of Sylleptic Center into Free Modules and Symmetric Local Modules}

\begin{Lemma} \label{lem:freeE3localmodule}
    Suppose $\mathfrak{S}$ is a sylleptic multifusion 2-category and $S$ is a separable symmetric algebra in $\mathfrak{S}$, then there is a sylleptic 2-functor $\mathfrak{S}^{3mp} \to \mathscr{Z}_2(\mathbf{Mod}_\mathfrak{S}(S))$.
\end{Lemma}

\begin{proof}
    For any object $x$ in $\mathfrak{S}$, we would like to assign it with half-syllepsis $(x \, \Box \, S, \widehat{\sigma}^{xS}_{-})$ in $\mathbf{Mod}_\mathfrak{S}(S)$. By Lemma \ref{lem:BraidedCentralizerAsMoritaDual}, this is equivalent to a $\int_{S^1} \mathbf{Mod}_\mathfrak{S}(S)$-module 2-functor on $\mathbf{Mod}_\mathfrak{S}(S)$; referring to the notations in Lemma \ref{lem:SendingFreeModulesToDrinfeldCenterOfModules}, it requires:
    \begin{itemize}
        \item 2-functor $\mathbf{Mod}_\mathfrak{S}(S) \to \mathbf{Mod}_\mathfrak{S}(S)$; \[(M,n^M,\nu^M,\rho^M) \mapsto (x \, \Box \, M, x \, \Box \, n^M, x \, \Box \, \nu^M, x \, \Box \, \rho^M).\]

        \item 2-natural equivalence $\varsigma^x_{M,N}: x \, \Box \, (M \, \Box_S \, N) \simeq (x \, \Box \, M) \, \Box_S \, N$ given on right $S$-modules $M$ and $N$.

        \item 2-natural equivalence $\vartheta^x_{M,N}: x \, \Box \, (M \, \Box_S \, N) \simeq M \, \Box_S \, (x \, \Box \, N)$ given on right $S$-modules $M$ and $N$.

        \item Invertible modifications \[\begin{tikzcd}[sep=40pt]
            {x \, \Box \, (M \, \Box_S \, N)}
                \arrow[r,"x \, \Box \, \widetilde{b}_{M,N}"]
                \arrow[d,"\varsigma^x_{M,N}"']
            & {x \, \Box \, (N \, \Box_S \, M)}
                \arrow[d,"\vartheta^x_{N,M}"]
            \\ {(x \, \Box \, M) \, \Box_S \, N}
                \arrow[r,"\widetilde{b}_{xM,N}"']
                \arrow[ur,Rightarrow,shorten <=20pt, shorten >=20pt,"\varphi^x_{M,N}"]
            & {N \, \Box_S \, (x \, \Box \, M)}
        \end{tikzcd},\] and \[
        \begin{tikzcd}[sep=40pt]
            {x \, \Box \, (M \, \Box_S \, N)}
                \arrow[r,"x \, \Box \, \widetilde{b}_{M,N}"]
                \arrow[d,"\vartheta^x_{M,N}"']
            & {x \, \Box \, (N \, \Box_S \, M)}
                \arrow[d,"\varsigma^x_{N,M}"]
            \\ {M \, \Box_S \, (x \, \Box \, N) }
                \arrow[r,"\widetilde{b}_{M,xN}"']
                \arrow[ur,Rightarrow,shorten <=20pt, shorten >=20pt,"\psi^x_{M,N}"]
            & {(x \, \Box \, N) \, \Box_S \, M}
        \end{tikzcd},\] given on right $S$-modules $M$ and $N$.
    \end{itemize}
    Here, 2-natural equivalences $\varsigma^x$ and $\vartheta^x$, which witness that the underlying 2-functor is $\mathbf{Mod}_\mathfrak{S}(S)$-bilinear, are induced in the same way as in Lemma \ref{lem:SendingFreeModulesToDrinfeldCenterOfModules}. Two additional invertible modifications are induced via the 2-universal property of $x \, \Box \, t_{M,N}: x \, \Box \, M \, \Box \, N \to x \, \Box \, (M \, \Box_S \, N)$ by the following two diagrams:
    \[\begin{tikzcd}
        {x \, \Box \, (M \, \Box_S \, N)}
            \arrow[ddd,"\varsigma^x_{M,N}"']
            \arrow[rrr,"x \, \Box \, \widetilde{b}_{M,N}"]
        & {}
        & {}
        & {x \, \Box \, (N \, \Box_S \, M)}
            \arrow[ddd,"\vartheta^x_{N,M}"]
        \\ {}
        & {xMN}
            \arrow[d,equal]
            \arrow[r,"1b"]
            \arrow[ul,"1t"']
        & {xNM}
            \arrow[d,"b1"'{name=A},bend right=20pt,shift right=10pt]
            \arrow[ur,"1t"]
        & {}
        \\ {}
        & {xMN}
            \arrow[ur,Rightarrow,"S",shorten <=10pt, shorten >=10pt,xshift=-10pt]
            \arrow[r,"b_2"']
            \arrow[dl,"t"]
        & {NxM}
            \arrow[dr,"t"']
            \arrow[u,"b1"'{name=B},bend right=20pt,shift right=10pt]
            \arrow[Rightarrow,"\sigma 1",shorten <=10pt, shorten >=10pt,from=A,to=B]
        & {}
        \\ {(x \, \Box \, M) \, \Box_S \, N}
            \arrow[rrr,"\widetilde{b}_{xM,N}"']
        & {}
        & {}
        & {N \, \Box_S \, (x \, \Box \, M)}
    \end{tikzcd},\]
    where \begin{itemize}
        \item the upper and lower quadrilaterals are both filled by coherence between braidings and relative tensor products,

        \item the left quadrilateral is from the definition of $\varsigma^x$,

        \item the right quadrilateral is from the definition of $\vartheta^x$,

        \item the middle-right bigon is filled by the currying of syllepsis $\sigma_{N,x} \, \Box \, M$ viewed as 2-isomorphism between $b_{x,N} \, \Box \, M$ and $b_{N,x} \, \Box \, M$,

        \item the middle-left square is filled by the coherence data $S_{x,M,N}$ for the braiding $b$;
    \end{itemize}
    and \[\begin{tikzcd}
        {x \, \Box \, (M \, \Box_S \, N)}
            \arrow[ddd,"\vartheta^x_{M,N}"']
            \arrow[rrr,"x \, \Box \, \widetilde{b}_{M,N}"]
        & {}
        & {}
        & {x \, \Box \, (N \, \Box_S \, M)}
            \arrow[ddd,"\varsigma^x_{N,M}"]
        \\ {}
        & {xMN}
            \arrow[r,"1b"]
            \arrow[ul,"1t"']
        & {xNM}
            \arrow[ur,"1t"]
        & {}
        \\ {}
        & {MxN}
            \arrow[r,"b"']
            \arrow[u,"b1"]
            \arrow[dl,"t"]
        & {xNM}
            \arrow[ul,Rightarrow,"R"',shorten <=10pt, shorten >=10pt]
            \arrow[dr,"t"']
            \arrow[u,equal]
        & {}
        \\ {M \, \Box_S \, (x \, \Box \, N)}
            \arrow[rrr,"\widetilde{b}_{M,xN}"']
        & {}
        & {}
        & {(x \, \Box \, N) \, \Box_S \, M}
    \end{tikzcd},\] where \begin{itemize}
        \item the upper and lower quadrilaterals are both filled by coherence between braidings and relative tensor products,

        \item the left quadrilateral is from the definition of $\vartheta^x$,

        \item the right quadrilateral is from the definition of $\varsigma^x$,

        \item the middle square is filled by the coherence data $R_{M,x,N}$ for the braiding $b$.
    \end{itemize}

    More concretely, the half-syllepsis on $x \, \Box \, S$ is constructed as \[\widehat{\sigma}^{xS}_M := \begin{tikzcd}[sep=30pt]
        {x \, \Box \, M}
            \arrow[rr,equal]
            \arrow[d,"x \, \Box \, \pmb{l}_M"']
        & {}
        & {x \, \Box \, M}
            \arrow[d,"x \, \Box \, \pmb{l}_M"]
        \\ {x \, \Box \, (S \, \Box_S \, M)}
            \arrow[d,"\varsigma^x_{S,M}"']
            \arrow[r,"x \, \Box \, \widetilde{b}_{S,M}"]
        & {x \, \Box \, (M \, \Box_S \, S)}
            \arrow[d,"\vartheta^x_{M,S}"]
            \arrow[r,"x \, \Box \, \widetilde{b}_{M,S}"]
            \arrow[ur,"x \, \Box \, \pmb{r}_M"]
            \arrow[ul,"x \, \Box \, \pmb{r}_M"']
        & {x \, \Box \, (S \, \Box_S \, M)}
            \arrow[d,"\varsigma^x_{S,M}"]
        \\ {(x \, \Box \, S) \, \Box_S \, M}
            \arrow[r,"\widetilde{b}_{xS,M}"']
            \arrow[ur,Rightarrow,shorten <=20pt, shorten >=20pt,"\xi^x_{S,M}"]
        & {M \, \Box_S \, (x \, \Box \, S)}
            \arrow[r,"\widetilde{b}_{M,xS}"']
            \arrow[ur,Rightarrow,shorten <=20pt, shorten >=20pt,"\zeta^x_{M,S}"]
        & {(x \, \Box \, S) \, \Box_S \, M}
    \end{tikzcd},\]given on right $S$-module $M$. The upper triangles witness the compatibility between braiding $\widetilde{b}$ and unitors $\pmb{l}, \pmb{r}$ in the braided monoidal 2-category $\mathbf{Mod}_\mathfrak{S}(S)$; or equivalently, they follow from the observation that $S$ is equipped with the canonical $\mathbb{E}_3$-local $S$-module structure, hence it always lies in the sylleptic center $\mathscr{Z}_2(\mathbf{Mod}_\mathfrak{S}(S))$ by Remark \ref{rmk:EmbeddingE3LocalModulesInSyllepticCenter}.

    By the naturality of monoidal product $\Box$, we can easily promote the above assignment to a 2-functor from $\mathfrak{S}$ to $\mathscr{Z}_2(\mathbf{Mod}_\mathfrak{S}(S))$. Finally, one needs to show that this 2-functor can be promoted to be a sylleptic 2-functor $\mathfrak{S}^{3mp} \to \mathscr{Z}_2(\mathbf{Mod}_\mathfrak{S}(S))$.
    
    The braided 2-functor structure on $\mathfrak{S} \to \mathscr{Z}_2(\mathbf{Mod}_\mathfrak{S}(S))$ is induced as follows. On the level of underlying objects, a detailed construction has been given in \cite[Proposition 3.9]{DY}. Take two objects $x$ and $y$ in $\mathfrak{S}$, to check that half-syllepses $\widehat{\sigma}^{xS}_-$, $\widehat{\sigma}^{yS}_-$ and $\widehat{\sigma}^{xyS}_-$ are compatible, it is enough to check the compatibility between invertible modifications $\varphi^x$, $\varphi^y$, $\varphi^{xy}$ and $\psi^x$, $\psi^y$, $\psi^{xy}$, respectively. The later compatibility conditions are straight-forward.

    Lastly, we check that the above data satisfies condition (\ref{eqn:Sylleptic2Functor}), hence it is upgraded to a sylleptic 2-functor $\mathfrak{S}^{3mp} \to \mathscr{Z}_2(\mathbf{Mod}_\mathfrak{S}(S))$.
\end{proof}

\begin{Theorem} \label{thm:E3LocalModulesAsE3Centralizer}
    Suppose $\mathfrak{S}$ is a sylleptic multifusion 2-category and $S$ is a separable symmetric algebra in $\mathfrak{S}$, then we have an equivalence of sylleptic multifusion 2-categories:
    $$ \mathbf{Mod}^{\mathbb{E}_3}_\mathfrak{S}(S) \simeq \mathscr{Z}_3(\mathfrak{S}^{3mp} \to \mathscr{Z}_2(\mathbf{Mod}_\mathfrak{S}(S))). $$
\end{Theorem}

\begin{proof}
    To go from left to right, one needs to show that free $S$-modules and $\mathbb{E}_3$ local $S$-modules centralize each other in $\mathscr{Z}_2(\mathbf{Mod}_\mathfrak{S}(S))$. Take any object $x$ in $\mathfrak{S}$ and any $\mathbb{E}_3$ local $S$-module $M$, we would like to show that \[\widetilde{\sigma}_{M,xS} \circ \widetilde{b}_{xS,M} = \widetilde{b}_{xS,M} \circ \widehat{\sigma}^{xS}_M.\] We save the diagrams of the proof to the Appendix \ref{sec:proofE3LocalModulesAsE3Centralizer}. By the 2-universal property of relative tensor products, it suffices to show that the diagram given in Step 1 is equal to that given in Step 9. In the first step, we apply the definition of half-syllepsis $\widetilde{\sigma}_{M,xS}$ from Remark \ref{rmk:EmbeddingE3LocalModulesInSyllepticCenter} to get Step 2. Then by monoidality of half-syllepsis (see Equation (\ref{eqn:braidedcentralizerproduct}) in Definition \ref{def:BraidedCentralizer}), we obtain Step 3. Since the braiding $b$ is assumed to be a adjoint 2-natural equivalence, we can pull the node $\xi$ down along the blue arrow, then take the top strand $b_{xS,M}$ close to the another strand $b_{xS,M}$ above $t$ along the red arrow, and apply the adjoint naturality to get Step 4. From Step 4 to Step 5, we again use the adjoint 2-natural equivalence of $b$ and its coherence data $S$, as depicted by the red arrow. Now in Step 5, we create a pair of invertible modification, denoted by $\Sigma$ and $\Sigma^{-1}$ as depicted in Step 6, from the definition of $\varsigma^x_{S,M}$, see Lemma \ref{lem:SendingFreeModulesToDrinfeldCenterOfModules}. Then we pull the red node $x \, \Box \, \sigma_{S,M}$ down, until it passes through the strand $1t$ along with its two legs, as shown by the red arrows. This will create two nodes labelled by $1\xi$ in Step 7. From Step 7 to Step 8, we use the fact that $S$ lies in the sylleptic center $\mathscr{Z}_2(\mathbf{Mod}_\mathfrak{S}(S))$ via its canonical $\mathbb{E}_3$ local $S$-module structure, see Remark \ref{rmk:EmbeddingE3LocalModulesInSyllepticCenter}. Lastly, we apply the definition of $\widehat{\sigma}^{xS}_M$ from Lemma \ref{lem:freeE3localmodule} to replace all blue nodes in Step 8 and get Step 9.

    Conversely, we would like to show that a half-syllepsis centralizing all free $S$-modules is induced by an $\mathbb{E}_3$ local $S$-module. Take a half-syllepsis $(M,\widehat{\sigma}^M_{-})$ from $\mathscr{Z}_2(\mathbf{Mod}_\mathfrak{S}(S))$, such that for any object $x$ in $\mathfrak{S}$ one has \[\widehat{\sigma}^M_{xS} \circ \widetilde{b}_{xS,M} = \widetilde{b}_{xS,M} \circ \widehat{\sigma}^{xS}_M.\] Using diagrams from Appendix \ref{sec:proofE3LocalModulesAsE3Centralizer} again, we still can go from Step 2 to Step 9. Then apply the above condition to Step 9, we obtain a diagram like Step 1, but with node $\widetilde{\sigma}_{M,xS}$ replaced by $\widehat{\sigma}^M_{xS}$. Since every right $S$-module is the condensate of some free right $S$-module, this implies that $M$ can be equipped with an $\mathbb{E}_3$ local $S$-module structure and $\widehat{\sigma}^M_{-} = \widetilde{\sigma}_{M,-}$ as invertible modifications.
\end{proof}

\begin{Theorem} \label{thm:E3CenterAndE3LocalModules}
    Suppose $\mathfrak{S}$ is a sylleptic fusion 2-category and $S$ is a separable symmetric algebra in $\mathfrak{S}$, then we have an equivalence of symmetric fusion 2-categories:
    $$ \mathscr{Z}_3(\mathfrak{S}) \simeq \mathscr{Z}_3(\mathfrak{S}^{3mp} \boxtimes \mathbf{Mod}^{\mathbb{E}_3}_\mathfrak{S}(S) \to \mathscr{Z}_2(\mathbf{Mod}_\mathfrak{S}(S))). $$
\end{Theorem}

\begin{proof}
    First, by Theorem \ref{thm:E3LocalModulesAsE3Centralizer}, $\mathbf{Mod}^{\mathbb{E}_3}_\mathfrak{S}(S)$ is equivalent to the sylleptic centralizer of the free $S$-modules $\mathfrak{S}^{3mp} \to \mathscr{Z}_2(\mathbf{Mod}_\mathfrak{S}(S))$, hence one has a canonical embedding of $\mathfrak{S}^{3mp}$ into its double centralizer $\mathscr{Z}_3(\mathbf{Mod}^{\mathbb{E}_3}_\mathfrak{S}(S) \to \mathscr{Z}_2(\mathbf{Mod}_\mathfrak{S}(S)))$.

    The desired equivalence is constructed via \[\mathscr{Z}_3(\mathfrak{S}^{3mp} \boxtimes \mathbf{Mod}^{\mathbb{E}_3}_\mathfrak{S}(S) \to \mathscr{Z}_2(\mathbf{Mod}_\mathfrak{S}(S))) \simeq \] \[\mathscr{Z}_3(\mathfrak{S}^{3mp} \to \mathscr{Z}_3(\mathbf{Mod}^{\mathbb{E}_3}_\mathfrak{S}(S) \to \mathscr{Z}_2(\mathbf{Mod}_\mathfrak{S}(S)))) \simeq \mathscr{Z}_3(\mathfrak{S}^{3mp}) \simeq \mathscr{Z}_3(\mathfrak{S}).\] The first equivalence follows from the observation that taking the centralizer of $\mathfrak{S}^{3mp} \boxtimes \mathbf{Mod}^{\mathbb{E}_3}_\mathfrak{S}(S)$ in $\mathscr{Z}_2(\mathbf{Mod}_\mathfrak{S}(S))$ is the same as first take the centralizer of $\mathbf{Mod}^{\mathbb{E}_3}_\mathfrak{S}(S)$ then take the centralizer of $\mathfrak{S}^{3mp}$. The second equivalence is the embedding of $\mathfrak{S}^{3mp}$ into its double centralizer \[\mathfrak{S}^{3mp} \simeq \mathscr{Z}_3(\mathbf{Mod}^{\mathbb{E}_3}_\mathfrak{S}(S) \to \mathscr{Z}_2(\mathbf{Mod}_\mathfrak{S}(S)))\] which we prove in Appendix \ref{sec:DoubleCentralizerTheorem}.
\end{proof}

\begin{Corollary} \label{cor:E3LocalModulesInNondegenerateSyllepticFusion2Cateory}
    Suppose $\mathfrak{S}$ is non-degenerate and $S$ is connected, then we have an equivalence of sylleptic fusion 2-categories $$ \mathscr{Z}_2(\mathbf{Mod}_\mathfrak{S}(S)) \simeq \mathfrak{S}^{3mp} \boxtimes \mathbf{Mod}_{\mathfrak{S}}^{\mathbb{E}_3}(S). $$
\end{Corollary}
    
\begin{Corollary} \label{cor:FactorizationOfNondegenerateSyllepticFusion2Cateory}
    For any non-degenerate sylleptic fusion 2-category $\mathfrak{S}$, we always have $\mathscr{Z}_2(\mathfrak{S}) \simeq \mathfrak{S}^{3mp} \boxtimes \mathfrak{S}.$
\end{Corollary}

\begin{Corollary} \label{cor:E2CenterOfModuleCategoryIsE3LocalModuleInE2Center}
    Any separable symmetric algebra $S$ in non-degenerate sylleptic fusion 2-category $\mathfrak{S}$ can be viewed as a separable symmetric algebra in $\mathscr{Z}_2(\mathfrak{S})$, with \[\mathbf{Mod}^{\mathbb{E}_3}_{\mathscr{Z}_2(\mathfrak{S})}(S) \simeq \mathscr{Z}_2(\mathbf{Mod}_\mathfrak{S}(S)).\]
\end{Corollary}

\appendix

\section{Proof of Theorem \ref{thm:E3LocalModulesAsE3Centralizer}}\label{sec:proofE3LocalModulesAsE3Centralizer}

\begin{figure}[!hbt]
    \begin{minipage}{.5\textwidth}
    \centering %131.25
    \includegraphics[width=40mm]{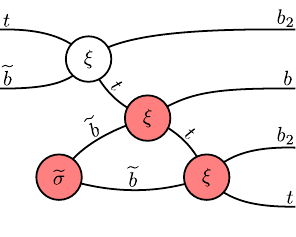}
    \caption{Step 1}
    \label{fig:AppendixA1}
    \end{minipage}%
    \begin{minipage}{.5\textwidth}
  \centering
    \centering %131.25
    \includegraphics[width=40mm]{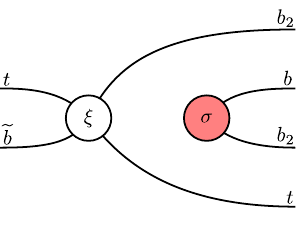}
    \caption{Step 2}
    \label{fig:AppendixA2}
    \end{minipage}
\end{figure}

\begin{figure}[!hbt]
    \begin{minipage}{.5\textwidth}
    \centering %131.25
    \includegraphics[width=56mm]{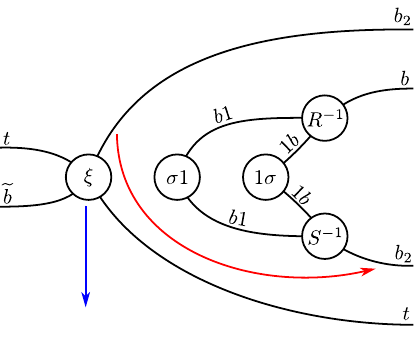}
    \caption{Step 3}
    \label{fig:AppendixA3}
    \end{minipage}%
    \begin{minipage}{.5\textwidth}
  \centering
    \centering %131.25
    \includegraphics[width=56mm]{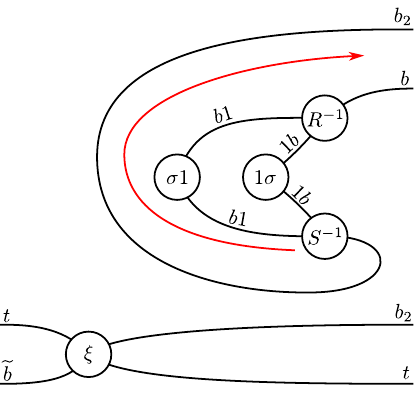}
    \caption{Step 4}
    \label{fig:AppendixA4}
    \end{minipage}
\end{figure}

\begin{figure}[hbtp]
    \begin{minipage}{.5\textwidth}
    \centering %131.25
    \includegraphics[width=56mm]{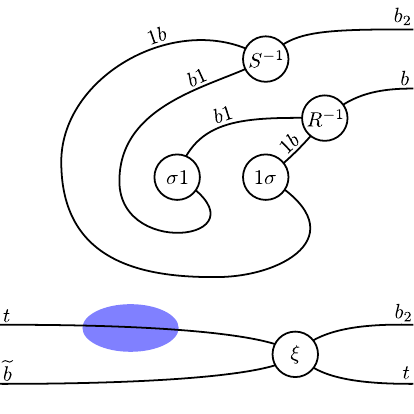}
    \caption{Step 5}
    \label{fig:AppendixA5}
    \end{minipage}%
    \begin{minipage}{.5\textwidth}
  \centering
    \centering %131.25
    \includegraphics[width=56mm]{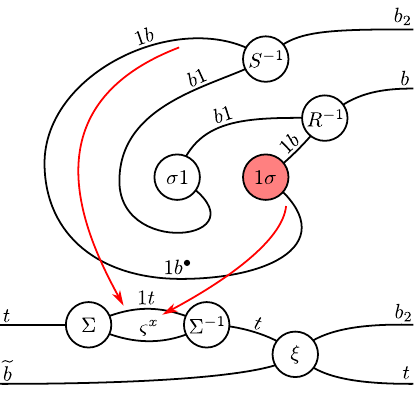}
    \caption{Step 6}
    \label{fig:AppendixA6}
    \end{minipage}
\end{figure}

\begin{figure}[hbtp]
    \centering %131.25
    \includegraphics[width=70mm]{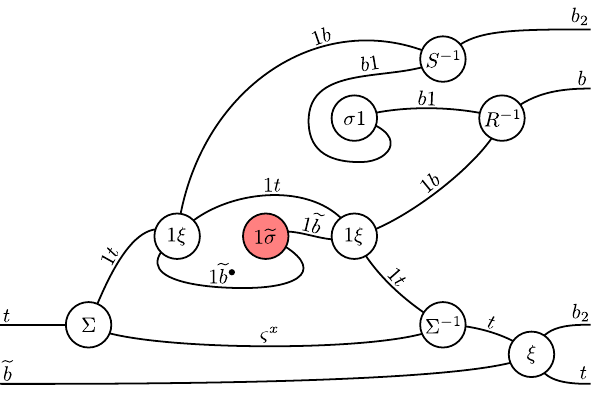}
    \caption{Step 7}
    \label{fig:AppendixA7}
\end{figure}

\begin{figure}[hbtp]
    \centering %131.25
    \includegraphics[width=70mm]{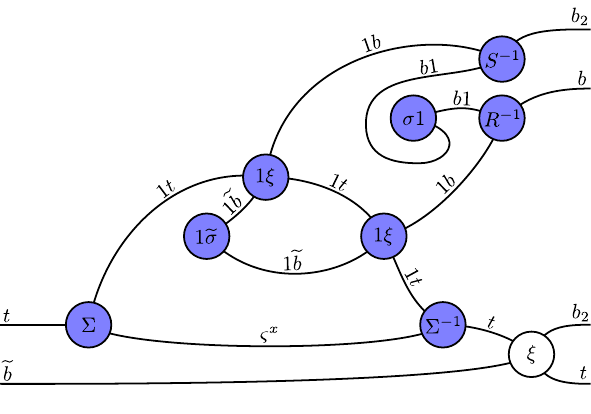}
    \caption{Step 8}
    \label{fig:AppendixA8}
\end{figure}

\begin{figure}[hbtp]
    \centering %131.25
    \includegraphics[width=40mm]{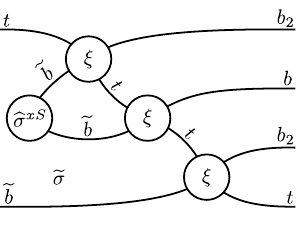}
    \caption{Step 9}
    \label{fig:AppendixA9}
\end{figure}

\newpage

\section{Double Centralizer Theorem} \label{sec:DoubleCentralizerTheorem}

In this section, we assume results about condensations between higher fusion categories \cite{GJF,JF,KZ20,KZ21}. This is laid upon the yet-to-be-established rigorous foundations of Cauchy completion of higher categories.\footnote{For a recent progress, see \cite{RZ25}.}

\begin{AppProposition}[{\cite[Theorem 3.43]{KZ20}}] \label{prop:HigherCenterAsCondensation}
    Suppose $\mathfrak{A}$ is an $\mathbb{E}_m$ monoidal $n$-category. The $n$-categorical $\mathbb{E}_m$ center\footnote{This notion is defined in \cite[Definition 3.39]{KZ20} following the universal property outlined in \cite{Lur17}.} $\mathscr{Z}_m(\mathfrak{A})$ is equivalent to $\Omega^m \mathbf{End}(\Sigma^m \mathfrak{A})$ as $\mathbb{E}_{m+1}$ monoidal $n$-categories.
\end{AppProposition}

More generally, one has the following proposition.

\begin{AppProposition}[{\cite[Theorem 4.8]{KZ21}}] \label{prop:HigherCentralizerAsCondensation}
    Suppose $\mathfrak{A}, \mathfrak{B}$ are $\mathbb{E}_m$ monoidal $n$-categories, $F:\mathfrak{A} \to \mathfrak{B}$ is an $\mathbb{E}_m$ monoidal $n$-functor. Then the $n$-categorical $\mathbb{E}_m$ centralizer\footnote{This notion is defined in \cite[Definition 4.1]{KZ21} following the universal property outlined in \cite{Lur17}.} $\mathscr{Z}_m(\mathfrak{A} \xrightarrow{F} \mathfrak{B})$ is equivalent to $\Omega^m( \mathbf{Fun}(\Sigma^m \mathfrak{A},\Sigma^m \mathfrak{B}), \Sigma^m F)$ as $\mathbb{E}_{m}$ monoidal $n$-categories. Here $(n+m)$-category $\mathbf{Fun}(\Sigma^m \mathfrak{A},\Sigma^m \mathfrak{B})$ is equipped with a distinguished object $\Sigma^m F$, hence its looping is well-defined.
\end{AppProposition}

Corollary \ref{cor:DrinfeldCenterAsBimoduleFunctors} has the following generalization.

\begin{AppProposition}[{\cite[Example 4.9]{KZ21}}] \label{prop:HigherDrinfeldCenterAsBimoduleFunctors}
    Suppose $\mathfrak{A}, \mathfrak{B}, F:\mathfrak{A} \to \mathfrak{B}$ is the same as above. Let $m=1$, then the $n$-categorical Drinfeld centralizer can be characterized as \[\mathscr{Z}_1(\mathfrak{A} \xrightarrow{F} \mathfrak{B}) \simeq \mathbf{Fun}_{\mathfrak{A} \boxtimes \mathfrak{A}^{1mp}}(\mathfrak{A},{}_{\langle F \rangle} \mathfrak{B}_{\langle F \rangle}) \simeq \mathbf{Fun}_{\mathfrak{A} \boxtimes \mathfrak{B}^{1mp}}({}_{\langle F \rangle} \mathfrak{B},{}_{\langle F \rangle} \mathfrak{B}). \]
\end{AppProposition}

Double Centralizer Theorem holds for general monoidal $n$-categories.

\begin{AppLemma}[{\cite[Corollary 4.23]{KZ21}}] \label{lem:DoubleDrinfeldCentralizer}
    Suppose $\mathfrak{A}, \mathfrak{B}, F:\mathfrak{A} \to \mathfrak{B}$ is the same as above and $m=1$. Then $\mathfrak{A}$ is canonically equivalent to its double Drinfeld centralizer: $\mathfrak{A} \simeq \mathscr{Z}_1(\mathscr{Z}_1(\mathfrak{A};\mathfrak{B});\mathfrak{B})$.
\end{AppLemma}

\begin{proof}
    By Proposition \ref{prop:HigherDrinfeldCenterAsBimoduleFunctors}, one needs to show that \[\mathfrak{A} \simeq \mathbf{Fun}_{\mathscr{Z}_1(\mathfrak{A};\mathfrak{B}) \boxtimes \mathfrak{B}^{1mp}}({}_{\langle G \rangle} \mathfrak{B},{}_{\langle G \rangle} \mathfrak{B}),\] where \[G:\mathscr{Z}_1(\mathfrak{A};\mathfrak{B}) \simeq \mathbf{Fun}_{\mathfrak{A} \boxtimes \mathfrak{B}^{1mp}}({}_{\langle F \rangle} \mathfrak{B},{}_{\langle F \rangle} \mathfrak{B}) \to \mathbf{Fun}_{\mathfrak{B}^{1mp}}(\mathfrak{B}, \mathfrak{B}) \simeq \mathfrak{B}\] is the monoidal $n$-functor induced via forgetting the left $\mathfrak{A}$-actions. Similarly, one obtains a monoidal $n$-functor $H:\mathscr{Z}_1(\mathscr{Z}_1(\mathfrak{A};\mathfrak{B});\mathfrak{B}) \to \mathfrak{B}$.
    
    By the definition of Morita equivalence, ${}_{\langle F \rangle} \mathfrak{B}_{\langle G \rangle}$ provides an invertible bimodule between $\mathfrak{A} \boxtimes \mathfrak{B}^{1mp}$ and $\mathscr{Z}_1(\mathfrak{A};\mathfrak{B})^{1mp}$. Equivalently, we can view $\mathfrak{B}$ as an invertible bimodule between $\mathfrak{A}$ and $\mathscr{Z}_1(\mathfrak{A};\mathfrak{B})^{1mp} \boxtimes \mathfrak{B}$. 
    
    By the same argument, ${}_{\langle G \rangle} \mathfrak{B}_{\langle H \rangle}$ provides an invertible bimodule between $\mathscr{Z}_1(\mathscr{Z}_1(\mathfrak{A};\mathfrak{B});\mathfrak{B})^{1mp}$ and $\mathscr{Z}_1(\mathfrak{A};\mathfrak{B}) \boxtimes \mathfrak{B}^{1mp}$. Equivalently, we can view $\mathfrak{B}$ as an invertible bimodule between $\mathscr{Z}_1(\mathscr{Z}_1(\mathfrak{A};\mathfrak{B});\mathfrak{B})$ and $\mathscr{Z}_1(\mathfrak{A};\mathfrak{B})^{1mp} \boxtimes \mathfrak{B}$. The $\mathscr{Z}_1(\mathfrak{A};\mathfrak{B})^{1mp} \boxtimes \mathfrak{B}$-actions on two copies of $\mathfrak{B}$ clearly agree with each other, so by the uniqueness of Morita dual, one gets $\mathfrak{A} \simeq \mathscr{Z}_1(\mathscr{Z}_1(\mathfrak{A};\mathfrak{B});\mathfrak{B})$.
\end{proof}

\begin{AppTheorem}
    Suppose $\mathfrak{A}, \mathfrak{B}, F:\mathfrak{A} \to \mathfrak{B}$ is the same as above. Then $\mathfrak{A}$ is canonically equivalent to its double $\mathbb{E}_m$ centralizer: $\mathfrak{A} \simeq \mathscr{Z}_m(\mathscr{Z}_m(\mathfrak{A};\mathfrak{B});\mathfrak{B})$.
\end{AppTheorem}

\begin{proof}
    One has \[\mathfrak{A} \simeq \Omega^{m-1} \Sigma^{m-1} \mathfrak{A} \simeq \Omega^{m-1} \mathscr{Z}_1(\mathscr{Z}_1(\Sigma^{m-1} \mathfrak{A};\Sigma^{m-1} \mathfrak{B});\Sigma^{m-1} \mathfrak{B}) \] 
    \[ \simeq \mathscr{Z}_m(\mathscr{Z}_m(\mathfrak{A};\mathfrak{B});\mathfrak{B}), \] where the first equivalence follows from the definition of Karoubi completion, the second equivalence comes from Lemma \ref{lem:DoubleDrinfeldCentralizer}, the last equivalence follows from Proposition \ref{prop:HigherCentralizerAsCondensation}.
\end{proof}

As corollaries, for $n=2$ and $m=1$, one gets the Double Centralizer Theorem for Drinfeld centralizer of fusion 2-categories. For $n=2$ and $m=2$, this provides a proof of Double Centralizer Theorem for braided centralizer of braided fusion 2-categories, which we used in the proof of Theorem \ref{thm:E2CenterAndE2LocalModules}. We see that this can be obtained from the Double Centralizer Theorem for Drinfeld centralizer of fusion 3-categories.

For $n=2$ and $m=3$, this provides the Double Centralizer Theorem for sylleptic centralizer of sylleptic fusion 2-categories, which we used in the proof of Theorem \ref{thm:E3CenterAndE3LocalModules}. Again, this can be obtained from the Double Centralizer Theorem for Drinfeld centralizer of fusion 4-categories, or the Double Centralizer Theorem for braided centralizer of braided fusion 3-categories.

As a bonus, for $n=1$ and $m=2$, the above theorem provides the Double Centralizer Theorem for Müger centralizer of braided fusion 1-categories, which is originally proved in \cite[Theorem 3.10]{DGNO} via Frobenius-Perron dimensions. In retrospect, we achieved many results for tensor categories using dimensional arguments, and yet a sensible categorification of dimension for higher fusion categories is still missing in the literature.

\begin{AppQuestion}
    Can one categorify the notion of Frobenius-Perron dimension for higher fusion categories?
\end{AppQuestion}

\bibliographystyle{alpha}
\bibliography{bibliography}

\end{document}